\documentclass{xbook}
\usepackage{amsthm,amsfonts,amssymb}

\theoremstyle{definition}

\theoremstyle{remark}

\numberwithin{section}{chapter}
\numberwithin{equation}{chapter}
\numberwithin{figure}{chapter}
\numberwithin{table}{chapter}

\begin{document}
\hoffset=-2truecm

\frontmatter
\title{The Based Ring of Two-Sided Cells of\\
Affine Weyl Groups of Type $\tilde A_{n-1}$}

\author{Nanhua Xi}
\address{Institute of Mathematics, Chinese Academy of Sciences, Beijing 100080,
 China}
\email{nanhua@math08.math.ac.cn}
\thanks{The   author was supported in part by Chinese National Sciences
Foundation.}


\date{June 25, 1999}

\begin{abstract}
In this paper we prove Lusztig's conjecture on based ring for an
affine Weyl group of type $\tilde A_{n-1}$.
\end{abstract}

\def\ga{\gamma}
\def\tt{\tilde T}
\def\medksip{\medskip}
\def\meskip{\medskip}
\newcommand{\A}{{\mathcal A}}
\maketitle

\tableofcontents

\setcounter{page}{6}

\chapter*{Introduction}

\def\bbC{\mathbb{C}}

The Kazhdan-Lusztig theory deeply increases our understanding of
Coxeter groups and their representations and their role in Lie
representation theory. The central concepts in Kazhdan-Lusztig theory
include Kazhdan-Lusztig polynomial, Kazhdan-Lusztig basis, cell, based
ring. Kazhdan-Lusztig polynomials play an essential role in
understanding certain remarkable representations in Lie theory, for
instances, the representations of quantum groups at roots of 1, the
modular representations of algebraic groups, the representations of
Kac-Moody algebras. Kazhdan-Lusztig basis and cells are very useful
 in
understanding structure and representations of Coxeter groups and
their Hecke algebras.

\medskip

The based ring of a two-sided cell of certain Coxeter groups is
defined through Kazhdan-Lusztig basis by Lusztig in [L5]. The based
ring is very interesting in understanding the concerned Hecke algebras
and their representations (see [L5, L6, L8]). Moreover we can
construct irreducible representations of the Hecke algebras if we know
the structure of the based ring explicitly (see [X3]). This fact is
remarkable since constructing computable irreducible modules of affine
Hecke algebras is in general difficult. Recently Lusztig introduced
periodic $W$-graphs for constructing finite dimensional modules of
affine Hecke algebras (see [L9]).

For the structure of the based ring of a two-sided cell in a Weyl
group or in an affine Weyl group, Lusztig has two nice
conjectures, one for Weyl groups (see [L7]) and the other for
affine Weyl groups that says the based ring is isomorphic to a
certain equivariant $K$-group (see [L8]). The conjecture for Weyl
groups is proved by Lusztig (unpublished). The conjecture for
affine Weyl groups is proved for rank 2 cases, for the lowest
two-sided cell of an affine Weyl group and for the second highest
two-sided cell of an affine Weyl group, see [X1, X3]. As an
application, the author gives a classification of irreducible
modules of the Hecke algebras of affine Weyl groups of rank 2 for
any non-zero parameter, see [X3]. (When the parameter is not a
root of 1, Kazhdan and Lusztig worked out the classification of
irreducible modules of an affine Hecke algebra, see [KL2].) Also
the author computed the irreducible modules of affine Hecke
algebras associated with  second highest two-sided  cell, see
[X3]. Recently Ram has developed an interesting combinatorical
approach to study representations of affine Hecke algebras, see
[R1, R2].

\medskip

In this paper we will prove Lusztig Conjecture on based ring for
an extended affine Weyl group associated with  the general linear
group $GL_n(\bbC)$ or the special linear group $SL_n(\bbC)$, which
is of type $\tilde A_{n-1}$.

\medskip

Let us here briefly explain the idea of the proof of the
conjecture for type $\tilde A_{n-1}$. For each two-sided cell of
the extended affine Weyl group associated with  $GL_n(\bbC)$, we
first show that the based ring of the two-sided cell is a matrix
algebra  over the based ring of the intersection of a left cell in
the two-sided cell and its inverse (a right cell). Then we show
that the based ring of the intersection is isomorphic to the
(rational) representation ring of a certain connected reductive
algebraic group over complex numbers. The main difficulty is to
establish a bijection between the intersection (of a left cell in
the two-sided cell and its inverse) and the set of isomorphism
classes of irreducible rational representations of the algebraic
group, and to show that the bijection leads to the isomorphism
between the based ring of the intersection and the representation
ring of the algebraic group. The bijection is established in
Chapter 5. In Chapters 6 and 7 we prove several formulas in the
based ring of the intersection, the corresponding formulas in the
representation ring of the concerned algebraic group are obvious.
Using the formulas established in Chapters 6 and 7, in Chapter 8
we show that the bijection in Chapter 5 is the right one. In
Chapter 8 we also show that Lusztig Conjecture on based ring is
true for the extended affine Weyl group associated with
$SL_n(\bbC)$ and can not be generalized to arbitrary extended
affine Weyl groups.

\medskip

The contents of the paper are as follows.

\medskip

In Chapter 1 we recall some basic definitions and facts about cells
and based rings. In section 1.4 we prove some simple properties about
the structure constants of a based ring that are important to our
proof of Lusztig Conjecture on based ring for an affine Weyl group of
type
$
\tilde A_{n-1} $. In Chapter 2 we discuss the structure of cells
in the extended affine Weyl group associated with  $GL_n(\bbC)$.
The cells in the extended affine Weyl group are explicitly
described by Shi and Lusztig (see [S, L3]). In section 2.3, we
show that the based ring of a two-sided cell of the extended
affine Weyl group is isomorphic to a matrix algebra  over the
based ring of the intersection of any given left cell in the
two-sided cell and its inverse (a right cell). Thus to understand
the based ring of the two-sided cell, we only need to work out the
structure of the based ring of the intersection of a left cell in
the two-sided cell and its inverse. In section 2.4 we discuss
chains and antichains defined by Shi, which are essential for our
bijection between the intersection of a left cell and its inverse
and the set of isomorphism classes of irreducible rational
representations of the corresponding algebraic group.

In Chapter 3 we give some discussion to canonical left cells.
Although we do not need the results in this chapter for the proof
of our main result of the paper, but for other types canonical
left cells maybe play a big role for some questions such as
Lusztig Conjecture on based ring and properties of Lusztig
bijection between the set of two-sided cells of an affine Weyl
group and the set of unipotent classes of the corresponding
algebraic group.

\medskip

In Chapter 4 we describe the reductive algebraic group
corresponding to a two-sided cell in the extended affine Weyl
group associated with $GL_n(\bbC)$ and recall some needed results
about the representation ring of a general linear group over
complex numbers.

In Chapter 5 we establish a bijection between the intersection of a
left cell and its inverse and the set of dominant weights of the
corresponding reductive algebraic group. We obtain this bijection by
means of antichains and this bijection is one key to our main result.

In Chapters 6 and 7 we work with the based ring of the
intersection of a left cell and its inverse. Believing in that
Lusztig Conjecture is true, motivated by some simple
multiplication formulas in the representation ring of the
concerned reductive algebraic group, we prove some multiplication
formulas in the based ring. Using the formulas in Chapters 6 and
7, in Chapter 8 we prove that the based ring of the intersection
is isomorphic to the representation ring of the corresponding
reductive algebraic group. This completes our proof of Lusztig
Conjecture on based ring for type $\tilde A_{n-1}$. In Chapter 8
we also discuss the based rings of extended affine Weyl groups
associated with  the projective linear group $PGL_n(\bbC)$ and
with the special linear group $SL_n(\bbC)$.

\medskip

It is expected that the explicit knowledge on the based rings will
have nice application to representation theory of the concerned affine
Hecke algebras. It is interesting that we can use the explicit
structure on based ring to find leading coefficients of some
Kazhdan-Lusztig polynomials, the details will appear elsewhere.

\medskip

I am very grateful to Professors V. Chari and J.-Y. Shi for
helpful discussions. I thank Professor Le Yang for his constant
encouragement.  I would like to thank Professor G. Lusztig for
many helpful comments. I am greatly indebted to the referee for
carefully reading, very valuable comments and for pointing out a
serious gap in an earlier version of the paper.  Part of the work
was done during my visit to Department of Mathematics at
University of California at Riverside, I thank the Department of
Mathematics for hospitality. The author was in part supported by
Chinese National Sciences Foundation.

\mainmatter

\chapter{Cells in Affine Weyl Groups}

%

\def\A{\mathcal A}

In this chapter we first recall some basic concepts such as cell,
based ring, etc., then in section 1.4 we give some discussions to star
operations. In section 1.1 we recall the definition of Kazhdan-Lusztig
basis of a Hecke algebra. In section 1.2 we recall the definition of
cell and of the $a$-function. In section 1.3 we recall some properties
about the integers $\gamma_{w,u,v}$ (defined through the structure
constants of Kazhdan-Lusztig basis), which are due to Lusztig.

The star operation was introduced in [KL1] and is a useful tool to
study cells of Coxeter groups. In section 1.4 we prove some
interesting results about the relations between star operations and
the integers $\gamma_{w,u,v}$. The main result in this section is
Theorem 1.4.5 which says that the integers $\gamma_{w,u,v}$ are
invariant under star operations. In
  section 1.5  we   recall the definition of based ring (due to
Lusztig, see [L5]) and Lusztig's nice conjecture about the
structure of the based ring of a two-sided cell of an affine Weyl
group. A few comments about the conjecture are given . In section
1.6 we discuss the relationship between the generalized star
operations (see [L4]) and the integers $\gamma_{w,u,v}$. The
results in section 1.6 should be helpful for understanding the
structure of the based ring of a two-sided cell of an affine Weyl
group of type $\tilde B_n,\ \tilde C_n,\ \tilde F_4$.

The basic references for this chapter are [KL1] and [L4-L8].

\smallskip

\goodbreak
\section{Hecke algebra}

\def\om{\omega}
\def\Om{\Omega}

Here Hecke algebras are defined for extended Coxeter groups.

\medskip

Let $(W',S)$ be a Coxeter system with $S$ the set of simple
reflections. Assume that a commutative group $\Omega$ acts on
$(W',S)$. Then we can consider the extended Coxeter group
$W=\Omega\ltimes W'$. The {\bf length function} $l$ on $W'$ and the
{\bf partial order} $\leq$ on $W'$ are extended to $W$ as usual, that
is, $l(\omega w)=l(w)$, and $\omega w\leq
\omega' u$ if and only if $\omega=\om'$ and $w\leq u$, where
$\om,\om'$ are in $\Om$ and $w,u$ are in $W'$.

\medskip

\def\bbZ{\mathbb{Z}}
Let $q$ be an indeterminate. The {\bf Hecke algebra} $H$ of $(W,S)$
over $\A=\bbZ[q,q^{-1}]$ with parameter $q^2$ is an associative
algebra over $\A$, with an $\A$-basis $\{T_w \ |\ w\in W\}$ and
relations (1) $(T_s-q^2)(T_s+1)=0$ if $s\in S$, (2) $T_wT_u=T_{wu}$ if
$l(wu)=l(w)+l(u)$.

Let $a\to\bar a$ be the involution of $\A$ defined by $\bar
q=q^{-1}$. Then we have a bar involution of $H$ defined by
$$\overline{\sum a_wT_w}=\sum \bar a_wT^{-1}_{w^{-1}},\qquad
a_w\in\A.$$ For each $w\in W$ there is a unique element $C_w$ in
$H$ such that $\bar C_w=C_w$ and $C_w=q^{-l(w)}\sum_{y\leq
w}P_{y,w}(q^2)T_y,$ where $P_{y,w}$ is a polynomial in $q$ of
degree $\leq \frac 12(l(w)-l(y)-1)$ if $l(w)>l(y)$ and
$P_{w,w}=1$.

\medskip

The basis $\{C_w\ |\ w\in W\}$ is the famous {\bf Kazhdan-Lusztig
basis} of the Hecke algebra $H$ and the polynomials $P_{y,w}$ are
the celebrated {\bf Kazhdan-Lusztig polynomials}. The basis and
the polynomials have deep relations with the geometry of Schubert
varieties when $W$ is the Coxeter group associated with  a
Kac-Moody algebra and the polynomials are essential in
understanding some difficult irreducible representations (such as
the finite dimensional irreducible representations of quantum
groups at roots of 1, some irreducible representations of
Kac-Moody algebras, irreducible rational representations of
algebraic groups over an algebraic closed field of positive
characteristic, etc.).

The basis also plays an important role in understanding the
representations of the Hecke algebra $H$, see [KL1, L4-L8]. Through
the basis Lusztig defined based ring for some Coxeter groups
(including Weyl groups and affine Weyl groups). In this paper we will
determine the structure of the based ring of an affine Weyl group of
type $\tilde A_{n-1}$. Applications to the representation theory of
the concerned Hecke algebras will appear elsewhere.

\medskip

If $y\leq w$ and $y\ne w$, we have $P_{y,w}=\mu(y,w)q^{\frac
12(l(w)-l(y)-1)}+$lower degree terms. We write $y\prec w$ if
$\mu(y,w)\ne 0$. We then set $\mu(y,w)=\mu(w,y)$. We shall write $y-w$
if $\mu(y,w)\ne 0$ or $\mu(w,y)\ne 0$. We have

\bigskip

\noindent (a) Assume that  $y\leq w$ and $s\in S$. If $sw\leq w$,
then $sy\le w$ and $P_{y,w}=P_{sy,w}$. If $ws\leq w$, then $ys\le
w$ and $P_{y,w}=P_{ys,w}$. See [KL1].

\def\ll{\underset L\leq}
\def\rl{\underset R\leq}
\def\lrl{\underset {LR}\leq}
\def\el{\underset L\sim}
\def\er{\underset R\sim}
\def\elr{\underset {LR}\sim}

\medskip

\section{Cell and $a$-function}

Cells of Coxeter groups are defined in [KL1] and are useful in
understanding the structure and representations of $W$ and of the
Hecke algebra $H$. Also we can systematically construct
representations of Hecke algebras by means of the cells (see [KL1, L4,
L6, L8]).

\medskip

Let us recall the definition of cells of Coxeter groups. For $w\in W$
set $$L(w)=\{s\in S\ |\ sw\leq w\}$$ $$R(w)=\{s\in S\ |\ ws\leq w\}.$$
Let $w$ and $u$ be elements in $W'$. We say that $w\ll u$ (resp. $w\rl
u;w\lrl u)$ if there exists a sequence $w=w_0,w_1,w_2,...,w_k=u$ in
$W'$ such that for $i=1,2,...,k$ we have $\mu(w_{i-1},w_i)\ne 0$ and
$L(w_{i-1})\not\subseteq L(w_i)$ (resp. $R(w_{i-1})\not\subseteq
R(w_i)$; $L(w_{i-1})\not\subseteq L(w_i)$ or $R(w_{i-1})\not\subseteq
R(w_i)$). Then for any $\om,\om'$ in $\Om$ we say that $\om w\ll\om'
u$ (resp. $w\om \rl u\om' $; $\om w\lrl\om' u$) if $ w\ll u$ (resp.
$w\rl u,\ w\lrl u$).

For elements $w,u$ in $W$ we write that $w\el u$ (resp. $w\er u$;
$w\elr u$) if $w\ll u\ll w$ (resp. $w\rl u\rl w$; $w\lrl u\lrl
w$). The relations $\ll,\rl,\lrl$ are preorders on $W$ and the
relations $\el,\er,\elr$ are equivalence relations on $W$. The
corresponding equivalence classes are called {\bf left cells,
right cells, two-sided cells} of $W$, respectively. The preorder
$\ll$ (resp. $\rl;\lrl$) induces a partial order on the set of
left (resp. right; two-sided) cells of $W$, denoted again by $\ll$
(resp. $\rl;\lrl$). The following properties (see [KL1] or [L4,
L5]) will be needed.

\bigskip
\def\no{\noindent}

\no(a) If $w\ll u$, then $R(w)\supseteq R(u)$. In particular, if $w\el
u$, then $R(w)=R(u)$.

\smallskip

\no (b) If $w\rl u$, then $L(w)\supseteq L(u)$. In particular, if $w\er
u$, then $L(w)=L(u)$.

\bigskip

For any $\xi,\xi'\in H$ and $w\in W$ we write $$\xi C_w=\sum_{v\in
W}h_{v}C_v,\qquad h_{v}\in\A ,$$
$$ C_w\xi'=\sum_{v\in
W}h'_{v}C_v,\qquad h'_{v}\in\A ,$$
$$\xi C_w\xi'=\sum_{v\in
W}h''_{v}C_v,\qquad h''_{v}\in\A .$$ Then

\bigskip

\noindent (c) $v\ll w$ (resp. $v\rl w$; $v\lrl w$) if $h_v\ne 0$
(resp. $h'_v\ne 0$; $h''_v\ne 0$).

\smallskip

\no(d)  Assume that $(W',S)$ is crystallographic and that there
exists a positive integer $a_0$ such that $a(v)\le a_0$ for all
$v\in W$ (see below for the definition of $a(v)$). If $h_v\ne 0$
(resp. $h'_v\ne 0$) and $w\not\el v$ (resp. $w\not\er v$), then
$w\not\elr v$.

\bigskip

The $a$-{\bf function} on $W$ introduced by Lusztig is a useful tool
to study cells of $W$ and is also necessary for defining based ring of
two-sided cells. Given $w,u$ in $W$ we write
$$C_wC_u=\sum_{v\in
W}h_{w,u,v}C_v,\qquad h_{w,u,v}\in\A .$$ For $v\in W$ we define
$a(v)$=the minimal non-negative integer $i$ such that
$q^{-i}h_{w,u,v}$ is in $\bbZ[q^{-1}]$ for all $w,u$ in $W$. If such
$i$ does not exist we set $a(v)=\infty$.

\medskip

Set $\tilde T_w=q^{-l(w)}T_w$ for $w\in W$. Then $\tilde T_w\in
C_w+q^{-1}\sum_{x\in W}\bbZ[q^{-1}] C_x$. Write
$$\tilde T_w\tilde T_u =\sum_{v\in
W}h'_{w,u,v}C_v,\qquad h'_{w,u,v}\in\A .$$ Then we also have
$a(v)$=the minimal non-negative integer $i$ such that
$q^{-i}h'_{w,u,v}$ is in $\bbZ[q^{-1}]$ for all $w,u$ in $W$. See
[L4].

\medskip

\section{Affine Weyl group}

\def\D{\mathcal D}

\def\bbC{\mathbb{C}}
We are mainly interested in extended affine Weyl groups and their
Hecke algebras (especially of type $\tilde A_{n-1})$. In this section
we recall some properties about the integers $\gamma_{w,u,v}$ defined
through the Kazhdan-Lusztig basis of an affine Hecke algebra. The
integers are actually structure constants of the based ring of the
extended affine Weyl group.

\medskip

Let $G$ be a connected reductive group over $\bbC$. Let $W_0, X,
P, R$ be the corresponding Weyl group, weight lattice, root
lattice, root system, respectively. Then $W'=W_0\ltimes P$ is an
{\bf affine Weyl group} and $W=W_0\ltimes X$ is an {\bf extended
affine Weyl group}. Let $S$ be a set of simple reflections of $W'$
such that $S\cap W_0$ generates $W_0$ and is a set of simple
reflections of $W_0$. We can find an abelian subgroup $\Om$ of $W$
such that $\om S=S\om$ for any $\om\in\Om$ and $W=\Om\ltimes W'$.

Lusztig proved that the number of left cells in $W$ is finite.
Moreover each left cell of $W$ contains a unique element of
$$\D=\{w\in
W'\ |\ 2\text{deg}P_{e,w}=l(w)-a(w)\},$$ where $e$ is the neutral
element of $W$. The set $\D$ is a finite set of involutions of $W'$.
The elements in $\D$ will be called {\bf distinguished involutions}.
See [L5].

\medskip

\def\ga{\gamma}
\def\Ga{\Gamma}
Let $w_0$ be the longest element of $W_0$. We have $a(w)\leq
l(w_0)=a(w_0)$ for all $w$ in $W$ (see [L4]). Thus for any
$w,u,v\in W$ we can define an integer $\gamma_{w,u,v}$ by the
condition $$q^{-a(v)}h_{w,u,v}-\ga_{w,u,v}\in
q^{-1}\bbZ[q^{-1}],$$ see 1.2 for the definition of $h_{w,u,v}$.
The following are some properties of $\ga_{w,u,v}$ (see [L5] for
(a)-(e) and [L4] for (f)).

\bigskip

\no (a) If $\gamma_{w,u,v}$ is not equal to 0, then $w\el u^{-1},\ u\el v$ and $w\er
v$. In particular we have $w\elr u\elr v$ if $\ga_{w,u,v}$ is not
equal to 0.

\medskip

\no(b) $\ga_{w,u,v}=\ga_{u,v^{-1},w^{-1}}$ and
$\ga_{u^{-1},w^{-1},v^{-1}}=\ga_{w,u,v}$.

\medskip

\no(c) Let $d$ be in $\D$. Then $\ga_{w,d,u}\ne 0$ if and only if $w=u$ and
$w\el d$. Moreover
$\ga_{w,d,w}=\ga_{d,w^{-1},w^{-1}}=\ga_{w,w^{-1},d}=1$.

\medskip

\no(d) $w\el u^{-1}$ if and only if $\ga_{w,u,v}$ is not equal to 0 for some $v$.

\medskip

\no(e) If $w\lrl u$ then $a(w)\ge a(u)$. In particular if $w\elr u$
then $a(w)=a(u)$. If $w\ll u$ (resp. $w\rl u$; $w\lrl u$) and
$a(w)=a(u)$, then $w\el u$ (resp. $w\er u$; $w\elr u$).

\medskip

\no(f) The positivity: $\ga_{w,u,v}\ge 0$ for all $w,u,v$ in $W$.

\medskip

\section{Star operation}

Let $(W,S)$ be as in section 1.3. There is a very useful operation on
$W$, called star operation, introduced by Kazhdan and Lusztig in
[KL1]. In this section we study the relations between star operation
and the integers $\gamma_{w,u,v}. $ The main result is Theorem 1.4.5.

\bigskip

\no{\bf 1.4.1.} Let $s$ and $t$ be in $S$   such that $st$
 has order 3, i.e. $sts=tst$. Define

$$D_L(s,t)=\{w\in W\  |\   L(w)\cap\{s,t\}\text{ has exactly one
element}\},$$

$$D_R(s,t)=\{w\in W\ |\ R(w)\cap\{s,t\}\text{ has
exactly one element}\}.$$

\medskip

If $w$ is in $D_L(s,t)$, then $\{sw,tw\}$ contains exactly one element
in $D_L(s,t)$, denoted by ${}^*w$, here $*=\{s,t\}$. The map:
$D_L(s,t)\to D_L(s,t)$, $w\to {}^*w$, is an involution and is called a
{\bf left star operation}. Similarly if $w\in D_R(s,t)$ we can define
the {\bf right star operation} $w\to w^*=\{ws,wt\}\cap D_R(s,t)$ on
$D_R(s,t)$, where $*=\{s,t\}$. The following are some properties
proved in [KL1].

\bigskip

 Let $*=\{s,t\}$. Denote by $<s,t>$ the subgroup of $W$ generated by
$s$ and $t$. Assume that $y,w$ are in $D_L(s,t)$. We have

\medskip

\no(a) If $yw^{-1}$ is not in $<s,t>$, then $y\prec w$ if and only if
${}^*y\prec{}^*w$. Moreover $\mu(y,w)=\mu({}^*y,{}^*w)$.

\medskip

\no(b) If $yw^{-1}$ is in $<s,t>$, then $y\prec w$ if and only if
${}^*w\prec{}^*y$. Moreover $\mu(y,w)=\mu({}^*w,{}^*y)=1$.

\medskip

\no(c) $y\er w$ if and only if ${}^*y\er{}^*w$.

\medskip

\no (d) $w\el {}^*w$.

\bigskip

 Let $*=\{s,t\}$. Assume that  $y,w$ are in $D_R(s,t)$. We have

\medskip

\no(e) If $y^{-1}w$ is not in $<s,t>$, then $y\prec w$ if and only if
${y}^*\prec{w}^*$. Moreover $\mu(y,w)=\mu({}y^*,{w}^*)$.

\medskip

\no(f) If $y^{-1}w$ is in $<s,t>$, then $y\prec w$ if and only if
${w}^*\prec{y}^*$. Moreover $\mu(y,w)=\mu({w}^*,{y}^*)=1$.

\medskip

\no(g) $y\el w$ if and only ${w}^*\el {y}^*$.

\medskip

\no(h) $w\er w^*$.

\medskip

Finally we have

\medskip

\no(i) If $\Ga\subseteq D_R(s,t)$ is a left cell of $W$, then
 $\Ga^*=\{w^*\ |\ w\in\Ga
\}$ is a left cell of $W$, here $*=\{s,t\}$.

\bigskip

The following lemma shows that $C_{w^*}$ and $C_{{}^*w}$ have nice
relationship with $C_w$.

\bigskip

\def\st{\stackrel}
\def\sc{\scriptstyle}
\noindent{\bf Lemma 1.4.2.} { \sl Let $s,t$ be in $ S$   such that
$st $ has order 3. Set $*=\{s,t\}$.}

\smallskip

\no (a) {\sl Assume that  $w$ is in $D_L(s,t)$ and $s\ {}^*w\leq {}^*w$. Then
$$C_sC_w=C_{ {}^*w}+\displaystyle\sum_{\st{\sc y-w}{\st{\sc sy<y}
{ty<y}}}
\mu(y,w)C_y.$$}

\no (b){ \sl Assume that  $w$ is in $D_R(s,t)$ and $w^*\ s\leq {w}^*$. Then
$$C_wC_s=C_{ {w}^*}+\displaystyle\sum_{\st{\sc y-w}{\st{\sc ys<y}
{yt<y}}}\mu(y,w)C_y.$$}

\bigskip

{\it Proof.} (a) Since $s\ {}^*w\leq {}^*w$, we have $sw\geq w$ and
$tw\leq w$. Thus (see [KL1])
$$C_sC_w=C_{sw}+\sum_{\st {\sc y\prec w}{ sy<y}}\mu{(y,w)}C_y.$$
If $w=tw'$ with $sw'\geq w'$ and $tw'\geq w'$, then ${}^*w=sw$. By 1.1
(a), for $y\prec w$ with $sy<y$ we have $ty<y$. In this case (a) is
true. If $w=tsw'$ with $sw'\geq w'$ and $tw'\geq w'$, then ${}^*w=sw'$
and $sw=stsw'$. Moreover, by 1.1 (a), if $y\ne {}^*w$, $y\prec w$ and
$sy<y$ then we have $ty<y$. In this case (a) is also true. (a) is
proved. The proof of (b) is similar.

\bigskip

\noindent{\bf Lemma 1.4.3.} {\sl Let $s,t,s',t'$ be in $ S$  such
that both $st$ and $s't'$ have order 3. Assume that $w$ is in $
D_L(s,t)\cap D_R(s',t')$. Set $*=\{s,t\}$ and $ \star=\{s',t'\}$.
Then}

\smallskip

\no(a) {\sl ${}^*w$ is in $ D_R(s',t')$ and $w^\star$ is in $ D_L(s,t)$.}

\no(b) {\sl ${}^*(w^\star)=(^*w)^\star$. We shall write
${}^*w^\star$ for ${}^*(w^\star)=(^*w)^\star$.}

\bigskip

\def\itp{{\it Proof.}\ }
\itp (a) By 1.4.1 (d),  $w\el {}^*w$, thus we have $R(w)=R(^*w)$ (see 1.3 (a)),
 so ${}^*w$ is in $ D_R(s',t')$.
Similarly we see that $w^\star$ is in $ D_L(s,t)$.

\medskip\smallskip

(b) By (a), $(^*w)^\star$ and ${}^*(w^\star)$ are well defined. Since
  $(^*w)^\star\er {}^*w$ and ${}^*(w^\star)\el w^\star$, we have

\medskip

\noindent(1) $L(^*w)=L((^*w)^\star)$ and $R(w^\star)=R(^*(w^\star))$.

\medskip

Thus we have

\medskip

  \noindent (2)  Both  $(^*w)^\star$ and ${}^*(w^\star)$ are in $D_L(s,t)\cap D_R(s',t')$.

\medskip

   It is no
harm to assume that $s\ {}^*w\leq {}^*w$ and $w^\star\ s'
\leq {w^\star}$. Then $sw\ge w$ and $ws'\ge w$. We also have
$sw^\star\ge w^\star$ and ${}^*ws'\ge{}^*w$ since $w^\star\er w$ and
${}^*w\el w$. Write
$$(C_sC_w)C_{s'}=\sum_{y\in W}h_yC_y,\qquad h_y\in\A.$$
By Lemma 1.4.2 and 1.2 (b-c), we get

\medskip

\noindent (3) $h_{(^*w)^\star}\ne 0$. Moreover, if
 $y\ne (^*w)^\star$ and  $h_y\ne 0$,
then $y$ is not in $D_L(s,t)$ or $y$ is not in $ D_R(s',t')$.

\medskip

Noting that $C_s(C_wC_{s'})=(C_sC_w)C_{s'}$, using Lemma 1.4.2 again
we see $h_{^*(w^\star)}$ is not equal to 0. Now using (2) and (3) we
get ${}^*(w^\star)=(^*w)^\star$. (b) is proved.

\medskip

The lemma is proved.

\bigskip

\no{\bf Proposition 1.4.4.}
{\sl Let $s,t,s',t'$ be as in Lemma 1.4.3. Set $*=\{s,t\}$ and
$\star=\{s',t'\}$. Suppose that $w$ is in $ D_L(s,t)$ and $u$ is in $
D_R(s',t')$. Let $v$ be in $W$ such that $v\el u$ and $v\er w$. Then}

\smallskip

\def\ds{\displaystyle\sum}
\no(a) {\sl We have $v\in D_L(s,t)\cap
D_R(s',t')$, so we can define ${}^*v^\star$.}

\no(b) {\sl $h_{w,u,v}=h_{{}^*w,u^\star, {}^*v^\star}$, see 1.2 for the
definition of $h_{w,u,v}$.}

\bigskip

\itp (a) is trivial.

\medskip

(b) We have ${}^*v^\star\el v^\star$ and ${}^*v^\star\er {}^*v$. By
definition we have

\medskip

\no(1) $C_wC_u=\ds_{z\in W}h_{w,u,z}C_z,\ \ h_{w,u,z}\in\A$.

\medskip

We may assume that $sw\ge w$ and $us'\ge u$. Using Lemma 1.4.2 we get

\medskip

\no(2)  $C_sC_w=C_{ {}^*w}+\displaystyle\sum_{\st{\sc y\ne{}^*w}
{\st{\sc sy<y}
{ty<y}}} h_yC_y,\ \ h_y\in\A.$

\medskip

\no(3) $C_uC_{s'}=C_{ {u}^\star}+\displaystyle\sum_{\st{\sc x\ne u^\star}
{\st{\sc xs'<y} {xt'<y}}}h'_xC_x,\ \ h'_x\in\A.$

\medskip

Write

\medskip

\no(4) $C_sC_zC_{s'}=\ds_{z'\in W}f_{z,z'}C_{z'},\ \ f_{z,z'}\in\A.$

\medskip

Then we have

\medskip

\no(5) $C_sC_wC_uC_{s'}=\ds_{z'\in W}(h_{{}^*w,u^\star,z'}+
\ds_{\st{\sc y\ne{}^*w}{\st{\sc sy<y}
{ty<y}}} \displaystyle h_y h_{y,u^\star,z'}$

$\qquad\qquad +
 \displaystyle\sum_{\st{\sc x\ne u^\star} {\st{\sc xs'<y}
{xt'<y}}} h'_xh_{{}^*w,x,z'}+
\ds_{\st{\sc y\ne{}^*w}{\st{\sc sy<y}
{ty<y}}} \displaystyle\sum_{\st{\sc x\ne u^\star} {\st{\sc xs'<y}
{xt'<y}}}h_yh'_xh_{y,x,z'})C_{z'}$.

\medskip

\no(6) $C_s\ds_{z\in W}h_{w,u,z}C_zC_{s'}=\ds_{z'\in W}\ds_{z\in
W}h_{w,u,z}f_{z,z'}C_{z'}$.

\medksip

Using (1) and comparing (5) with (6), we get

\medskip

\no (7) $\qquad h_{{}^*w,u^\star,z'}+\ds_{\st{\sc y\ne{}^*w}{\st{\sc sy<y}
{ty<y}}} \displaystyle h_y h_{y,u^\star,z'}+
 \displaystyle\sum_{\st{\sc x\ne u^\star} {\st{\sc xs'<y}
{xt'<y}}} h'_xh_{{}^*w,x,z'}+\ds_{\st{\sc y\ne{}^*w}{\st{\sc sy<y}
{ty<y}}} \displaystyle\sum_{\st{\sc x\ne u^\star} {\st{\sc xs'<y}
{xt'<y}}}h_yh'_xh_{y,x,z'}$

$=\ds_{z\in W}h_{w,u,z}f_{z,z'}$.

\meskip

Using 1.2 (a-b) we see

\medskip

\no(8) Assume that  $z'\er {}^*w$ and $z'\el u^\star$. Then  $$\ds_{\st{\sc y\ne{}^*w}{\st{\sc sy<y}
{ty<y}}} \displaystyle h_y h_{y,u^\star,z'}+
 \displaystyle\sum_{\st{\sc x\ne u^\star} {\st{\sc xs'<y}
{xt'<y}}} h'_xh_{{}^*w,x,z'}+\ds_{\st{\sc y\ne{}^*w}{\st{\sc sy<y}
{ty<y}}}
\displaystyle\sum_{\st{\sc x\ne u^\star} {\st{\sc xs'<y}
{xt'<y}}}h_yh'_xh_{y,x,z'}=0.$$

\medskip

When $h_{w,u,v}=0$ we must have $h_{{}^*w,u^\star,{}^*v^\star}=0$.
Otherwise, noting that ${}^*v^\star\er{}^*w$ and ${}^*v^\star\el
u^\star$, by (7) and (8) we have $h_{{}^*w,u^\star,{}^*v^\star}
=\ds_{z\in W}h_{w,u,z}f_{z,{}^*v^\star}$ and
$h_{w,u,z}f_{z,{}^*v^\star}\ne 0$ for some $z$. Assume that
$h_{w,u,z}f_{z,{}^*v^\star}\ne 0$. By 1.2 (c) and 1.3 (e) we have
$z\er w$ and $z\el u$. By the proof of Lemma 1.4.3 we then have
$f_{z,{}^*z^\star}=1$, and $z''\not\er{}^*z$,\ \ $z''\not\el
z^\star$ if $f_{z,z''}\ne 0$ and $z''\ne {}^*z^\star$. Thus
${}^*v^\star={}^*z^\star$ and $v=z$. This contradicts that
$h_{w,u,v}= 0$. Therefore $h_{{}^*w,u^\star,{}^*v^\star}=0$
whenever $h_{w,u,v}=0$.

Now suppose that $h_{w,u,v}\ne 0$. As above we have $z=v$ whenever
$h_{w,u,z}f_{z,{}^*v^\star}\ne 0,$ and $f_{v,{}^*v^\star}=1$.
  Using (7) and (8) we then get $h_{{}^*w,u^\star,{}^*v^\star}=
  h_{w,u,v}$ in this case.
We proved (b).

\medskip

The Proposition is proved.

\bigskip

\noindent{\bf Theorem 1.4.5.} {\sl Suppose that $w$ is in
$D_L(s,t)\cap D_R(s'',t'')$ and $u$ is in $D_L(s'',t'')\cap
D_R(s',t')$. Set $*=\{s,t\},\#=\{s'',t''\},$ and
$\star=\{s',t'\}$. Let $v$ be in $W$ such that $v$ is in
$D_L(s,t)\cap D_R(s',t')$. Then we have

\smallskip

\centerline{$\ga_{w,u,v}=\ga_{{}^*w^\#,{}^\#u^\star,{}^*v^\star}.$}}

\bigskip

\itp By Lemma 1.4.3 and Proposition 1.4.4 (a), the elements
${}^*w^\#,{}^\#u^\star$ and ${}^*v^\star$ are well defined. If
either $v\not\er w$ or $v\not\el u$, then both sides of the wanted
equality are 0. Now suppose that  $v \er w$ and $v\el u$.
According to Proposition 1.4.4 (b) we have
$h_{w,u,v}=h_{{}^*w,u^\star, {}^*v^\star}$. Therefore
$\ga_{w,u,v}=\ga_{{}^*w,{}u^\star,{}^*v^\star}.$ It is easy to see
that $(^*w)^{-1}=(w^{-1})^*$, \ ${}^\#(w^{-1})^*=(^*w^\#)^{-1}$,
and ${}^\#(({}^*w)^{-1}) =(^*w^\#)^{-1}$. Thus we have
$$\begin{array}{rl}  \ga_{w,u,v}&=
\ga_{{}^*w,{}u^\star,{}^*v^\star}\qquad\qquad\quad \text{using 1.3
(b)}\\[3mm] &=\ga_{{}u^\star,({}^*v^\star)^{-1},(^*w)^{-1}}\qquad
\text {using Prop. 1.4.4}
\\[3mm]
&=\ga_{{}^\#u^\star,({}^*v^\star)^{-1},{}^\#((^*w)^{-1})}  \\[3mm]
&=\ga_{{}^\#u^\star,({}^*v^\star)^{-1},(^*w^\#)^{-1}}
\qquad \text{using 1.3 (b)}\\[3mm]
&=\ga_{{}^*w^\#,{}^\#u^\star,{}^*v^\star} .
\end{array}$$
The theorem is proved.

\bigskip

\noindent{\bf Proposition 1.4.6.} {\sl Let $W$ be an extended
affine Weyl group.}

\smallskip

\no(a) {\sl Let $I$ be a subset of $S$ such that the subgroup $W_I$ of $W$
generated by $I$ is finite. Then the longest element $w_I$ is a
distinguished involution.

In (b) and (c) $d$ is a distinguished involution.}

\no(b) {\sl For any $\om\in
\Om$, the element $\om d\om^{-1}$ is a distinguished involution.}

\no (c) {\sl Suppose $s,t\in S$ and $st $ has order 3. Then $d\in D_L(s,t)$ if and
only if $d\in D_R(s,t)$. If $d\in D_L(s,t)$, then ${}^*d^* $ is a
distinguished involution.}

\bigskip

{\it Proof.} (a) is well known, see for example [L5].

\medskip

(b) Since $C_\om C_wC_{\om^{-1}}=C_{\om w\om^{-1}}$ and
$C_{\om}C_{\om^{-1}}=1$, for any $w,u,v$ in $W$ we have $$P_{\om
u\om^{-1},\om w\om^{-1}}=P_{u,w}$$ and $$h_{\om w\om^{-1},\om
u\om^{-1},\om v\om^{-1}}=h_{w,u,v}.$$ In particular we have $a(\om
w\om^{-1})=a(w)$ and $P_{e,\om w\om^{-1}}=P_{e,w}$ for any $w\in
W$. Noting that $l(\om w\om^{-1})=l(w)$ for any $w$ in $W$, we see
$$\text{deg}P_{e,\om d\om^{-1}}=\text{deg}P_{e,d}=\frac
12(l(d)-a(d))=\frac12(l(\om d\om^{-1})-a(\om d\om^{-1})).$$ By
definition, $\om d\om^{-1}$ is a distinguished involution.

\medskip

(c) Let $d'$ be the distinguished involution of the left cell
containing ${}^*d^*$. We have
$$({}^*d^*)^{-1}={}^*(d^{-1})^*={}^*d^*,$$ so ${}^*d^*$ is an
involution. Using 1.3 (c-d) we get $\ga_{{}^*d^*,d',{}^*d^*}=1$. Using
Theorem 1.4.5 we then have $\ga_{ d,{}^*{d'}^*,d}=1$. Using 1.3 (b) we
get $\ga_{{}^*{d'}^*,d,d}=1$. Applying 1.3 (c) we have $d={}^*{d'}^*$.
Therefore $d'={}^*d^*$.

\medskip

The proposition is proved.

\medskip

\section{Based ring}

\def\jg{J_{\Ga\cap\Ga^{-1}}}
\def\gg{\Ga\cap\Ga^{-1}}
Following [L5] we define the based rings. Let $J$ be the free
$\bbZ$-module with a basis $\{t_w\ |\ w\in W\}$. The multiplication
$t_wt_u=\sum_{v\in W}\ga_{w,u,v}t_v$ defines an associative ring
structure on $J$, see [L5]. The ring $J$ is called the {\bf based
 ring}
of $W$, its unit is $\sum_{d\in\D}t_d$. According to 1.3(a), for
each two-sided  cell $\bold c$ of $W$, the $\bbZ$-submodule
$J_{\bold c}$ of $J$, spanned by all $t_w$ $(w\in\bold c)$, is a
two-sided ideal of $J$. The ideal $J_{\bold c}$ is in fact an
associative ring with unit $\sum_{d\in\D\cap \bold c}t_d$. The
ring $J_{\bold c}$ is called the {\bf based ring of the two-sided
cell} $\bold c$. For a left cell $\Ga$ of $W$, the
$\bbZ$-submodule $\jg$ of $J$, spanned by all $t_w$ $( w\in
\Ga\cap \Ga^{-1}$), is also an associative ring, its unit is
$t_d$, here $d$ is the unique distinguished involution in $\Ga$.

\bigskip

\noindent{\bf Proposition 1.5.1.}  {\sl Let $\Ga$ be a left cell of
$W$.}

\smallskip

\no(a) {\sl For any $\om\in\Om$, $\Ga'=\om\Ga\om^{-1}$ is a left cell of $W$.
Moreover the map $w\to \om w\om^{-1}$ induces an isomorphism between
the based rings $\jg$ and $J_{\Ga'\cap{\Ga'}^{-1}}$.}

\no(b) {\sl  Let $s,t\in S$ be such that $st$ has order 3. Suppose
$\Ga\subset D_R(s,t)$. Then $\jg\simeq J_{  \Ga^*\cap
(\Ga^*)^{-1}}$, here $*=\{s,t\}$.}

\bigskip

\itp (a) Obviously $\Ga'$ is a left cell of $W$. From the proof of  Prop.
1.4.6 (b) we see $\ga_{\om w\om^{-1},\om u\om^{-1},\om
v\om^{-1}}=\ga_{w,u,v}$ for any $w,u,v\in W$. Therefore the map $w\to
\om w\om^{-1}$ induces an isomorphism between the based rings $\jg$
and $J_{\Ga'\cap{\Ga'}^{-1}}$.

\medskip

(b) Suppose $w\in\gg$. Then ${}^*w^*$ is in $\Ga^*\cap
(\Ga^*)^{-1}$. The map $w\to {}^*w^*$ defines a bijection between
$\gg$ and $\Ga^*\cap (\Ga^*)^{-1}$. According to Theorem 1.4.5 we
see that the map $t_w\to t_{{}^*w^*}$ defines a ring isomorphism
between $\jg$ and $J_{\Ga^*\cap (\Ga^*)^{-1}}$.

\medskip

The proposition is proved.

\bigskip

The based rings $J_{\bold c}$ are very interesting in
understanding the structure and representations of Hecke algebras
of $W$, see [L5-L8]. If we know the structure of $J_{\bold c}$
explicitly we can construct modules of affine Hecke algebras in a
computable way, see [L8, X3].

Lusztig has a nice conjecture concerning the structure of $J_{\bold
c}$. Assume that $G$ is connected and has a simply connected derived
group. There is a natural bijection between the set of two-sided cells
of $W$ and the set of unipotent classes of $G$, see [L8]. Assume that
$u$ is an element in the unipotent class corresponding to a two-sided
cell {\bf c} of $W$. Denote by $F_{\bold c}$ a maximal reductive
subgroup of the centralizer $C_G(u)$ of $u$ in $G$. Lusztig
conjectured that there exists a finite $F_{\bold c}$-set $Y$ and a
bijection $\pi:
\bold c\to$ the set of isomorphism classes of irreducible $F_{\bold
c}$ vector bundles on $Y\times Y$ such that $t_w\to \pi(w)$ defines a
ring isomorphism (preserving the unit element) between $J_{\bold c}$
and $K_{F_{\bold c}}(Y\times Y)$ and $\pi(w^{-1})=\widetilde
{\pi(w)}$, see [L8] for the conjecture and for the definition of
$\widetilde {\pi(w)}$.

\medskip

\def\bc{\bold c}
When $F_{\bold c}$ is connected (or equivalently $C_G(u)$ is
connected), $F_{\bold c}$ must act on $Y$ trivially. In this case
$|Y|$ is the number of left cells contained in $\bold c$ and
$K_{F_{\bold c}}(Y\times Y)$ is isomorphic to the $|Y|\times|Y|$
matrix algebra  $M_{|Y| }(R_{F_\bc})$ over the rational
representation ring $R_{F_\bc}$ of $F_{\bold c}$. Let Irr$F_{\bold
c}$ be the set of isomorphism classes of rational irreducible
representations of $F_\bc$. Then Lusztig Conjecture says that
there is a bijection $$\pi:\bold c\to
 \{(V,i,j)\ |\ V\in \text{Irr}F_{\bold c}, \ 1\le i,j\le
|Y|\}$$ such that (1) the map $t_w\to \pi(w)$ defines a ring
isomorphism between $J_\bc$ and $M_{|Y|}(R_{F_\bc})$, where we
identify $\pi(w)=(V,i,j)$ with the matrix whose $(i,j)$-entry is $V$
and other entries are 0, (2) $\pi(w^{-1})=(V^*,j,i)$ if
$\pi(w)=(V,i,j)$, here $V^*$ is the dual of $V$.

Lusztig Conjecture has been proved for the following cases, (1)
$\bc$ is the lowest two-sided cell of $W$, (that is, $\bold c$ is
the two-sided cell of $W$ containing the longest element $w_0$ of
$W_0$,) see [X1], (2) rank 2 cases, see [X3], (3) the case
$a(\bc)=1$, see [X3].

\medskip

When $G=GL_n(\bbC)$, each $F_{\bold c}$ is connected. The purpose
of this paper is to show that Lusztig Conjecture is true for the
extended affine Weyl group associated with  $GL_n(\bbC)$.

\medskip

\section{Star operation, II}

The star operation, introduced in [KL1], was generalized in [L4]. We
are interested in the relationship between the (generalized) star
operation and the structure constants $\gamma_{w,u,v}$ of the based
ring of an extended affine Weyl group. In this section we show that
Theorem 1.4.5 remains true for the generalized star operation. We also
give a few other identities about the constants $\gamma_{w,u,v}$. The
results in this section are not used in the sequent chapters but
should be useful for understanding the based ring of an arbitrary
extended affine Weyl group.

Let $(W,S)$ be as in section 1.3. Assume that $s$ and $t$ are in
$S$ and $st$ has order $m$. Denote by $U$ the subgroup of $W$
generated by $s$ and $t$. Each coset $Uw$ can be decomposed into
four parts: one consists of the unique element $x$ of minimal
length in the coset, one consists of the unique element $y$ of
maximal length in the coset, one consists of the $m-1$ elements
$sx,\ tsx,\ stsx,\ ...,$ one consists of the $m-1$ elements $tx,\
stx,\ tstx,\ ...$. The last two subsets are called left strings
(related to $\{s,t\}$) and shall be regarded as sequences (as
above) rather than subsets. Similarly we define right strings
(related to $\{s,t\}$). We have (see [L4])

\medskip

\noindent (a) A left string in $W$ is contained in a left cell of $W$
and a right string in $W$ is contained in a right cell of $W$.

\medskip

Assume that $w$ is in a left (resp. right) string (related to
$\{s,t\}$) of length $m-1$ and is the $i$th element of the left (resp.
right) string, we define ${}^*w$ (resp. $w^*$) to be the $(m-i)$th
element of the string, where $*=\{s,t\}$. We have

\bigskip
\def\no{\noindent}

\noindent{\bf Lemma 1.6.1.} {\sl Let $w$ be in $W$ such that $w$ is
in a left string related to $*=\{s,t\}$ and is also in a right string
related to $\star=\{s',t'\}$. Then}

\smallskip

\no(a) {\sl ${}^*w $  is in a right string related to $\{s',t'\}$
and $w^\star$ in a left string related to $\{s,t\}$.}

\no(b) {\sl ${}^*(w^\star)=(^*w)^\star$. We shall write ${}^*w^\star$ for
${}^*(w^\star)=(^*w)^\star$.}

\bigskip

The proof is similar to that of Lemma 1.4.3 although more complicated.

\bigskip
\def\ga{\gamma}

The following theorem is a generalization of Theorem 1.4.5

\bigskip

\noindent{\bf Theorem 1.6.2.} {\sl Let $w,u,v$ be in $W$ such that (1) $w$
is in a left string related to $*=\{s,t\}$ and also in a right string
related to $\#=\{s',t'\}$, (2) $u$ is in a left string related to
$\#=\{s',t'\}$ and also in a right string related to
$\star=\{s'',t''\}$, (3) $v$ is in a left string related to
$*=\{s,t\}$ and also in a right string related to $\star=\{s'',t''\}$.
Then}

\smallskip

\centerline{$\ga_{w,u,v}=\ga_{{}^*w^\#,{}^\#u^\star,{}^*v^\star}.$}

\bigskip

The proof is similar to that of Theorem 1.4.5.

\bigskip

\noindent{\bf 1.6.3.} In this subsection we assume that $s$ and $t$
are in $S$ and $st$ has order 4. Let $w,u,v$ be in $W$ such that
$l(ststw)=4+l(w)$ and $l(ststv)=4+l(v)$. As in the proof of Lemma
1.4.2 we have
$$C_tC_{sw}=C_{ tsw}+\displaystyle\sum_{\st{\sc y-w}{\st{\sc ty<y}
{sy<y}}}
\mu(y,w)C_y,$$
$$C_tC_{stw}=C_{ tstw}+C_{tw}+\displaystyle\sum_{\st{\sc y-w}{\st{\sc ty<y}
{sy<y}}}
\mu(y,w)C_y,$$
$$C_tC_{stsw}=C_{ tsw}+\displaystyle\sum_{\st{\sc y-w}{\st{\sc ty<y}
{sy<y}}}
\mu(y,w)C_y.$$

Considering the products $C_tC_{sw}C_u,\ C_tC_{stw}C_u,\
C_tC_{stsw}C_u,$ and using the associativity of multiplication of the
Hecke algebra $H$, we get

\bigskip

\no(a) $h_{tsw,u,tv}=h_{sw,u,stv},$

\medskip

\no(b) $h_{tsw,u,tsv}=h_{sw,u,sv}+h_{sw,u,stsv},$

\medskip

\no(c) $h_{tsw,u,tstv}=h_{sw,u,stv},$

\medskip

\no(d) $h_{tstw,u,tv}+h_{tw,u,tv}=h_{stw,u,stv},$

\medskip

\no(e) $h_{tstw,u,tsv}=h_{stw,u,stsv},$

\medskip

\no(f) $h_{tstw,u,tstv}+h_{tw,u,tstv}=h_{stw,u,stv}.$

\medskip

In particular, we have

\medskip

\no(a') $\gamma_{tsw,u,tv}=\gamma_{sw,u,stv},$

\medskip

\no(b') $\gamma_{tsw,u,tsv}=\gamma_{sw,u,sv}+\gamma_{sw,u,stsv},$

\medskip

\no(c') $\gamma_{tsw,u,tstv}=\gamma_{sw,u,stv},$

\medskip

\no(d') $\gamma_{tstw,u,tv}+\gamma_{tw,u,tv}=\gamma_{stw,u,stv},$

\medskip

\no(e') $\gamma_{tstw,u,tsv}=\gamma_{stw,u,stsv},$

\medskip

\no(f') $\gamma_{tstw,u,tstv}+\gamma_{tw,u,tstv}=\gamma_{stw,u,stv}.$

\medskip

One may compare the above equalities with the equalities in [L4,
(10.4.2)]. Using 1.3 (b) and the above equalities we can get more
equalities for the structure constants of the based ring $J$, we omit
them. When $st$ has higher order, there exist similar equalities. We
omit the discussion.

\chapter{Type $\tilde A_{n-1}$}


\def\ta{\tilde A_{n-1}}
\def\gl{GL_n(\bbC)}
From now on we will concentrate on type $\tilde A_{n-1}$. The
cells of an affine Weyl group of type $\ta$ have been described in
[S, L3]. In this chapter we first recall some facts about the
cells, then we derive some new results for our purpose. In section
2.1 we recall an alternative definition (due to Lusztig) for the
extended affine Weyl group $W$ associated with  $\gl$. In section
2.2 we recall the description of cells of $W$ in [S, L3]. For this
description and later use we slightly refine the definition of
chain and antichain in [S]. More precisely we will define d-chain
(resp. d-antichain) and r-chain (resp. r-antichain). In this
section we also consider the intersection of left cells and right
cells. The intersections are important to our purpose. In section
2.3 we show that the based ring of a two-sided cell of $W$ is
isomorphic to a matrix algebra  over the based ring
$J_{\Ga\cap\Ga^{-1}}$ for any left cell $\Ga$ in the two-sided
cell (see Theorem 2.3.2). This is the main result of this chapter.
Thus to understand Lusztig's conjecture on the structure of the
based ring of the two-sided cell we only need to understand the
structure of $J_{\Ga\cap\Ga^{-1}}$. We will do this in Chapters
5-8.

For later use in section 2.4 we work out some results about chains
and antichains. The results will be needed in Chapter 5 for
defining the required bijection between $\Ga\cap\Ga^{-1}$ and
Irr$F_\bc$, see Chapter 5. In section 2.5 we prove a result about
star operation.

\medskip

\section{The affine Weyl group associated  with $GL_n(\bbC)$}

In this section we recall a definition of Lusztig for affine Weyl
group of type $\tilde A_{n-1}$, see [L2].

\medskip

Suppose that $G=GL_n(\bbC)$ is the general linear group over $\bbC$ of
degree $n$. Then the Weyl group $W_0$ of $G$ is isomorphic to the
symmetric group of $n$ letters. Let $T$ be the subgroup of $G$
consisting of all diagonal matrices in $G$. Then $T$ is a maximal
torus of $G$ and the weight lattice $X$ is the character group of $T$.
Let $\tau_i\in X$ be the homomorphism $T\to
\bbC^*,\
\text{diag} (a_1,...,a_n)\to a_i$. We have $X=<\tau_1,...,\tau_n>$ and
 the extended affine Weyl group associated with  $GL_n(\bbC)$ is
 $W=W_0\ltimes X$. We
identify $W_0$ with the group of all the $n\times n$ permutation
matrices (by a permutation matrix, we mean a  monomial matrix
whose nonzero entries are all equal to 1). Let $s_i$
$(i=1,2,...,n-1)$ be the simple reflection of $W_0$ obtained from
the identity matrix by interchanging  the $i$th and the $(i+1)$th
rows. Then $s_i\tau_i=\tau_{i+1}s_i$ and $s_i\tau_j=\tau_{j}s_i$
if $j\ne i,i+1$.

\medskip

Another realization of $W$, due to Lusztig, is as follows.
Consider the permutations $\sigma: \bbZ\to\bbZ$ that satisfy
$\sigma(i+n)=\sigma(i)+n$ and $ \sum_{i=1}^{n}(\sigma(i)-i)\equiv
0\text{(mod }n).$ All such permutations form a permutation group
$W_*$ of $\bbZ$.

\bigskip

\noindent{\bf Lemma 2.1.1:} {\sl $W$ is isomorphic to $W_*$.}

\bigskip

{\it Proof.} We can define a $W$-action on $\bbZ$ as follows, \[
s_i(j)
= \left\lbrace
           \begin{array}{c l}
              j & \text{if $j\not\equiv i,i+1\text{(mod }n)$},\\
              j+1 & \text{if $j\equiv i\text{(mod }n)$},\\
              j-1 & \text{if  $j\equiv i+1\text{(mod }n);$}
           \end{array}
         \right. \]

\[
\tau_i(j)
= \left\lbrace
           \begin{array}{c l}
              j & \text{if $j\not\equiv i\text{(mod }n)$},\\
              j+n & \text{if $j\equiv i\text{(mod }n)$}.
           \end{array}
         \right. \]

It is easy to check that this action defines an isomorphism between
$W$ and $W_*$. The lemma is proved.
\bigskip

We shall identify $W$ and $W_*$ using the above isomorphism.

\bigskip

\noindent{\bf Lemma 2.1.2.} {\sl The permutation
$\om:\bbZ\to\bbZ, i\to i+1$ for all $i$, is in $W$.}

\bigskip

{\it Proof.} It is clear from the definition.

\bigskip

\no{\bf 2.1.3.} The following are some simple properties of the
extended affine Weyl group $W$.

\medskip

\no(a) The center of $W$ is generated by $\om^n$.

\medskip

\no(b) $\om s_i=s_{i+1}\om$, where $s_0$ is defined by
$s_0(0)=1,s_0(1)=0,s_0(i)=i$ if $i\not\equiv 0,1\text{(mod }n)$,
and $s_i=s_j$ if $i\equiv j\text{(mod }n)$.

\medskip

\no(c) Let $\Om$ be the subgroup of $W$ generated by $\om$ and
$W'$ the subgroup of $W$ generated by all $s_i$. Then $\Om$ is an
infinite cyclic group and $W'$ is an affine Weyl group of type
$\tilde A_{n-1}$. Moreover we have $W\simeq \Om\ltimes W'$.

\medskip

\no(d) We have $\tau_i=\om s_{i-2}s_{i-3}\cdots s_0 s_{n-1}\cdots
s_i$ if $i\geq 2$ and $\tau_1=\om s_{n-1}s_{n-2}$ $\cdots s_1.$
Note that the length of $\tau_i$ is $n-1$.

\medskip

\no(e) As a special case of the length formula in [IM], we have
$$l(w)=\sum_{1\leq i<j\leq
n}\left|\left[\frac{w(j)-w(i)}n\right]\right|,$$ where [h] is the
integer part of h (recall that $0\le h-[h]<1)$. For a direct proof
of the formula, see [S, Lemma 4.2.2].

\medskip

\noindent(f) Let $w\in W$. Then  $w(k)<w(k+1)$ if and
only if $w\le ws_k$, see [S, Corollary 4.2.3].

\medskip

\def\l{\lambda}
\section{Cells}

The cells in $W$ have a beautiful combinatorial description,
conjectured by Lusztig and proved by Shi. In this section we recall
the description of cells of $W$ in [S, L3].

\medskip

Following Shi we define chains and antichains. More precisely we
refine slightly his definition by defining d-chains, d-antichains,
r-chains and r-antichains.

\medskip

Let $w$ be an element of $W$ and  $j_1,j_2,\cdots,j_k$ be
integers. We call $j_1,j_2,\cdots,j_k$ a {\bf d-chain} of $w$ of
length $k$ and $w(j_1),...,w(j_k)$ an {\bf r-chain} of $w$ of
length $k$ if

\no(1) $j_1<j_2<\cdots<j_k$,

\no(2) $j_i\not\equiv j_{i'}\text{(mod }n)$ whenever $1\le i\ne
i'\le k$,

\no(3) $w(j_1)>w(j_2)>\cdots
>w(j_k)$,

\no here d-chain means domain chain and r-chain means range chain.
 Thus a d-chain of $W$ essentially is a subset of $\bbZ$ whose
natural order is reversed by $w$.

\smallskip

A {\bf d-chain (resp. r-chain) family set } of $w$ of index $q$ is
a subset $Y$ of $\bbZ$ such that (1) elements in $Y$ are
noncongruent to each other  modulo $n$, (2) $Y$ is a disjoint
union of $q$ d-chains (resp. r-chains)  $A_1,...,A_q$ of $w$. We
also call $\{A_1,...,A_q\}$ a {\bf d-chain (resp. r-chain) family}
of $W$.

\medskip

We call $j_1,j_2,\cdots,j_k$ a {\bf d-antichain} of $w$ of length
$k$  and $w(j_1),...,w(j_k)$ an  {\bf r-antichain} of $w$ of
length $k$ if

\no(1) $j_k-n<j_1<j_2<\cdots<j_k$, (then $j_i\not\equiv
j_{i'}$(mod $n$) whenever $1\le i\ne i'\le k$,)

\no(2) $w(j_k)-n<w(j_1)<w(j_2)<\cdots <w(j_k)$.

\smallskip

A {\bf d-antichain (resp. r-antichain) family set} of $w$
 of
index $q$ is a subset $Y$ of $\bbZ$ such that (1) elements in $Y$
are noncongruent to each other   modulo  $n$, (2) $Y$ is a
disjoint union of $q$ d-antichains (resp. r-antichains)
$A_1,...,A_q$ of $w$.  We also call $\{A_1,...,A_q\}$ a {\bf
d-antichain (resp. r-antichain) family} of $W$.

\medskip

Obviously if $Y$ is an r-chain (resp. r-antichain) family set of
$w$ of index $q$, then $Y$ is a d-chain (resp. d-antichain) family
set of $w^{-1}$ of index $q$.

\medskip

We shall regard $W$ as the permutation group $W_*$ of $\bbZ$.
Following Lusztig we associate a partition $\lambda$ of $n$ with
an element $w$ of $W$ as follows. Let $d_i$ be the maximal one
among the cardinalities of all d-chain family sets  of $w$ of
index $i$. Then $d_1\leq d_2\leq\cdots \leq d_n=n$. According to
[G, Th. 1.5], $d_1\geq d_2-d_1\geq d_3-d_2\geq\cdots\geq
d_n-d_{n-1}$. We call $$\l(w)=(d_1,d_2-d_1,...,d_n-d_{n-1})$$ the
{\bf partition associated with  $w$}. According to [S, L3], $w\elr
u$ if and only if $\l(w)=\l(u)$. Moreover the number of left cells
in a two-sided cell corresponding to a partition $\l$ of $n$ is
$n_\mu= {n!}/({\mu_1!\mu_2!\cdots \mu_{r'}!})$, where
$(\mu_1,...,\mu_{r'})$ is the dual partition of $\l$, see [S].

The dual partition of $\l$ can be defined through antichains of
$w$. Let $d'_i$ be the maximal one among the cardinalities of all
d-antichain family sets of $w$ of index $i$. Then $d'_1\leq
d'_2\leq\cdots \leq d'_n=n$. According to [G, Th. 1.5], $d'_1\geq
d'_2-d'_1\geq d'_3-d'_2\geq\cdots\geq d'_n-d'_{n-1}$. The
partition $$\mu(w)=(d'_1,d'_2-d'_1,...,d'_n-d'_{n-1})$$ is the
partition dual to the partition $\l(w)$, see [G, Th. 1.6].

\medskip

\def\Gal{\Ga_\l}
\def\Gali{\Ga_{\l,i}}
\def\Phili{\Phi_{\l,i}}

Let $\bold c$ be the two-sided cell of $W$ corresponding to a
partition $\l=(\l_1,...,\l_r)$ of $n$. Denote by $w_\l$ the
longest element of the subgroup of $W$ generated by all $$s_1,\
...,\ s_{\l_1-1},\ s_{\l_1+1},\ ...,\ s_{\l_1+\l_2-1},\ ...,\
s_{\l_1+\cdots+\l_{r-1}+1},\ ...,\ s_{\l_1+\cdots+\l_{r}-1}.$$
Then $w_\l\in \bc$ and there is a unique left cell in $\bold c$
containing $w_\l$, denoted by $\Ga_\l$. We shall write
$\Ga_{\l,i}$ for the left cell containing $\om^i w_\l\om^{-i}$.
Set $\Phili=\Gali^{-1}$, this is a right cell contained in $\bold
c$.

Denote by $u_\l$ the longest element of the subgroup of $W$ generated
by all
$$s_1,\ ...,\ s_{\l_r-1},\ s_{\l_r+1},\ ...,\ s_{\l_r+\l_{r-1}-1},\ ...,\
s_{\l_r+\cdots+\l_{2}+1},\ ...,\ s_{\l_r+\cdots+\l_{1}-1}.$$ Then
$u_\l\in \bc$ and there is a unique left cell in $\bold c$ containing
$u_\l$, denoted by $\Delta_\l$. We shall write $\Delta_{\l,i}$ for the
left cell containing $\om^i u_\l\om^{-i}$. Set
$\Psi_{\l,i}=\Delta_{\l,i}^{-1}$. Note that the intersection of the
set $H_\l$ in [S, 9.3] and $\bold c$ is contained in the union
$\bigcup_{1\leq i\leq n}\Delta_{\l,i}^{-1}$.
\bigskip

\noindent{\bf Lemma 2.2.1.} {\sl Assume that  $w\in \bc$. Then through a succession of left star
operations and of right star operations on $w$ we can get an element
in $\Gali\cap\Phi_{\l,j}$ for some integers $i,j$.}

\bigskip

\itp It follows from Prop. 9.3.7, Theorem 1.6.3 (i) and Lemma 18.3.2 in [S].

\bigskip

\noindent{\bf Corollary 2.2.2.} {\sl Assume that  $\Ga$ is a left cell in $\bold c$ and
$\Phi$ is a right cell in $\bold c$. Then $\Ga\cap\Phi$ can be
obtained by applying a succession of left star operations and of right
star operations on some $\Gali\cap\Phi_{\l,j}$.}

\bigskip

Since $\Gali=\Gal\om^{-i}$, the map $w\to \om^j w\om^{-i}$ defines a
bijection between $\Gal\cap\Gal^{-1}$ and $\Gali\cap\Phi_{\l,j}$. Thus
it is very fundamental to understand the properties of
$\Gal\cap\Gal^{-1}$.

\bigskip

\noindent{\bf Proposition 2.2.3.} {\sl Let $\Ga$ be a left cell in $\bold c$. Then
there is a bijection $\phi: \Gal\cap\Gal^{-1}\to \Ga\cap\Ga^{-1}$ such
that $t_w\to t_{\phi(w)}$ defines an isomorphism between
$J_{\Gal\cap\Gal^{-1}}$ and $J_{\Ga\cap\Ga^{-1}}.$ }

\bigskip

\itp Using Corollary 2.2.2 and Lemma 18.3.2 in [S], we can find $i$ such that $\om^i(\Ga\cap\Ga^{-1})\om^{-i}$ is
obtained from
$\Gal\cap\Gal^{-1}$ by applying a succession of left star operations
and the corresponding right star operations. Using Prop. 1.5.1 we see
that the assertion is true.

\medskip

\section{The based ring $J_{\bc}$}

In this section we will show that the based ring $J_\bc$ of a
two-sided cell $\mathbf c$ of $W$ is a matrix algebra  over
$J_{\Ga\cap\Ga^{-1}}$ for any left cell $\Ga$ in $\bc$. This is
the main result of this chapter and is also one of the key steps
of our proof of Lusztig Conjecture for the structure of $J_\bc$.

\medskip

Recall that for the two-sided cell $\mathbf c$ corresponding to a
partition $\l$ of $n$, we have a unique left cell $\Ga_\l$ containing
$w_\l$. We number the left cells in $\mathbf c$ as $
\Ga_\l=\Ga_1,\Ga_2,...,\Ga_{n_\mu}$,
where $n_\mu=n!/(\mu_1!\cdots\mu_{r'}!)$ and
$\mu=(\mu_1,\mu_2,...,\mu_{r'})$ is the dual partition of
$\l=(\l_1,...,\l_r)$. Let $\Phi_i=\Ga_i^{-1}$ be the right cell
corresponding to $\Ga_i$. We would like to define a bijection
$\phi_{ij}$ between $A_{ij}=\Phi_i\cap\Ga_j$ and $
A_{11}=\Phi_1\cap\Ga_1=\Ga_1^{-1}\cap\Ga_1$.

\bigskip

There are several cases. First we suppose $j=1$. Let $l\geq 1$ be
such that $\om^lw_\l \om^{-l}=w_\l$ but $\om^hw_\l \om^{-h}\ne
w_\l$ for $1\leq h\leq l-1$.

\smallskip

\no(a) If $\Phi_i=\om^h\Phi_1$ for some $0\leq h\leq l-1$, then obviously
$\om^hw\to w$ defines a bijection from $A_{i1}$ to $A_{11}$.

\smallskip

\no (b) Now let $\Phi_i$ be arbitrary. Using  [S, Prop.9.3.7], we
can find $h_i$ between $0$ and $l-1$ such that $\Phi_i$ is
obtained from $\Phi =\om^{h_i}\Phi_1$ by a sequence of left star
operations. Such $h_i$ and sequence of left star operations are
usually not unique. We fix such an $h_i$ and the sequence of left
star operations. Then $\Phi_i\cap\Ga_1$ is obtained from $\Phi
\cap\Ga_1$ by applying the sequence of left star operations. This
of course establishes a bijection between $\Phi_i\cap\Ga_1$ and
$\Phi \cap\Ga_1$. Using (a) we then get a bijection between
$\Phi_i\cap\Ga_1$ and $A_{11}=\Ga^{-1}_1\cap\Ga_1$.

\smallskip

\no(c) Now for any $\Phi_i$ we have fixed an $\om^{h_{i}}$ and a
sequence of left star operations. Then $\Ga_i$ can be obtained
from $\Ga_1\om^{-h_i}$ by applying the corresponding sequence of
right star operations. By this way and using (a) we get a
bijection $\phi_{ij}$ between $A_{ij}$ and $A_{11}$, cf. Lemma
2.2.1.

\bigskip

The following are some properties of the bijection $\phi_{ij}:
A_{ij}\to A_{11}$.

\bigskip

\no{\bf Lemma 2.3.1.} (a) {\sl Let $d_i$ be the distinguished involution
in $A_{ii}$. Then $\phi_{ii}(d_i)=w_\l$.}

\no (b) {\sl Note that $A_{ij}^{-1}=A_{ji}$. For $x\in A_{ji}$ we have
$\phi_{ij}(x^{-1})=(\phi_{ji}(x))^{-1}$.}

\bigskip

\itp Using Prop. 1.4.6 we see that (a) is true. (b) follows from
$(^*x^\star)^{-1}={}^\star(w^{-1})^*$ and $(\om x\om^{-1})^{-1}=
\om x^{-1}\om^{-1}$.

\bigskip

\def\bc{\mathbf c}

We shall use $E(t_w,i,j)$ for any square matrix whose $(i,j)$-entry is
$t_w$ and other entries are 0.

\bigskip

\noindent{\bf Theorem 2.3.2.} {\sl Let $\bc$ be the two-sided cell
of $W$ corresponding to a partition $\l$ of $n$ and $\mu$ the dual
partition of $\l$.}

\smallskip

\no(a){\sl The map $t_w\to t_{\phi_{ii}(w)}$ induces a ring isomorphism
 from $J_{\Ga_i\cap\Ga_i^{-1}}$ to $J_{\Ga_\l\cap\Ga_\l^{-1}}$.}

\no(b) {\sl The based ring $J_{\Ga_\l\cap\Ga_\l^{-1}}$ is commutative.}

\no(c) {\sl The map $$t_w\to E(t_{\phi_{ij}(w)},i,j), \qquad
 w\in A_{ij},$$
defines an isomorphism from the based ring $J_\bc$ to
$M_{n_\mu}(J_{\Ga_\l^{-1}\cap\Ga_\l})$, the $n_\mu\times n_\mu$
matrix algebra over the ring $J_{\Ga_\l^{-1}\cap\Ga_\l}$.}

\bigskip

\itp  (a) follows from Prop. 1.5.1 and its proof.

\medskip

(b) Let $\Ga_\bc$ be the canonical left cell in $\bc$, i.e.,
$\Ga_\bc$ is a left cell in $\bc$ and $R(w)\subseteq\{s_0\}$. By
[LX] we know that $J_{\Ga_\bc\cap\Ga_\bc^{-1}}$ is commutative.
Using (a) we see that (b) is true.

\medskip

(c) Note that for any $\omega,
\omega',\omega''$ in $\Omega$ and $x,y,z$ in $\bold c$ we have
$$\ga_{x,y,z}=\ga_{\omega x\omega'',{\omega''}^{-1}y\omega',\omega
z\omega'}.$$ Let $x\in A_{ij}, y\in A_{jk}$ and $z\in A_{ik}$. Using
Theorem 1.4.5 repeatedly
  we get
$$\ga_{x,y,z}=\ga_{\phi_{ij}(x),\phi_{jk}(y),\phi_{ik}(z)}.$$
Combining this and 1.3 (a) we see that (c) is true.

\medskip

The theorem is proved.

\bigskip

\no{\bf Remark:} For $\Phi_i$ if we choose different $h_i$ and/or
different sequence of left star operations then we would usually
get a different isomorphism in Theorem 2.3.2 (c).

For instance, let $n=3$ and $\bc$ the two-sided cell of $W$
corresponding to the partition $\l=(2,1)$. We have $s_\l=s_1$. Let
$\Phi_i$ ($i=0,1,2)$ be the right cell of $W$ containing $s_i$.
Then $\bc$ is the union of $\Phi_0,\Phi_1$ and $\Phi_2$. We have
$\Phi_0=\om^2\Phi_1$ and $\Phi_2=\om\Phi_1$. Let $*=\{s_1,s_0\}$
and $\star=\{s_1,s_2\}$. We also have $\Phi_0={}^*\Phi_1$ and
 $\Phi_2={}^\star\Phi_1.$ In this example it is easy to see that
 different $h_i$ and/or
different sequences of left star operations usually lead to
different isomorphisms in Theorem 2.3.2 (c).

\bigskip

Let $\l$ be a partition of $n$ and $\bc$ the two-sided cell of $W$
corresponding to $\l$. Let $\mu$ be the dual partition of $\l$ and
$u$ a unipotent element of $GL_n(\bbC)$ whose Jordan blocks are
given by the partition $\mu$. Denote by $F_\l$ a maximal reductive
subgroup of the centralizer of $u$ in $GL_n(\bbC)$. According to
Theorem 2.3.2 (c) and Lemma 2.3.1, to prove the conjecture of
Lusztig for $J_\bc$ we only need to prove the following special
case of the conjecture.

\bigskip

\no{\bf Conjecture 2.3.3.}  There is a bijection
  $\pi:
\Ga_\l^{-1}\cap\Ga_\l\to\text{Irr}F_\l$ such that

\smallskip

 \no (a) The map $t_w\to\pi(w)$ defines a ring isomorphism from
 $J_{\Ga_\l^{-1}\cap\Ga_\l}$ to $R_{F_\l}$.

\no (b) $\pi(w^{-1})=\pi(w)^*$ for any $w\in \Ga_\l^{-1}\cap\Ga_\l$. (Recall that $\pi(w)^*$ is the dual
of $\pi(w)$.)

\bigskip

We will prove this conjecture in Chapter 8. To define the map $\pi$ in
2.3.3 (a) we need some properties of antichains. In the following
section we give some discussions to chains and antichains.

\medskip

\section{Chains and antichains}

In Chapetr 5 we will define the map $\pi$ in Conjecture 2.3.3  by
means of r-antichains. In this section we prove some results about
chains and antichains which will be used later. In this section
$w$ stands for an element of $W$.

\bigskip

\noindent{\bf Lemma 2.4.1.} {\sl If $i$ and $j$ are in a d-chain
(resp. a d-antichain) of $w$, then any d-antichain (resp. d-chain) of
$w$ contains at most one of $i+an,j+bn$ for any given integers $a,b$.}

\bigskip

\itp Suppose that $i<j$ are in a d-chain of $w$. Then $w(i)>w(j)$.

\medskip

Assume that
 $i+an$ and $j+bn$ are in a d-antichain of $w$ for some integers $a,b$.
 If $$j+(b-1)n<i+an<j+bn,$$ then we have $b-a\le 0$ since $i<j$.
 Thus $$w(i+an)=w(i)+an>w(j+bn)=w(j)+bn$$  since  $w(i)>w(j)$. In this
 case
 $i+an$ and $j+bn$ can  not be  in the same  d-antichain of $w$. If
 $$i+an>j+bn>i+(a-1)n,$$ then $b-a\le -1$  since  $i<j$.
 Then   we
 have $$w(j+bn)=w(j)+bn<w(i+an-n)=w(i)+(a-1)n$$ since $w(i)>w(j)$. In
 this case   $i+an$ and $j+bn$ are also not in the same d-antichain of $w$.

\medskip

 Suppose that $i<j$ are in a d-antichain of $w$,
 then $j-n<i$ and $w(j)-n<w(i)<w(j)$.
 Assume that    $i+an$ and $j+bn$ are
 in a d-chain of $w$ for some integers $a,b$.
 If $i+an<j+bn$, then $b-a\ge 0$  since  $j-n<i$.
 Thus $$w(i+an)=w(i)+an<w(j+bn)=w(j)+bn$$  since  $w(i)<w(j)$. In this case
   $i+an$ and $j+bn$ can not be in the same d-chain of $w$. If
 $i+an>j+bn$, then $a-b\ge 1$  since  $i<j$. Then   we
 have $$w(j+bn)=w(j)+bn<w(i+an)=w(i)+an$$ since $w(i)>w(j)-n$. In
 this case both $i+an$ and $j+bn$ are also not in the same d-chain of $w$.

\medskip

 The lemma
 is proved.

 \bigskip

 \noindent{\bf Corollary 2.4.2.} {\sl If $w(i)$ and $ w(j)$ are in an r-chain (resp.
 r-antichain) of $w$, then any d-antichain (resp. d-chain)
  of $w^{-1}$ contains at most one
 of $w(i)+an, w(j)+bn$ for any given integers $a,b$.}

\bigskip

\itp Since $w(i)$ and $w(j)$ are in a d-chain (resp. d-antichain)
 of $w^{-1}$, the assertion
follows from Lemma 2.4.1.

\bigskip

\def\no{\noindent}
\noindent{\bf Lemma 2.4.3.} {\sl If $j_1<j_2<\cdots<j_k$ is a
d-antichain of $w$ of length $k$, then
$j_{k}-n<j_1<j_2<\cdots<j_{k-1}$ is also a d-antichain of length
$k$. (We shall say that $$j_1<j_2<\cdots<j_k$$ and
$$j_i-an<j_{i+1}-an<\cdots<j_{k}-an<j_1-(a-1)n<\cdots<
j_{i-1}-(a-1)n$$ are {\bf equivalent antichains} for any integer
$a$ and $1\le i\le k$.)}

\bigskip

\itp Since $$w(j_k-n)=w(j_k)-n>w(j_{k-1}-n)=w(j_{k-1})-n$$ and
$$j_k-n>j_{k-1}-n,$$ we see that the lemma is true.

\bigskip

\noindent{\bf Proposition 2.4.4.} {\sl Let $w\in W$ and
$\mu=(\mu_1,...,\mu_{r'})$ be the dual partition of $\l(w)$. Given
any $1\le k\le \l_1$ and any consecutive $n$ integers
$Y=\{a_1,a_1+1...,a_1+n-1\}$, we can find a d-antichain family set
of $w$ of index $k$ that is included in $Y$ and  has cardinality
$\mu_1+\cdots+\mu_k$. }

\bigskip

\itp By definition we have a d-antichain family $Z$ of $w$ of
index $k$ that has cardinality $\mu_1+\cdots+\mu_k$.  Using Lemma
2.4.3 we see that any d-antichain in $Z$ is equivalent to a
d-antichain of $w$ in $Y$. Noting that equivalent d-antichains
have the same congruence classes  modulo  $n$, we see that the
required d-antichain family set exists.  The proposition is
proved.

\bigskip

Of course this assertion is false for  chain family set. For
instance, if $w(i)=n-i+1+2(i-1)n$ for $1\le i\le n$, then
$\l(w)=(n)$ since $n<n-1+n<n-2+2n<\cdots <1+(n-1)n$ is a d-chain
of $w$ of length $n$. But $\{1,2,...,n\}$ is not a d-chain of $w$.

\bigskip

Let $w\in W$ and $\l(w)=(\l_1,...,\l_r)$. Assume that the dual
partition of $\l(w)$ is $\mu(w)=(\mu_1,...,\mu_{r'}).$ We call a
subset $Y$ of $\bbZ$ a {\bf complete d-chain (resp. d-antichain)
set} of $w$ if

\no(1) $Y$ contains $n$ numbers and any two numbers in Y are
noncongruent  modulo $n$,

\no(2) $Y$ is a disjoint union of $r$ d-chains (resp. $r'$
d-antichains) of $w$ of length $\l_1,...,\l_r$ (resp.
$\mu_1,...,\mu_{r'}$) respectively.

Similarly we define {\bf complete r-chain sets} and {\bf complete
r-antichain sets} of $w$.

\bigskip

\no{\bf Definition 2.4.5.} Let $w\in W$ and
$\mu=(\mu_1,...,\mu_{r'})$ the dual partition of $\l(w)$. We say
that $\{A_1,...,A_k\}$ is a {\bf complete d-antichain family} of
$w$ if the following four conditions are satisfied:

\no(1) all $A_i$ are d-antichains of $w$,

\no(2) $k=r'$,

\no(3) the cardinality of $A_i$ is $\mu_i$ for all $i$,

\no(4) the union of all $A_i$ is $\{1,2,...,n\}$. (Necessarily $A_i\cap
A_j=\emptyset$ if $i\ne j$.)

We say that $\{w(A_1),...,w(A_{r'})\}$ is a {\bf complete
r-antichain family} of $w$ if $\{A_1,...,A_{r'}\}$ is a complete
d-antichain family of $w$.

\bigskip

\noindent {\bf 2.4.6. Example:} Let $\l=(\l_1,...,\l_r)$ be a
partition of $n$. Then $w_\l$ has a complete d-antichain family.
In fact, for given $1\le i\le \l_1$, let $1\le h_i\le r$ be such
that $i\le \l_{h_i}$ but $i\ge \l_{h_i+1}$ (we set $\l_{r+1}=0)$.
Then define
$$A_{i}=\{i,\l_1+i,\l_1+\l_2+i,...,\l_1+\cdots+\l_{h_i-1}+i\}.$$
It is easy to see that $A_1,...,A_{\l_1}$ form a complete
d-antichain family of $w_\l$.

\bigskip

 In general, an element in $W$ does not have a
complete d-antichain family. For example, let $n=6$ and
$(w(1),w(2), w(3),w(4),w(5),w(6))=(6,3,10,7,8,11)$. Then
$w\in\Gamma_\l$, here $\l=(2,2,1,1),$ and $w$ does not have any
complete d-antichain family. (I am grateful to the referee for
providing a similar example.)

 \medskip

\section{Star operations for $W$}

In this section we give a result (see Lemma 2.5.2) about the star
operation. The following result is a partial generalization of [S,
Corollary 4.2.3].

\bigskip

\noindent{\bf Lemma 2.5.1.} {\sl Let $s\in S$ be a simple reflection and $w\in
W$. Suppose $a<b$ and $w(a)>w(b)$. Then $sw(a)>sw(b)$ if $sw\geq w$.}

\bigskip

\itp Write $$a=a_1+a_2n\text{\ \  and\ \ }b=b_1+b_2n,$$
 where $1\leq a_1,b_1\leq n$
and $ a_2,b_2\in \bbZ$. Set $$w(i)=p_i+r_in,\ \ \ \text{where }1\leq
p_i\leq n,\
\ r_i\in
\bbZ.$$

\medskip

Assume that $s=s_k$ for some $1\leq k\leq n-1$. Choose $1\leq
j_k,j_{k+1}\leq n$ such that $$w(j_k)=k+\xi n,\ \ w(j_{k+1})=k+1+\zeta
n,\ \ \
\xi,\zeta\in\bbZ.$$ If $i\ne j_k,j_{k+1}$, then $sw(i)=w(i)$ and
$p_i\ne k,k+1$. Thus $(p_i-k)(p_i-k-1)>0$. Therefore if $1\leq i<j\leq
n$ and either $i$ or $j$ is not in $ \{j_k,j_{k+1}\},$ then
$$\left|\left[\frac{w(j)-w(i)}n\right]\right|=\left|\left
[\frac{sw(j)-sw(i)}n\right]\right|.$$

\medskip

Now $sw(j_k)=k+1+\xi n$ and $sw(j_{k+1})=k+\zeta n$. If $j_k>j_{k+1},$
 we must have $\xi\geq \zeta+1 $ since $sw\geq w$. If $j_k<j_{k+1}$, we
have $\xi\geq\zeta$ for the same reason.

\smallskip

If neither $a_1$ nor $b_1$ is contained in $\{j_k,j_{k+1}\}$, then we
have $sw(a)=w(a)>w(b)=sw(b).$ If $a_1$ is either $j_k$ or $j_{k+1}$
and $b_1$ is not in $\{ j_k,j_{k+1}\}$, then $sw(a)$ is either equal
to $$k+1+\xi n+a_2n=w(a)+1$$ or $$k+\zeta n+a_2n=w(a)-1,$$ and
$$sw(b)=w(b)=p_{b_1}+r_{b_1}n+b_2n.$$ When $sw(a)=w(a)-1$, we must have
$\zeta+a_2>r_{b_1}+b_2$ or $\zeta+a_2=r_{b_1}+b_2$ and $k+1>p_{b_1}$
 since $w(a)>w(b)$. But $p_{b_1}\ne k,k+1$, so we
have $sw(a)>sw(b).$ Similarly we see $sw(a)>sw(b)$ if $a_1\ne
j_k,j_{k+1}$ and $b_1$ is either $j_k$ or $j_{k+1}$.

\smallskip

Now suppose $\{a_1,b_1\}=\{j_k,j_{k+1}\}$. If
$a_1=j_k<j_{k+1}=b_1$, then $$sw(a)=w(a)+1>w(b)-1=sw(b).$$ If
$a_1=j_{k+1}<j_k=b_1$ then $$sw(a)=k+\zeta n+a_2n=w(a)-1$$ and
$$sw(b)=k+1+\xi n+b_2n=w(b)+1.$$ we must have $\xi\geq \zeta+1$
since $sw\ge w$ . Then $\xi+b_2>\zeta+a_2$ since $b>a$, but this
is impossible since
  $w(a)>w(b).$ Thus we proved the lemma if $s=s_k,
k=1,2,...,n-1$.

\medskip

Assume that $s=s_0$. Choose $1\leq j_1,j_{n}\leq n$ such that
$$w(j_1)=1+\xi n,\ w(j_{n})=n+\zeta n,\ \ \ \xi,\zeta\in\bbZ.$$ If $i$ is
not in $\{ j_1,j_{n}\}$, then $sw(i)=w(i)$ and $p_i\ne 1,n$. Thus
$(p_i-1)(p_i-n)<0$. Therefore if $1\leq i<j\leq n$ and
$\{i,j\}\ne\{j_1,j_{n}\},$ then
$$\left|\left[\frac{w(j)-w(i)}n\right]\right|=\left|\left
[\frac{sw(j)-sw(i)}n\right]\right|.$$

\medskip

Note that $sw(j_1)=\xi n$ and $sw(j_{n})=n+1+\zeta n$. If $j_1>j_{n}$
we must have $\xi\leq \zeta+1 $ since $sw\geq w$. If $j_1<j_{n}$, we
have $\xi\leq\zeta$ for the same reason.

\smallskip

If neither $a_1$ nor $b_1$ is contained in $\{j_1,j_{n}\}$, then we
have $sw(a)=w(a)>w(b)=sw(b).$ If $a_1$ is either $j_1$ or $j_{n}$ and
$b\ne j_1,j_{n}$, then $sw(a)$ is either $$\xi n+a_2n=w(a)-1$$ or
$$n+1+\zeta n+a_2n=w(a)+1,$$ and $$sw(b)=w(b)=p_{b_1}+r_{b_1}n+b_2n.$$ We
must have $\xi+a_2>r_{b_1}+b_2$ if $sw(a)=w(a)-1$ since $w(a)>w(b)$.
But $p_{b_1}\ne 1,n$, so we have $sw(a)>sw(b).$ Similarly we see
$sw(a)>sw(b)$ if $a_1\ne j_1,j_{n}$ and $b_1$ is either $j_1$ or
$j_{n}$.

\smallskip

Now suppose $\{a_1,b_1\}=\{j_1,j_{n}\}$. If $a_1=j_n<j_{1}=b_1$,
then $sw(a)=w(a)+1>w(b)-1=sw(b)$. If $a_1=j_{1}<j_n=b_1$ we must
have $\xi+a_2>\zeta+b_2$ since $w(a)>w(b)$, but this is impossible
since $\xi\leq \zeta $ for the reason of  $sw\ge w$ and $a_2\leq
b_2$.

\medskip

The lemma is proved.

\bigskip

\noindent{\bf Lemma 2.5.2.} {\sl Let $i\in\bbZ$ and $*=\{s_i,s_{i+1}\}$.
We have}

\no (a) {\sl Suppose that $w(i)$ is between $w(i+1)$ and $w(i+2)$,
then $w$ is in $D_R(s_{i},s_{i+1}).$ Moreover,} $$w^*(a)=\cases
w(a),\quad&\text{\sl if }a\not\equiv i+1,i+2\text{(mod }n),\\
w(i+2),\quad&\text{\sl if }a=i+1,\\ w(i+1),\quad&\text{\sl if
}a=i+2.\endcases$$

\no(b) {\sl Suppose that $w(i+2)$ is between $w(i)$ and $w(i+1)$,
then $w$ is in $D_R(s_{i},s_{i+1}).$ Moreover,} $$w^*(a)=\cases
w(a),\quad&\text{\sl if } a\not\equiv i,i+1\text{(mod }n),\quad\ \
\\ w(i+1),\quad&\text{\sl if }a=i,\quad\ \ \\
 w(i),\quad&\text{\sl
if }a=i+1.\quad\ \ \endcases$$

\bigskip

\itp (a) If $w(i+1)<w(i)<w(i+2)$, using 2.1.3 (f) and Lemma 2.5.1
repeatedly we see $w=us_i$ for some $u$ with $us_i\ge u$ and
$us_{i+1}\ge u$. If $w(i+1)>w(i)>w(i+2)$, using 2.1.3 (f) and
Lemma 2.5.1 repeatedly we see $w=us_is_{i+1}$ for some $u$ with
$us_i\ge u$ and $us_{i+1}\ge u$. Therefore $w\in
D_R(s_{i},s_{i+1}).$ In both cases we have $w^*=ws_{i+1}$ so that
$w^*(a)=w(a)$ if $a\not\equiv i+1,i+2\text{(mod }n)$, and
$w^*(i+1)=w(i+2),\ w^*(i+2)=w(i+1)$.

\medskip

The proof of (b) is similar. The lemma is proved.


\chapter{Canonical Left Cells}


\def\bbN{\mathbb N}

In section 2.3 we proved that $J_{\Ga\cap\Ga^{-1}}$ is commutative
for any left cell $\Gamma$ of the extended affine Weyl group $W$
associated with  $GL_n(\bbC)$, using the fact that it is true for
canonical left cells. In this chapter we give some discussion to
canonical left cells in $W$. Although we do not really need the
results here for our main purpose, the discussion here maybe is
helpful for understanding other types. In other types the based
ring $J_{\Ga\cap\Ga^{-1}}$ is not commutative in general, one may
see this from [X3, Chapter 11],  but it is always commutative when
$\Ga$ is a canonical left cell (see [LX]). Lusztig conjectured
that $J_{\Ga\cap\Ga^{-1}}\simeq R_{F_\bc}$ if $\Ga$ is a canonical
left cell in a two-sided cell $\bc$ of an extended affine Weyl
group (see [L8]). If this is true, maybe $J_{\Ga\cap\Ga^{-1}}$ is
a key to understand $J_\bc$. One more reason to consider canonical
left cells is that the intersection $ {\Ga\cap\Ga^{-1}}$ has a
good presentation if $\Gamma$ is a canonical left cell (see [LX]).
This fact maybe can be used to prove that the bijection between
the set of two-sided cells of an extended affine Weyl group and
the set of unipotent classes in the corresponding algebraic group
(see [L8]) preserves partial orders.

\medskip

In section 3.1 we recall some reduced expressions for fundamental
weights, which can be found in [L2]. In section 3.2 we determine
the right cell containing a given dominant weight. In section 3.3
we discuss the shortest elements $m_x$ in the double cosets
$W_0xW_0$ for all dominant weights $x$. It is not so easy to
describe the partition associated with  $m_x$. In section 3.4 we
describe the distinguished involution in a canonical left cell of
$W$.

\medskip

\section{The dominant weights}

Let $W$ be as in Chapter 2, that is, $W$ is the extended affine
Weyl group associated with  $GL_n(\bbC)$. Define
$x_i=\tau_1\tau_2\cdots\tau_i$. Recall that $X$ is the subgroup of
$W$ generated by all $\tau_i$. Then the set $X^+=\{x\in X\ |\
l(w_0x)=l(w_0)+l(x)\}$ of dominant weights in $W$ consists of the
elements $x_1^{a_1}x_2^{a_2}\cdots x_{n}^{a_n}$,
$a_1,a_2,...,a_{n-1}\in\bbN $ and $a_n\in\bbZ$. For further
investigation we need some reduced expression for $x_i$, see [L2],

\bigskip

$$\begin{array}{rl}
x_1&=\om s_{n-1}s_{n-2}\cdots s_2s_1\\
 x_2&=\om^2 (s_{n-2}s_{n-3}\cdots s_2s_1)(s_{n-1}s_{n-2}
 \cdots s_2)\\
&\qquad\qquad...\\
 x_i&=\om^i (s_{n-i}s_{n-i-1}\cdots s_1)(s_{n-i+1}s_{n-i}\cdots s_2)\cdots
(s_{n-1}s_{n-2}\cdots s_{i})\\
 &\qquad\qquad  ...\\
 x_{n-2}&=\om^{n-2}(s_2s_1)(s_3s_2)\cdots (s_{n-1}s_{n-2})\\
 x_{n-1}&=\om^{n-1}s_{1}s_2\cdots s_{n-2}s_{n-1}\\
 x_n&=\om^n.
\end{array}$$

\bigskip

In particular, the length of $x_i$ is $(n-i)i$. From the reduced
expressions we see

\medskip

\no(a) $l(xy)=l(x)+l(y)$ if $x$ and $y$ are in $X^+$.

\medskip

\no(b) $x_is_i\le x_i$ for $1\le i\le n-1$ and $s_jx_i=x_i s_j\ge x_i$ if
$i\ne j$.

\medskip

\no(c) $l(x_is_is_{i+1}\cdots s_{n-1})=l(x_i)-(n-i)$ for $1\le i\le n-1$.

\medskip

\section{The right cell containing $x\in X^+$}

In this section we describe the partition associated with  a
dominant weight $x$ in $W$. This is equivalent to determine the
right cell containing $x$. Let $\Ga_\bc$ be the unique left cell
contained in a two-sided cell $\bc$ of $W$ with $R(\Ga_{\bold
c})\subseteq\{s_0\}$. For convenience we set $\Phi_{\bold
c}=\Ga_{\bold c}^{-1}$. Obviously we have

\bigskip

\no(a) If $x$ is a dominant weight in $W$ then  $x$ is contained in some $\Phi_\bc$.

\bigskip

Let $\l$ be a partition of $n$ and $\mu=(\mu_1,\mu_2,...,\mu_{r'})$ be
the dual partition of $\l$. Let $\bc$ be the two-sided cell
corresponding to $\l$.

\bigskip

\noindent{\bf Lemma 3.2.1.} {\sl Keep the notation above.
Let $\nu_1,...,\nu_{r'}$ be a permutation of $\mu_1,...,\mu_{r'}$. Set
$x_\nu=x_{\nu_1}x_{\nu_1+\nu_2}\cdots x_{\nu_1+\cdots+\nu_{r'}}$. Then
$x_\nu\in\Phi_\bc$.}

\bigskip

\itp Set $\nu_{10}=0$ and $ \nu_{1i}=\nu_1+\cdots+\nu_i$ for
$i=1,2,...,{r'}$. Define
$V_i=\{\nu_{1,i-1}+1,\nu_{1,i-1}+2,...,\nu_{1i}\}$ for
$i=1,2,...,{r'}$. For any $a\in V_i$ we have $x_\nu(a)=a+({r'}-i+1)n$.
Thus $V_i$ is a d-antichain of $x_\nu$.

If $a\in V_i$, $b\in V_j$ and $i<j$, then $x_\nu(a)>x_\nu(b)$.
Using Lemma 2.4.2 we see that $a+kn$ and $b+ln$ can not be in the
same d-antichain of $x_{\nu}$ for any integers $k,l$. Therefore
$V_1,V_2,...,V_{r'}$ must form a complete d-antichain family of
$x_\nu$ (see 2.4.5 for definition) and the corresponding partition
is $\mu$. So $x_\nu$ is in $\Phi_\bc$.

\def\medksip{\medskip}
\bigskip

\noindent{\bf Lemma 3.2.2.}  {\sl Let $I$ be a subset of
$N=\{1,2,...,n-1\}$. Then}

\smallskip

\no(a) {\sl The elements $\prod_{i\in I}x_i^{a_i}$ and
$ x_I=\prod_{i\in I}x_i$ are contained in the same right cell of $W$
if all $a_i\ge 1$. }

\no(b) {\sl  The elements $x_n^k\prod_{i\in I}x_i^{a_i} (\ a_i\ge
1,k\in\bbZ)$ are contained in the right cell containing $x_I$.}

\bigskip

\itp We can find a partition $\mu=(\mu_1,...,\mu_{r'})$ of $n$ and
a permutation $\nu_1,...,\nu_{r'}$ of $\mu_1,...,\mu_{r'}$ such
that $x_Ix_n=x_\nu$ (see Lemma 3.2.1 for definition). Let
$V_1,...,V_{r'}$ be as in the proof of Lemma 3.2.1. Then we see
that $V_1,...,V_{r'}$ form a complete d-antichain family for both
$\prod_{i\in I}x_i^{a_i}$ and $x_I$ if all $a_i\ge 1$. Therefore
(a) is true. (b) follows from (a) since $x_n$ is in the center of
$W$. The lemma is proved.

\bigskip

The two lemmas show that it is easy to determine the right cell
containing a given dominant weight.

\bigskip

\no {\bf Examples:}

\no (1) $x_1x_2\cdots x_{n-1}$ is contained in the lowest two-sided
 cell of $W$.

\no(2) $x_i$ and $x_{n-i}$ are contained in the same two-sided cell (we
set $x_0=e$, the neutral element of $W$). The partition corresponding
to the two-sided cell is (2,2,...,2,1,...,1), where 2 appears
min$\{i,n-i\}$ times.

\medskip

\section{The elements $m_x$}

It is well known that $W$ is the union of the double cosets $W_0xW_0$
$(x\in X^+$). For the double coset $W_0xW_0$ we have a unique element
$m_x$ of minimal length. The elements $m_x\ (x\in X^+)$ are
interesting since we have (see [LX])

\bigskip

\no(a) $\bigcup \Ga_\bc\cap\Ga_\bc^{-1}=\{m_x\ |\ x\in X^+\},$
where $\bold c$ runs through the set of two-sided cells of $W$.

\medskip

\no(b) $m_x\ne m_y$ if $x,y$ are in $X^+$ and $x\ne y$.

\medskip

\no(c) Let $w\in W$. Then $L(w)=R(w)\subseteq\{s_0\}$
 if and only if $w=m_x$ for some $x\in X^+$.

\bigskip

Since $x_n=\om^n$ is in the center of $W$, we have

\medskip

\no(d) $m_{xx_n^a}=x_n^am_x=\om^{an}m_x$ if $x\in X^+$ and $a$ is an
integer.

\bigskip

In this section we work out an explicit form for $m_x$.
Unfortunately the author is unable to describe the left cell
containing $m_x$. Fix a subset $I$ of $N=\{1,2,...,n-1\}$. Recall
$x_I=\prod_{i\in I}x_i$. We set $m_I=m_{x_I}$.

\def\a{\alpha}
\def\lan{\langle}
\def\ran{\rangle}

\bigskip

\noindent{\bf Lemma 3.3.1.} {\sl Let $I$ be  a subset of $N$ and $x=\prod_{i\in
I}x_i^{a_i}$. Then $m_x=xx_I^{-1}m_I$ if $a_i\ge 1$ for all $i$.}

\bigskip

{\it Proof.} Let $y=x^{-1}$, $y_I=x^{-1}_I$, and $m'_I=m^{-1}_I$.
Then $l(s_im'_I)=1+l(m'_I)$ for $1\le i\le n-1$. Obviously we have
$$ m'_Iy_I^{-1}ys_i\ge m'_Iy_I^{-1}y$$ for $i=1,2,...,n-1$. Thus
we only need to show that $l(s_im'_Iy_I^{-1}y)=1+l(m'_Iy_I^{-1}y)$
for $1\le i\le n-1$. Let $w\in W_0$ be such that $y_I=w^{-1}m'_I$.
Recall that $R$ is the root system of $W$. Let $R^+$ be the set of
positive roots in $R$ and $R^-=-R^+$. We have (see [IM])

\def\b{\beta}
$$l(s_iwy_I)=\sum_{\underset {\alpha\in R^+}{s_iw(\a)\in
R^-}}|\lan y_I,\a^\vee\ran+1|+\sum_{\underset {\alpha\in
R^+}{s_iw(\a)\in R^+}}|\lan
y_I,\a^\vee\ran|,$$
$$l(wy_I)=\sum_{\underset {\alpha\in
R^+}{w(\a)\in R^-}}|\lan y_I,\a^\vee\ran+1|+\sum_{\underset {\alpha\in
R^+}{w(\a)\in R^+}}|\lan y_I,\a^\vee\ran|.$$

If $s_iw\ge w$, we can find a unique $\beta\in R^+$ such that
$s_iw(\beta)\in R^-$ and $w(\beta)\in R^+$. Then by the formulas
above we have $|\lan y_I,\beta^\vee\ran+1|>|\lan
y_I,\beta^\vee\ran|$ since $l(s_iwy_I)>l(wy_I)$. Since $\lan
y_I,\beta^\vee\ran\le 0$, we necessarily have $\lan
y_I,\b^\vee\ran=0$. This implies that $\lan y,\b^\vee\ran=0$. By
the length formulas for $l(s_im'_Iy_I^{-1}y)=l(s_iwy)$ and for
$l(m'_Iy_I^{-1}y)=l(wy)$ we see
$l(s_im'_Iy_I^{-1}y)=1+l(m'_Iy_I^{-1}y)$.

If $s_iw\le w$, we can find a unique $\beta\in R^+$ such that
$s_iw(\beta)\in R^+$ and $w(\b)\in R^-$. Then by the formulas
above we have $|\lan y_I,\beta^\vee\ran+1|<|\lan
y_I,\beta^\vee\ran|$ since $l(s_iwy_I)>l(wy_I)$. Since $\lan
y_I,\beta^\vee\ran\le 0$, we necessarily have $\lan
y_I,\b^\vee\ran<0$. This implies that $\lan y,\b^\vee\ran<0$. By
the length formulas for $l(s_im'_Iy_I^{-1}y)$ and for
$l(m'_Iy_I^{-1}y)$ we see
$l(s_im'_Iy_I^{-1}y)=1+l(m'_Iy_I^{-1}y)$.

The lemma is proved.

\bigskip

By the lemma above, to get an explicit form for $m_x$ we only need
to get an explicit form for $m_I$. We need some notation. For
$1\leq i\le j\le n-1$ we define $s_{ij}=s_is_{i+1}\cdots s_{j}$.

\bigskip

\noindent{\bf Lemma 3.3.2.} {\sl Suppose $I=\{a_k>a_{k-1}>\cdots >a_1\}$ is a
subset of $N$.
 Then we have
$$m_I=x_Ix_{a_1}^{-1}w_kw_{k-1}\cdots w_{2}\om^{a_1},$$
where $w_i= s_{a_i,n-1} s_{a_i-1,n-2}s_{a_i-2,n-3}\dots
s_{a_{i-1}+1,n+a_{i-1}-a_i}$ for $i=k,k-1,...,2.$}

\bigskip

\itp It is easy to see that $w=x_{a_1}^{-1}w_kw_{k-1}\cdots
w_{2}\om^{a_1}$ is in $W_0$. Therefore $z=x_Iw$ is in the double
coset $W_0x_IW_0$. Using the reduced expressions in 3.1 and using
3.1 (a) and 3.1 (b), we see $l(z)=l(x_I)-l(w)$. Thus
$L(z)\subseteq L(x_I)\subseteq \{s_0\}$. We also need show that
$R(z)\subseteq \{s_0\}$. We can write down explicitly the action
of $z$ on $\{1,2,...,n\}$, which is, $$\begin{array}{rl}
 1&\to a_{k}+1 \\ 2&\to
a_k+2\\ &...\\ n-a_k&\to n\\ n-a_k+1&\to a_{k-1}+1+n\\ n-a_k+2&\to
a_{k-1}+2+n\\&...\\n-a_{k-1}&\to a_k+n\\ & ...\\ n-a_i+1&\to
a_{i-1}+1+\xi_in\\n-a_i+2&\to a_{i-1}+2+\xi_in\\&...\\
n-a_{i-1}&\to a_i+\xi_in\\ & ... \\ n-a_1+1&\to
1+\xi_1n\\n-a_1+2&\to 2+\xi_1n\\&...\\n &\to a_1+\xi_1n
\end{array}$$
where $\xi_i=k-i+1$. In other words, we have
$z(n-a_i+j)=a_{i-1}+j+\xi_in$ for $1\le i\le k$, $1\le j\le
a_i-a_{i-1}$ (we set $a_0=0$ and $a_{k+1}=n$). Thus $z(i)>z(j)$ if
$n\ge i>j\ge 1$. Using Lemma 2.5.1 we see $R(z)\subseteq
\{s_0\}$.

The lemma is proved.

\bigskip

\no{\bf Examples:}

\no (1) $m_N=x_Nw_0$.

\no(2) $m_{x_i}=\om^i$ for all $i$. Thus $m_{x_i^k}=x_i^{k-1}\om^i$ if
 $k\ge 1$.

 \no(3) Let $n\ge q\ge p\ge 1$. If $I=\{q,q-1,...,p+1,p\}$, then
$$m_I=x_qx_{q-1}\cdots x_{p+1}s_{q,n-1}s_{q-1,n-1}\cdots
s_{p+1,n-1}\om^p.$$

\bigskip

\noindent{\bf Lemma 3.3.3.} {\sl Let $I$ be a subset of $ \{1,2,...,n-1
\}$.
Suppose that $x=\prod_{i\in I}x_i^{a_i}$ and $a_i\ge 2$ for all $i$ in
$I$. Then $m_x\er x_I$.}

\bigskip

\itp We have $m_x=xx_I^{-1}m_I$. Now $x\er xx_I^{-1}\er x_I$ since
$a_i\ge 2$ for all $i$, we necessarily have $m_x\er x_I$ since $x\rl
m_x\rl xx_I^{-1}$. The lemma is proved.

\bigskip

Let $\mu$ be the dual partition of $\l$. By Lemma 3.3.3 we know
$m_{x_\mu^2}=x_\mu m_{x_\mu}$ is contained in the two-sided cell $\bc$
corresponding to $\l$. Using Lemma 2.5.2 and the table in the proof of
Lemma 3.3.2, it is easy to find the sequence of right star operations
and $i$ for passing $\Ga_\bc$ to $\om^i\Ga_\l\om^{-i}$.

\medskip

\section{The distinguished involutions}

Now we determine the distinguished involution in a canonical left
cell. The involutions look complicated.

\bigskip

\noindent{\bf Lemma 3.4.1.} {\sl Suppose that $m_x\ (x\in X^+)$ is a distinguished
involution and $$x=x_1^{a_1}x_2^{a_2}\cdots
x_{n-2}^{a_{n-2}}x_{n-1}^{a_{n-1}}x_n^{a_n}.$$ Then}

\smallskip

\no(a) {\sl $a_i=a_{n-i}$ for all $1\le i\le n-1$.}

\smallskip

\no(b) $a_1+2a_2+\cdots+(n-1)a_{n-1}+na_n=0$.

\bigskip

\itp (a) According to the proof of [LX, Theorem 3.5], we have
$m_x^{-1}=m_{w_0x^{-1}w_0}$. The required assertion then follows from
$m_x^{-1}=m_x$, $w_0x^{-1}w_0=x_1^{a_{n-1}} x_2^{a_{n-2}}\cdots
x_{n-2}^{a_{2}}x_{n-1}^{a_{1}}x_n^{a_n}$ and 3.3 (b).

\medskip

\no(b) Since   $m_x\in W'$, we have $x\in W'$. Using the reduced
expressions in 3.1 we see that (b) is true. Recall that $W'$ is
the subgroup of $W$ generated by all simple reflections
$s_0,s_1,...,s_{n-1},$ see 2.1.3 (c).

\medskip

The lemma is proved.

\bigskip

\noindent{\bf Theorem 3.4.2.} {\sl The following elements are all distinguished
involutions contained in canonical left cells.}

\smallskip

\no(a)\qquad\qquad {\sl $\om^{-(a-1)n}m_{x_{i_1}x_{i_2}\cdots x_{i_{a-1}}x_{n-i_{a-1}}\cdots
x_{n-i_2}x_{n-i_1}}$,

\no where $1\le i_1< i_2\cdots< i_{a-1}$ satisfies that
$i_{h+1}-i_h-(i_{h-1}-i_{h-2})\ge 0$ for $h=1,...,a-2$, (we set
$i_0=i_{-1}=0,)$ and $n-2i_{a-1}-(i_{a-2}-i_{a-3})\ge 0$. The
associated partition is
$$(a,...,a,a-1,...,a-1,a-2,...,a-2,...,a-h,...,a-h,...,2,...,2,1,...,1),$$
where $a$ appears $i_1$ times, $a-h$ appears
$i_{h+1}-i_h-(i_{h-1}-i_{h-2})$ times for $h=1,2,...,a-2$, and 1
appears $n-2i_{a-1}-(i_{a-2}-i_{a-3})$ times.}

\medskip

\no(b) {\sl \qquad\qquad $\om^{-2pn}m_{x_{i_1}^2x_{i_2}^2
\cdots x_{i_p}^2x_{n-i_p}^2\cdots x_{n-i_2}^2x_{n-i_1}^2}$,

\smallskip

\no where $1\le i_1< i_2\cdots< i_{p}$ satisfies
that $i_{h+1}-i_h-(i_{h}-i_{h-1})\ge 0$ for $h=1,...,p-1$, (we set
$i_0=i_{-1}=0,)$ and $ n-2i_{p}-(i_{p}-i_{p-1})\ge 0$. The associated
partition is
$$(a,...,a,a-2,...,a-2,...,a-2h,...,a-2h, ...,3,...,3,1,...,1),$$
where $a$ appears $i_1$ times, $a-2h$ appears
$i_{h+1}-i_h-(i_{h}-i_{h-1})$ times for $h=1,2,...,p-1$, and
 1 appears
$n-2i_{p}-(i_{p}-i_{p-1})$ times, where $a=2p+1$.}

\medskip

\no(c) {\sl \qquad\qquad $\om^{-(p+r)n}m_{x_{i_1}^2x_{i_2}^2
\cdots x_{i_p}^2x_{i_{p+1}}\cdots x_{i_{r}}x_{n-i_{r}}\cdots
x_{n-i_{p+1}}x_{n-i_p}^2\cdots x_{n-i_2}^2x_{n-i_1}^2}$,

\smallskip

\no where $1\le i_1< i_2\cdots< i_{p}<i_{p+1}<\cdots<i_r$
satisfies that $i_{h+1}-i_h-(i_{h}-i_{h-1})\ge 0$ for $h=1,...,p$,
(we set $i_0=i_{-1}=0,)$ and $i_{h+1}-i_h-(i_{h-1}-i_{h-2})\ge 0$
for $p+2\le h\le r-1$, and $ n-2i_{r}-(i_{r-1}-i_{r-2})\ge 0$ if
$r\ge p+2$, $n-2i_{p+1}\ge 0$ if $r=p+1$. The associated partition
is
$$(a,...,a,a-2,...,a-2p,...,a-2p,a-2p-1,...,a-2p-1,...,2,...,2,1,...,1),$$
where $a$ appears $i_1$ times, $a-2h$ appears
$i_{h+1}-i_h-(i_{h}-i_{h-1})$ times for $h=1,2,...,p$, $a-2p-1$
appears $i_{p+2}-i_{p+1}$ times if $r\ge p+2$, $a-2p-h$ appears
$i_{p+h+1}-i_{p+h}-(i_{p+h-1}-i_{p+h-2})$ times for $h=2,3,...,
a-2p-2$, 1 appears $n-2i_{r}-(i_{r-1}-i_{r-2})$ times if $r\ge
p+2$ and appears $n-2i_{p+1}$ times if $r=p+1$, where $a=p+r+1$.}

\bigskip

\itp (a) According to the proof of  3.3.2, the action of $$m=\om^{-(a-1)n}m_{x_{i_1}x_{i_2}\cdots x_{i_{a-1}}x_{n-i_{a-1}}\cdots
x_{n-i_2}x_{n-i_1}}$$ on $\{1,2,...,n\}$ is given by the following
table,
$$\begin{array}{rl} 1&\to n-i_1+1-(a-1)n\\
 2&\to n-i_1+2-(a-1)n\\&...\\i_1&\to n-(a-1)n\\
i_1+1&\to n-i_2+1-(a-2)n\\ i_1+2&\to n-i_2+2-(a-2)n\\...\\i_2&\to
n-i_1-(a-2)n\\
  &...\\
i_j+1&\to n-i_{j+1}+1-(a-1-j)n \\i_j+2&\to
n-i_{j+1}+2-(a-1-j)n\\&...\\i_{j+1}&\to n-i_{j}-(a-1-j)n\\
 & ...\\
i_{a-1}+1&\to i_{a-1}+1 \\i_{a-1}+2&\to i_{a-1}+2\\&...\\n-i_{a-1}&\to
n-i_{a-1}\\
 & ...
\\
n-i_j+1&\to i_{j-1}+1+(a-j)n \\n-i_j+2&\to
i_{j-1}+2+(a-j)n\\&...\\n-i_{j-1}&\to i_{j}+(a-j)n\\
 & ...\\
 n-i_1+1&\to 1+(a-1)n\\n-i_1+2&\to 2+(a-1)n\\&...\\n
&\to a_1+(a-1)n
\end{array}$$

\medskip

Set $i_0=0$. We define $$\eta_{jk}=\cases i_{j-1}+k&\text{ if
}1\le j\le a-1,\ 1\le k\le i_{j}-i_{j-1}\\ i_{a-1}+k&\text{ if }j=
a,\ 1\le k\le n-i_{a-1} \\ n-i_{2a-j}+k&\text{ if } a<j\le 2a-1,\
1\le k\le i_{2a-j}-i_{2a-1-j}.\endcases$$ Then we have
$$m(\eta_{jk})=\eta_{2a-j,k}+(j-a)n$$ for $1\le j\le 2a-1$ and all
$k$.

\medskip

We now define an order among all the pairs $(j,k)$ for which
$\eta_{jk}$ are defined. We say $(j,k)<(j',k')$ if one of the
following three cases happens, (1) $j$ is even and $j'$ is odd,
(2) both $j,j'$ have the same parity and $k<k'$, (3) both $j,j'$
have the same parity and $j>j'$, $k=k'$. It is easy to check that
this defines a total order on the set of such pairs $(j,k)$.
Arranging the pairs in the increasing order, then we have a
corresponding arrangement $\xi_1,\xi_2,...,\xi_n$ for all
$\eta_{jk}$ (i.e. if $\xi_l=\eta_{jk}$ and $\xi_{l'}=a_{j'k'}$,
then $l<l'$ if and only if $(j,k)<(j',k')$). For a fixed $j$ we
denote by $c_j$ the cardinality of the set consisting of all
$\eta_{jk}$. Set $$h=\ds_{1\le j\le 2a-3}(a-1-\left[\frac
j2\right])c_j.$$

\medskip

We shall use $\psi_j(w)$ for $w^*$ if $w\in D_R(s_j,s_{j+1})$ and
$*=\{s_j,s_{j+1}\}$. We also write $\psi_{ij}(w)$ ($i\le j$) for
$\psi_i\psi_{i+1}\cdots\psi_{j-1}\psi_j(w)$ if the latter is
defined. Using Lemma 2.5.2 we see that $m_1=\psi_{i_2+1,n-1}(m)$
is well defined. Moreover $l(m_1)=l(m)-(n-i_2-1)$ and
$l(\psi_{i_2}(m_1))=l( m_1)+1$. Set $m_0=m$. Suppose for
$j=1,...,k-1$ ($k\ge 2$) we have defined
$m_{j}=\psi_{p_j,n+j-2}(m_{j-1})$ with
$l(m_{j})=l(m_{j-1})-(n+j-p_j-1)$ and
$l(\psi_{p_j-1}(m_{j}))=l(m_{j})+1$. Then we define
$m_{k}=\psi_{p_k,n+k-2}(m_{k-1})$, here $p_k$ is chosen so that
$\psi_{p_k,n+k-2}(m_{k-1})$ is well defined and
$l(m_{k})=l(m_{k-1})-(n+k-p_k-1)$ and
$l(\psi_{p_k-1}(m_{k}))=l(m_{k})+1$. By the conditions
$i_{h+1}-i_h-(i_h-i_{h-1})\ge 0$ and
$n-2i_{a-1}-(i_{a-2}-i_{a-3})\ge 0$, such $p_1,p_2,...,p_k$ exist.
Continuing this process we finally get an element $m'$ whose
action on $\{h+1,h+2,...,h+n\}$ is given by $m'(h+l)=m(\xi_l)+qn$,
where $q=(a-1-\left[\frac j2\right])$ if $\xi_l=\eta_{jk}$. Thus
we have $$m'(h+l)=\eta_{2a-j,k}+(j-a)n+(a-1-\left[\frac
j2\right])n$$ if $\xi_l=\eta_{jk}$. Therefore we have
 $${m'}^{-1}(\eta_{jk})=h+l-(2a-j-1-\left[\frac {2a-j}2\right])n,$$
 where $l$ is defined by $\xi_l=\eta_{2a-j,k}$.
Applying the same sequence of right star operations
$\psi_{p_1,n-1},\psi_{p_2,n},...$ (in the same order) to ${m'}^{-1}$
we get an element $m''$ whose action on $\{h+1,h+2,...,h+n\}$ is the
same as the action of $w_I$ on, here $I$ is a suitable subset of
$\{h+1,h+2,...,h+n\}$, and $w_I$ is the longest element of the
subgroup generated by all $s_j$ ($j\in I$). Moreover, $\l(w_I)$ is the
partition in (a). Using Prop. 1.4.6 we see that (a) is true.

\medskip

(b) Using Lemma 3.3.1 and the proof of Lemma 3.3.2 we can write down
explicitly the action on $\{1,2,...,n\}$ of the element in (b). Then
as the proof of (a) we can see that the element in (b) is a
distinguished involution and its associated partition is the given one
in (b).

\medskip

The proof for part (c) is similar.

\medskip The theorem is proved.

\bigskip

\no {\bf Examples:}

\no The following elements are distinguished involutions, (1)
$\om^{-n}m_{x_ix_{n-i}}$ ($i\le\frac n2$), (2)
$\om^{-2n}m_{x_2^2x_{n-2}^2}$, (3) $\om^{-3n}
m_{x_1^2x_ix_{n-i}x_{n-1}^2}\ (i\le\frac n2)$, (4)
$\om^{-(n-1)n}$\break$\times m_{x_1^2x_2^2\cdots x_{n-1}^2}$. The
associated partitions are, (1) (2,...,2,1,...,1), where $2$
appears $i$ times, (2) (3,3,1,...,1), (3) (4,2,...,2,1,...,1),
where 2 appears $i-2$ times, (4) (n).


\chapter{The Group $F_{\lambda}$ and Its Representation}


\def\bbN{\mathbb N}
\def\vare{\varepsilon}

For later use we discuss the group $F_\l$ and its representations in
this chapter. In section 4.1 we give an explicit description for the
group $F_\l$. In section 4.2 we give some facts about the
representations of $F_\l$. For completeness we also supply a few
proofs although the facts are well known.

\medskip

\section{The group $F_\lambda$}

Let $W$ be as in Chapter 2, that is, $W$ is the extended affine
Weyl group associated with  $GL_n(\bbC)$. For each two-sided cell
$\mathbf c$ of $W$ we have a corresponding partition $\l$ of $n$.
Let $\mu=(\mu_1,\mu_2,...,\mu_{r'})$ be the dual partition of
$\l$. Let $u$ be a unipotent element in $GL_n(\bbC)$ whose Jordan
blocks are determined by the partition $\mu$. Choose $r'\ge
j_1>j_2>\cdots>j_p\ge 1$ such that

\no(1) $\mu_{j_1}<\mu_{j_2}<\cdots<\mu_{j_p}$,

\no(2) for any $1\le i\le r'$ we have $\mu_i=\mu_{j_k}$ for some $k$.

\no Let
$$n_k=\#\{i \ |\ 1\le i\le r'\text{ and } \mu_i=\mu_{j_k}\}.$$ Let
$C_G(u)$ be the centralizer of $u$ in $G=GL_n(\bbC)$. Then the maximal
reductive subgroup $F_\l=F_{\mathbf c}$ of $C_G(u)$ is isomorphic to
$GL_{n_1}(\bbC)\times GL_{n_2}(\bbC)\times\cdots\times
GL_{n_p}(\bbC)$. We shall identify $F_\l$ with $GL_{n_1}(\bbC)\times
GL_{n_2}(\bbC)\times\cdots\times GL_{n_p}(\bbC)$.

\bigskip

\no{\bf Examples:}

\no (1) If $\l=(n)$, then $\mu=(1,...,1)$. In this case
$F_\l=GL_n(\bbC)$.

\no(2) If $\l=(1,...,1)$, then $\mu=(n)$. In this case
$F_\l\simeq\bbC^*$.

\no(3) If $\l=(2,1,...,1)$, then $\mu=(n-1,1)$. Thus we have $F_\l\simeq
\bbC^*\times\bbC^*$ if $n\ge 3$.

\bigskip

\def\zd{\bbZ_{\text{dom}}}

 We define $\zd^n$ to be the set $\{(r_1,r_2,...,r_n)\in\bbZ^n\ |\
  r_1\ge r_2\ge\cdots\ge r_n \}$. The elements in $\zd^n$ will be
  called dominant elements in $\bbZ^n$. It is well known that the set
  Irr$GL_n(\bbC)$ of isomorphism classes of irreducible rational
  representations of $GL_n(\bbC)$ is one to one corresponding to
  $\zd^n$.

Thus the set Irr$F_\l$ of isomorphism classes of irreducible rational
  representations of $F_\l$ is one to one corresponding to
  $$\text{Dom($F_\l)=\zd^{n_1}\times\zd^{n_2}\times
  \cdots\times\zd^{n_p}$}.$$
  For an element
  $\vare=(\vare_{11},...,\vare_{1n_1},...,\vare_{p1},...,\vare_{pn_p})$
  in Dom$(F_\l)$ we also call $\vare_{ij}$ the ($i,j)$-component of
  $\vare$.

\medskip

\section{The representation ring of $F_\lambda$}

For any algebraic group $G$ we use $R_G$ for its (rational)
representation ring. We are interested in the representation ring
$R_{F_\l}$. Since $F_\l$ is isomorphic to a direct product of the
general linear groups $GL_{n_i}(\bbC)$ $(1\le i\le p)$, we see
that $R_{F_\l}$ is isomorphic to the tensor product (over $\bbZ)$
of the representation rings $R_{GL_{n_i}(\bbC)}\ (1\le i\le p)$.
Thus we are reduced to understand $R=R_{GL_n(\bbC)}$.

The structure of $R$ is known. For our purpose we list some properties
of $R$. We shall identify $X$ (resp. $X^+$) with $\bbZ^n$ (resp.
$\zd^n$) by identifying $\tau_i$ with the element in $\bbZ^n$ whose
$i$th component is 1 and other components are 0. For any $x\in X^+$,
we denote by $V(x)$ an irreducible representation of $GL_n(\bbC)$ with
highest weight $x$.

\bigskip
\def\no{\noindent}

\no(a) As a ring $R$ is generated by $V(x_1),...,V(x_{n-1}),$
$V(x_n),V(x_n^{-1})$, see 3.1 for definition of $x_i$.

\medskip

\no(a') As a ring $R$ is generated by $V(x_n^{-1})$ and $V(x_1^a)$ for $a=1,2,...,n$.

\medskip

\no(b) $R$ is a free $\bbZ$-module. The elements $V(x)\ (x\in X^+)$
form a $\bbZ$-basis of $R$, and the elements
$$V(x_1)^{a_1}V(x_2)^{a_2}\cdots V(x_{n-1})^{a_{n-1}}V(x_n)^{a_n},\
a_1,...,a_{n-1}\in\bbN,\ a_n\in\bbZ,$$ also form a $\bbZ$-basis of
$R$.

\medskip

\no(b') The elements $$V(x_1)^{a_1}V(x_1^2)^{a_2}\cdots
V(x_1^{n-1})^{a_{n-1}}V(x_n)^{a_n},\ a_1,...,a_{n-1}\in\bbN,\
a_n\in\bbZ,$$ form a $\bbZ$-basis of $R$.

\medskip

\no(c) For $1\le i\le n$, the dimension of $V(x_i)$ is $n\choose i$.

\medskip

\no(d) Each weight space of $V(x_i)$ has dimension one and the weight
set of $V(x_i)$ consists of all $\tau_{k_1}\tau_{k_2}\cdots\tau_{k_i}$
$(1\le k_1<k_2<\cdots<k_i\le n)$.

\medskip

\no(e) Let $x\in X^+$. Then $$V(x_i)\otimes V(x)\simeq\oplus_\tau
V(\tau x),$$ where $\tau$ runs through the set of weights of
$V(x_i)$ such that $\tau x\in X^+$.

\bigskip

{\it Proof.} Let
$\delta=\tau_1^{n-1}\tau_2^{n-2}\cdots\tau_{n-1}$. By (d),
$x\tau\delta\in X^+$ for any weight $\tau$ of $V(x_i)$. Using the
tensor product formula in [H, Ex. 24.9] we see that (e) is true.

\bigskip

\no (f) Let $K$ be a ring. Assume that  $K$ is a free
$\bbZ$-module with a basis $\{t_x\ |\ x\in X^+\}$. If for all $i$
we have $t_{x_i}t_x= \sum_\tau t_{\tau x}$, here $\tau$ runs
through the set of weights of $V(x_i)$ with $\tau x\in X^+$, then
we have a ring isomorphism between $R$ and $K$ defined by $V(x)\to
t_{x}$.

\bigskip

\itp It follows from (a), (b) and (e).

\bigskip

\no(g) Let $x\in X^+$ and $a\in\bbN$. Then $$V(x_1^a)\otimes
V(x)\simeq\oplus_\tau V(\tau x),$$ where $\tau$ runs through the
set of weights $\tau_1^{a_1}\tau_2^{a_2}\cdots\tau_n^{a_n}$ of
$V(x_1^a)$ with \break
$\tau_i^{a_i}\tau_{i+1}^{a_{i+1}}\cdots\tau_n^{a_n}x\in X^+$ for
$i=1,2,...,n$. (Note that each weight space of $V(x_1^a)$ has
dimension one and the weight set of $V(x_1^a)$ consists of all
$\tau_{1}^{a_1}\tau_2^{a_2}\cdots\tau_n^{a_n},$
$a_1,...,a_n\in\bbN$ and $a_1+a_2+\cdots+a_n=a$.)

\bigskip

\itp It follows from Littlewood-Richardson rule. When $a=2$, we also
can see (g) using [H, Ex. 24.9].

\bigskip

Recently Littelmann generalized the rule to arbitrary types, see [Li].

\bigskip

\no(h) Let $K$ be a ring. Assume that  $K$ is a free $\bbZ$-module
with a basis $\{t_x\ |\ x\in X^+\}$. If for all positive integers
$a$ we have $t_{x_1^a}t_x= \sum_\tau t_{\tau x}$, here $\tau$ runs
through the set of weights
$\tau_1^{a_1}\tau_2^{a_2}\cdots\tau_n^{a_n}$ of $V(x_1^a)$ with
$\tau_i^{a_i}\tau_{i+1}^{a_{i+1}}\cdots\tau_n^{a_n}x\in X^+$ for
$i=1,2,...,n$, then we have a ring isomorphism between $R$ and $K$
defined by $V(x)\to t_{x}$.

\bigskip

\itp It follows from (a'), (b') and (f).

\chapter{A Bijection Between $\Ga_\lambda\cap\Ga_\lambda^{-1}$ And Irr$F_\lambda$}



In this chapter we will establish a bijection between
$\Ga_\l\cap\Ga_\l^{-1}$ and the set of dominant weights of $F_\l$.
This bijection gives rise to a bijection between
$\Ga_\l\cap\Ga_\l^{-1}$ and Irr$F_\l$, that is in fact the
bijection $\pi$ in Conjecture 2.3.3. We use some r-antichains of
elements in $\Ga_\l\cap\Ga_\l^{-1}$ to define the bijection. To do
this we show that each element in $\Ga_\l\cap\Ga_\l^{-1}$ has a
complete r-antichain family, see 2.4.5 for definition.  The
contents of this chapter is as follows. In section 5.1 we show
that each element in $\Ga_\l\cap\Ga_\l^{-1}$ has a complete
r-antichain family, see Theorem 5.1.12. In section 5.2  we define
the map $\varepsilon$ from $\Ga_\l\cap\Ga_\l^{-1}$ to the set of
dominant weights of $F_\l$ by means of r-antichains of elements in
$\Ga_\l\cap\Ga_\l^{-1}$ and give some discussions to the map. The
main result of this chapter is Theorem 5.2.6, which says that the
map $\varepsilon$ is bijective. In section 5.3 we prove the main
theorem by describing elements of $\Ga_\l\cap\Ga_\l^{-1}$. In
section 5.4 we give some simple properties of elements in
$\Ga_\l\cap\Ga_\l^{-1}$. In section 5.5 we give some elements in
$\Ga_\l\cap\Ga_\l^{-1}$, most of them  correspond to fundamental
weights. In next 3 chapters we will show that this bijection
provides the isomorphism between $J_{\Ga_\l\cap\Ga_\l^{-1}}$ and
the rational representation ring $R_{F_\l}$.

\medskip

\section{r-antichains of elements in
$\Ga_\lambda\cap\Ga_\lambda^{-1}$}

\def\L{\Lambda}

Let $\l$ be a partition of $n$. Recall that we have defined the
element $w_\l$ in \S 2.2 and $\Ga_\l$ is the left cell of $W$
containing $w_\l$. In this section we show that each element in
$\Ga_\lambda\cap\Ga_\lambda^{-1}$ has a complete r-antichain
family. This is a rather delicate result since it does not hold
even for some elements in $\Gamma_\l$, see the example in 2.4.6.

\medskip

 Let $\l=(\l_1,...,\l_r)$ and
$\mu=(\mu_1,...,\mu_{r'})$ be its dual. Let $n_1,...,n_p$ be as in
4.1. The numbers $n_1,...,n_p$ can be defined in another way.
Choose $0=r_0<r_1<r_2<\cdots<r_k=r$ such that
$$\l_{r_i+1}=\l_{r_i+2}=\cdots=\l_{r_{i+1}}>\l_{r_{i+1}+1}=
\cdots=\l_{r_{i+2}}$$ for $i=0,1,...,k-2.$ Then $k=p$ and
$n_i=\l_{r_i}-\l_{r_{i+1}}$ for all $i=1,2,...,p$ (we understand
that $\l_{r_{p+1}}=0$). 

\bigskip

\noindent{\bf Lemma 5.1.1.} {\sl Let  $w\in \Ga_\l$. Define
$e_i=\l_1+\l_2+\cdots+\l_i$ for $1\leq i\leq r$ and $e_0=0$. Then}

\smallskip

\no(a) {\sl $w(e_i+1)>w(e_i+2)>\cdots>w(e_i+\l_{i+1})$ for
$i=0,1,...,r-1$.}

\smallskip

\no(b) {\sl  Given $0\leq h< i\leq r-1$ and  $1\leq j\leq \l_{i+1}$,
 we have  $$w(e_i+j)>w(e_{h}+\l_{h+1}-\l_{i+1}+j).$$
 In particular, we have
$$w(e_i+j)>w(e_{i-1}+\l_i-\l_{i+1}+j).$$}

\no(c) {\sl $w(e_{i-1}+k+n)>w(e_j+k)$ for $j=i,i+1,...,r-1$ and $1\le k\le\l_{j+1}$.}

\bigskip

\itp Since $w\in\Ga_\l$, $w=w_1w_\l$ for some $w_1$  with
$l(w_1w_\l)=l(w_1)+l(w_\l)$. Using Lemma 2.5.1 we see that (a) is
true.

\medskip

(b) Suppose that the assertion is not true. Then
  $$w(e_i+j)<w(e_{h}+\l_{h+1}-\l_{i+1}+j).$$ Thus
$$\begin{array}
{rl} &w(1)>w(2)>\cdots >w(\l_1),\\ &\ \ ...\\
&w(e_{h-1}+1)>w(e_{h-1}+2)>\cdots >w(e_{h-1}+\l_{h}),\\
&w(e_{h}+1)>w(e_{h}+2)>\cdots> w(e_{h}+\l_{h+1}-\l_{i+1}+j)\\
&\ \ \ \ >w(e_i+j)>\cdots
>w(e_i+\l_{i+1}),\end{array}$$ provides an r-chain family set of $w$ of
index $h+1$. The cardinality of the set is $e_{h+1}+1$ since the
lengths of the $h+1$ r-chains of $w$ in the r-chain family set are
$\l_1,...,\l_{h},\l_{h+1}+1$ respectively. This contradicts that
$\l(w)=\l$. (b) is proved.

\medskip

Similarly we prove (c). The lemma is proved.

\bigskip

\noindent{\bf 5.1.2.} Given $1\le j\le r$ and $1\le k\le \l_j$ we
set $$a_{jk}=e_{j-1}+k$$ and $$\L_j=\{a_{j1},\ a_{j2},\ ...,\
a_{j\l_j}\}.$$ We have  $w_\l(\L_j)=\L_j$ for all $j=1,2,...,r$.

Before going further we give some discussions to complete
r-antichian families.

Let $w$ be in $\Ga_\lambda\cap\Ga_\lambda^{-1}$ and assume that
$w$ has a complete d-antichain family (see 2.4.5 for definition).
By the argument for  Prop. 2.4.4,  we may decompose the set $
\{1,2,...,n\}$ into $r'$ d-antichains $A_1,A_2,$ $ ...,A_{r'}$ of
$w$ whose lengths are $\mu_1,...,\mu_{r'}$ respectively. By Lemma
5.1.1 and Lemma 2.4.1, any intersection $A_i\cap\L_j$ contains at
most one element.
 Thus we have

\bigskip

\no{\bf 5.1.2 (a)} The d-antichain $A_l$ of $w$ containing
$a_{ij}=e_{i-1}+j\ (1\le i\le r,\ 1\le j\le \l_i)$ has length $\ge
i$.

\bigskip

 The d-antichains $A_1,A_2,...,A_{r'}$
form a complete d-antichain family $Z$ of $w$ and $w(A_1),$ $
w(A_2),...,w(A_{r'})$ form a complete r-antichain family of $w$,
see \S 2.4.5 for definition.

\medskip

\def\dom{\text{Dom}(F_\l)}
Let $Z$ be  a complete r-antichain family  of $w$ and
    $\mu_{j_1}<\mu_{j_2}<\cdots<\mu_{j_p}$ be as in
section 4.1. Then we have $\mu_{j_i}=r_{i}$. Thus $Z$ contains
$n_i$ r-antichains of $w$ of length $r_i$. Let
$B_{i1},...,B_{in_i}$ be the r-antichains in $Z$ of length $r_i$.
Let $$b_{r_i,j}+c_{r_i,j}n>b_{r_i-1,j}+c_{r_i-1,j}n>\cdots
>b_{1,j}+c_{1,j}n$$ be elements in $B_{ij}$, where $1\le b_{k,j}\le n$ and
$c_{k,j}\in\bbZ$ for all $1\le k\le r_i$ and $1\le j\le n_i$.

\bigskip

\noindent{\bf Lemma 5.1.3.} {\sl Let $w$ be in
$\Ga_\l\cap\Ga_\l^{-1}$ and assume that $w$ has a complete
r-antichain family.
 Keep the notation in 5.1.2. For $1\le j\le n_i$,
 let $C_{ij}$ be the set consisting
of all $b_{k,j} \ (1\le k\le r_i) $. Then $C_{ij}$ contains
exactly one element of $\L_k$ if $1\le k\le r_i$ and contains no
element of $\L_k$ if $k>r_i$.}

\bigskip

\itp Using Corollary 2.4.2, any d-chain of $w^{-1}$ contains at
 most one element of
  $C_{ij}$. But $w^{-1}$
  is in $\Ga_\l$, so $\L_h$ is a d-chain of $w^{-1}$ for $h=1,...,r$.
  Thus $\L_h\cap C_{ij}$ contains at most one element for any $h,i,j$.
Note that $C_{ij}$ contains $r_i$ elements,
$n_i=\l_{r_{i}}-\l_{r_{i+1}}$
  and $h\le r=r_p$. Therefore $C_{pj}$ contains exactly one element of $\L_h$
  for all $h$ and $1\le j\le n_p$. This forces that each
  $C_{p-1,j}$ ($1\le j\le n_{p-1})$
  contains exactly one element of $\L_h$ if $1\le h\le r_{p-1}$ and
  contains no element if $h> r_{p-1}$. Inductively we see that
  $C_{ij}$ contains exactly one
element of $\L_k$ if $1\le k\le r_i$ and contains no element of
$\L_k$ if $k>r_i$. The lemma is proved.

\bigskip

\noindent{\bf Lemma 5.1.4.} {\sl Let $w$ be in
$\Ga_\l\cap\Ga_\l^{-1}$ and assume that $w$ has a complete
r-antichain family. Keep the notation in 5.1.2. Given $1\le i\le
p$ and $1\le j\le n_i$,
 if $b_{r_i,j}$ is in $ \L_k$, then}

\smallskip

 \no(a) {\sl $b_{r_i-h,j}$ is in $\L_{k-h}$ for all
 $1\le h< k$, and  $b_{r_i-k-h,j}$ is in $\L_{r_i-h}$ if $0\le h<r_i-k$.}

\no(b) {\sl $c_{r_i,j}=\cdots =c_{r_i-k+1,j}$, and
$c_{r_i-k,j}=\cdots =c_{1,j}=c_{r_i,j}-1$.}

\bigskip

 \itp (a) Let $b_q$ be the unique number in $C_{ij}\cap \L_q$ for
 $1\le q\le r_i$. Then $b_{k}=b_{r_i,j}$. We claim that
 $b_h=b_{r_i-h',j}$ for some $0\le h'<k$ if $1\le h\le k$. Otherwise,
 we have $b_h=b_{r_i-h',j}$ for some $1\le h\le k$ and $h'\ge k$. Then we
 can find some $k< q\le r_i,\ 1\le q'<k-1$ such that
 $b_q=b_{r_i-q',j}$. Since $$b_{r_i-q',j}=b_q>b_{r_i,j}=b_k>
 b_{r_i-h',j}=b_{h}$$ and

 $$b_{r_i,j}+c_{r_i,j}n>b_{r_i-q',j}+c_{r_i-q',j}n>b_{r_i-h',j}+
 c_{r_i-h',j}n,$$

\no we see that $$c_{r_i,j}>c_{r_i-q',j}\ge c_{r_i-h'}.$$ Thus
 $$b_{r_i,j}+c_{r_i,j}n-b_{r_i-h',j}-c_{r_i-h',j}n>n.$$ This is
 impossible since $B_{ij}$ is an r-antichain of $w$. So $b_h=b_{r_i-h',j}$ for
 some $0\le h'<k$ if $1\le h\le k$.

\medskip

For $1\le q<q'\le k$, we then have $0\le h,h'\le k-1$ such that
$b_{q}=b_{r_i-h,j}$ and $b_{q'}=b_{r_i-h',j}$. We claim that
$h>h'$. Otherwise, we have $h<h'$. Then $c_{r_i-h,j}>c_{r_i-h',j}$
since $b_{q}<b_{q'}$ and
$b_{q}+c_{r_i-h,j}n>b_{q'}+c_{r_i-h',j}n$. Thus
$b_{k}+c_{r_i,j}n>b_{q'}+c_{r_i-h',j}n+n$. This is impossible
since $B_{ij}$ is an r-antichain of $w$. So we have $h>h'$. Thus
we have $b_q=b_{r_i-k+q,j}$ if $1\le q\le k$. Similarly we see
$b_q=b_{q-k,j}$ if $k<q\le r_i$.

\medskip

(b) The assertion follows from the fact that $B_{ij}$  is an
r-antichain and (a).

\medskip

The lemma is proved.

\bigskip

Now we are going to show that each element in
$\Ga_\lambda\cap\Ga_\lambda^{-1}$ has a complete r-antichain
family. The way is long.

\bigskip

\no{\bf Definition 5.1.5.} We say that $w\in W$ is positive if
$w(m)$ is positive for any positive integer $m$.

\bigskip

\no{\bf Lemma 5.1.6.} {\sl Let
$w\in\Ga_\lambda\cap\Ga_\lambda^{-1}$ be positive. If
$w(1)+\cdots+w(n)=1+\cdots+n$, then $w=w_\l$.}

\bigskip

\itp By Lemma 2.5.1 we only need to show that $w(\L_j)=\L_j$ for
all $j=1,...,r.$ If the result is not true, we can find $1\le j\le
r$ such that $w(\L_i)=\L_i$ for $i=1,2,...,j-1$ and
$w(\L_j)\ne\L_j$. Then there exists $a$ in $\L_j$ such that
$w(a)=b$ is in   $\L_k$ for some $r\ge k>j$.

Note that $w^{-1}$ is also positive since  $1\le w(m)\le n$
whenever $1\le m\le n$. By Lemma 5.1.1 we have
$$w^{-1}(b)=a>w^{-1}(a_{k-1,\l_{k-1}})>\cdots>w^{-1}(a_{1,\l_1})>0.$$
Since $a$ is in $\L_j$, by our assumption on $w$, all the $k$
elements in the above sequences are contained in the union of
$\L_1,...,\L_j$. Now $k>j$, we can find some $1\le i\le j$ such
that $\L_i$ contains at least two elements in the above sequence.
This contradicts Corollary 2.4.2 since the above sequence is an
r-antichain of $w^{-1}$ and $\L_i$ is a d-chain of $w$.

Therefore $w(\L_j)=\L_j$ for all $j=1,2,...,r$. The lemma is
proved.

\bigskip

\no{\bf Lemma 5.1.7.} {\sl Let
$w\in\Ga_\lambda\cap\Ga_\lambda^{-1}$ be positive. Let
$a_{ij}\in\L_i$ be such that $w(a_{ij})>n$ and $w(a)\le n$ if
$a_{ij}<a\le n$. Then $w(a)=w_\l(a)$ if $a_{ij}<a\le n$.}

\bigskip

\itp Assume $a_{ij}<a\le n$ and $a\in\L_k$. We first show that
$w(a)$ is in $\L_k$. We use descend induction on $k$.

When $k=r$, by Lemma 5.1.1 and assumption on $a$ we get $$n\ge
w(a)>w(a_{r-1,\l_{r-1}})>\cdots>w(a_{1,\l_1})>0.$$ Using Corollary
2.4.2 we can see that $w(a)$ is in $\L_r$ and $w(a_{k,\l_k})$ is
in $\L_k$ for $k=1,2,...,r-1$.

Now suppose that $w(b)$ is in $\L_l$ if $b\in\L_l$ for some
$k<l\le r$. Using Lemma 5.1.1 we get   $$n\ge
w(a)>w(a_{k-1,\l_{k-1}})>\cdots>w(a_{1,\l_1})>0.$$ Since
$w(\L_l)=\L_l$ whenever $k<l\le r$, using Corollary 2.4.2 we see
that $w(a)$ is necessarily in $\L_k$.

\medskip

Now we show $w(a)=w_\l(a)$ if $a_{ij}<a\le n$. If $a\in\L_l$ for
some $i<l\le r$, we must have $w(a)=w_\l(a)$ since $w(\L_l)=\L_l$
and $w\in\Gamma_\l$.

Assume $a_{ij}<a\le n$ and $a\in\L_i$.  We only need to show that
the two sets $$A=\{w(a)\ |\ a_{ij}<a<a_{{i+1},1} \}$$ and
$$B=\{w_\l(a)\ |\ a_{ij}<a<a_{{i+1},1}\}$$ are equal. If this is
not true, then we can find $a_{ij}<a,b<a_{{i+1},1}$ such that
$w(a)\not \in B$ and $w_\l(b)\not\in A$. Since $w(\L_l)=\L_l$ if
$l>i$, we have $w(a')=w_\l(b)+cn$ for some $1\le a'\le a_{ij}$ and
integer $c$. Since $w$ is positive, $c$ is non-negative. We have
$w(a)>w_\l(b)$ since both $w(a)$ and $w_\l(b)$ are in $\L_{i}$ and
each number in $\L_{r_i}-B$ is greater than any number in $B$. But
$$w^{-1}(w(a))=a>a_{ij}\ge a'-cn=w^{-1}(w_\l(b)).$$ This
contradicts that $w^{-1}\in\Ga_\l$.

\medskip

The lemma is proved.

\bigskip

Let $w\in\Ga_\lambda\cap\Ga_\lambda^{-1}$ and $1\le a\le n$.
Assume that $w(a)>n$ and $w(b)\le n$ for all $a<b\le n$. By Lemma
5.1.1 we see that $a$ is $\L_{r_i}$ for some $r_i$, see the
beginning of this section for the definition of $r_i$.

\bigskip

\no{\bf Lemma 5.1.8.} {\sl Let $w\in\Ga_\lambda$ be positive and
$a_{r_ij}\in\L_{r_i}$.  We can find $\xi_k\in\L_k$ for
$k=1,2,...,r_i-1$ such that
$$w(a_{r_ij})>w(\xi_{r_{i-1}})>\cdots>w(\xi_2)>w(\xi_1)>w(a_{r_ij})-n.$$}

\bigskip

\itp Since $w\in\Ga_\l$ and $\l_{r_i}>\l_{r_i+1}$, we can find a
d-antichain family set $A$ of $w$ of index $\l_{r_i}$ with
cardinality $\mu_1+\cdots+\mu_{\l_{r_i}}$. By Prop. 2.4.4 we may
assume that $A$ is included in $\{1,2,...,n\}$. Let
$D_1,...,D_{\l_{r_i}}$ be the d-antichains of $w$ in $A$. Since
$\L_k$ is a d-chain of $w$ for any $k$, using Lemma 2.4.1 we see
that $A$ contains at most $\l_{r_i}$ elements of $\L_k$ for any
$k$. Since $A$ has cardinality $\mu_1+\cdots+\mu_{\l_{r_i}}$, $A$
contains exactly $\l_{r_i}$ elements of $\L_k$ if $1\le k\le r_i$
and contains all elements in $\L_k$ if $r_i\le k\le r$.  Thus each
d-antichain $D_h$ ($1\le h\le \l_{r_i}$) contains exactly one
element of $\L_k$ if $1\le k\le r_i$ and we can find a d-antichain
$D_l$ ($1\le l\le \l_{r_i}$) that contains $a_{r_ij}$. For $1\le
k<r_i-1$, let $\xi_k$ be the unique element in $D_l\cap \L_k$.
Then the chosen $\xi_k$'s satisfy our requirement. The lemma is
proved.

\bigskip

\no{\bf 5.1.9.} Let $w\in\Ga_\l\cap\Ga_\l^{-1}$ be positive and
$a_{r_ij}\in\L_{r_i}$. Assume that $w(a_{r_ij})>n$ and $w(b)\le n$
for all $a_{r_ij}<b\le n$. We want to define a positive element
$u$ in $\Ga_\l\cap\Ga_\l^{-1}$  with
$u(1)+\cdots+u(n)=w(1)+\cdots+w(n)-n$.

Set $w(a_{k,\l_k+1})=-\infty$ for all $k$. Choose $1\le
j_1\le\l_1$ such that $ w(a_{1j_1})>w(a_{r_ij})-n$ but
$w(a_{1,j_1+1})<w(a_{r_ij})-n$ if $j_1+1\le\l_1$. Then choose
$2\le j_k\le\l_k$ for $k=1,2,...,r_i-1$ such that
$$w(a_{k,j_k})>w(a_{k-1,j_{k-1}})>w(a_{k,j_{k}+1}) $$
 for $k=2,3,...,r_i-1$.  According to Lemma 5.1.8,
$a_{kj_k}$ ($1\le k\le r_i-1$) exists. Finally let $j_{r_i}=j$.
For simplicity we set $$a_k=a_{kj_k}\ \ \text{ for }k=1,...,r_i.$$

\bigskip

\noindent{\bf Lemma 5.1.10.} {\sl Keep the notation above. Then }
$w(a_{r_i})>w(a_{r_i-1})$.

\bigskip

\itp According to Lemma 5.1.8 we can find  $\xi_k$ ($1\le k\le
r_i-1$) in $\L_k$ such that
$w(a_{r_i})>w(\xi_{r_i-1})>\cdots>w(\xi_1)>w(a_{r_i})-n$. By the
definition of $a_k$ we have $w(\xi_k)\ge w(a_k)$ for
$k=1,2,...,r_i-1$. Therefore $w(a_{r_i})>w(a_{r_i-1})$. The lemma
is proved.

\bigskip

Keep the notation in 5.1.9. Now we define $u$ by $$u (a)=\cases
w(a_{r_i})-n& \text{if } a=a_1,\\w(a_{k-1}) & \text{if } a=a_k,\
2\le k\le r_i,\\ w(a)& \text{if } a\not\equiv a_{k}\text{(mod
}n),\ \text{ for all }1\le k\le r_i.\endcases $$ Clearly $u$ is
positive and $u(1)+\cdots+u(n)=w(1)+\cdots+w(n)-n$.

\bigskip

\noindent{\bf Lemma 5.1.11.} {\sl Let $u$ be as above. Then $u$ is
in $\Ga_\l\cap\Ga_\l^{-1}$.}

\bigskip

\itp We shall prove (1) $l(u)=l(uw_\l)+l(w_\l),$ (2)
$l(u^{-1})=l(u^{-1}w_\l)+l(w_\l),$ (3) $\l(u)= \l$.

\medskip

 Write $$w(a_k)=b_k+c_kn,\qquad
1\le b_k\le n,\ \ c_k\in \bbZ.$$ Then $b_k\in\L_{q_k}$ for some
$1\le q_k\le r_i$, since by Lemma 5.1.7, $w(\L_q)=\L_q$ if
$r_i<q\le r$. We also have $c_{r_i}\ge 1$ since $w(a_{r_i})>n$.

By Lemma 5.1.10 and the definition of $w(a_k)$, we see that
$w(a_1),...,w(a_{r_i})$ form an r-antichain of $w$. Note that
$w\in\Ga_\l\cap\Ga_\l^{-1}$. By Lemma 5.1.7,
$w(\L_q)=w_\l(\L_q)=\L_q$ whenever $q>r_i$. Using Corollary 2.4.2
we see that each $\L_q\ (1\le q\le r_i)$ contains exactly one
$b_k$ $(1\le k\le r_i)$. As in the proof of Lemma 5.1.4  we get

\medskip

\no($*)$ If $b_{r_i}$ is in $ \L_k$, then (a) $b_{r_i-h}$ is in
$\L_{k-h}$ for all
 $1\le h< k$, and  $b_{r_i-k-h}$
 is in $\L_{r_i-h}$ if $0\le h<r_i-k$;
 (b)  $c_{r_i}=\cdots =c_{r_i-k+1}$, and $c_{r_i-k}=\cdots
=c_{1}=c_{r_i}-1$.

\medskip

\def\wi{w^{-1}}

\medskip

(1) By the definition of $u$ and using Lemma 2.5.1, to check
$l(u)=l(uw_\l)+l(w_\l),$ we only need to verify that
$w(a_{r_i-1})>w(a_{r_i,j+1})$. If $b_{r_i}$ is not in $\L_1$, then
$c_{r_i-1}=c_{r_i}\ge 1$. In this case we have
$$w(a_{r_i-1})>w(a_{r_i,j+1})=w_\l(a_{r_i,j+1}).$$ If $b_{r_i}$ is
in $\L_1$, then $b_{r_i-1}$ is in $\L_{r_i}$ and
$c_{r_i-1}=c_{r_i}-1$. Since
$$w^{-1}(b_{r_i-1})=a_{r_i-1}-(c_{r_i}-1)n<\wi(w(a_{r_i,j+1}))
=a_{r_i,j+1},$$ $w^{-1}\in\Ga_\l$ and
$w(a_{r_i,j+1})=w_\l(a_{r_i,j+1})$ is in $\L_{r_i}$ (cf. Lemma
5.1.7), we have $b_{r_i-1}>w(a_{r_i,j+1})$. Therefore
$$w(a_{r_i-1})=b_{r_i-1}+(c_{r_i}-1)n>w(a_{r_i,j+1}),$$ We have
proved that $l(u)=l(uw_\l)+l(w_\l).$

\medskip

\def\at{&\text}
(2) Now we show that $l(u^{-1}w_\l)=l(u^{-1})-l(w_\l)$. We have
$$u^{-1}(b)=\cases a_{1}-(c_{r_i}-1)n\at{if } b=b_{r_i},
 \\a_{k+1}-c_kn\at{if }b=b_k,\ 1\le k\le r_i-1,\\w^{-1}(b)\at{if
 }b\not\equiv b_k\text{(mod }n), \ \text{ for all }1\le k\le r_i.\endcases$$

\def\ui{u^{-1}}
\def\wi{w^{-1}}
Let $b>b'$ be in the same $\L_q$ for some $q$. We need to show
that $u^{-1}(b)<\ui(b')$. Note that $w^{-1}(b)<w^{-1}(b')$ since
$w^{-1}\in\Ga_\l$. When both $b,b'$ are different from any $b_k$,
we have $$\ui(b)=w^{-1}(b)<w^{-1}(b')=\ui(b').$$ If $b>b'=b_k$ for
some $k$, then $b\ne b_h$ for any $1\le h\le r_i$ since $b,b'$ are
in the same $\L_q$. Then we have
$$\ui(b')=\ui(b_k)>w^{-1}(b_k)>w^{-1}(b)=\ui(b).$$

\medskip

Now suppose that $b_k=b>b'$. Then $b'\ne b_h$ for any $1\le h\le
r_i$ since $b,b'$ are in the same $\L_q$. Since $\wi\in\Ga_\l$ we
have $$w^{-1}(b_k)=a_k-c_kn<\wi(b')=a-cn,$$ where $1\le a\le n$
and $c\in\bbZ$. Thus $c_k\ge c$.

\medskip

Assume that $c_k>c$. When $1\le k\le r_i-1$, we have
$$\ui(b_k)=a_{k+1}-c_kn<a-cn=\wi(b')=\ui(b')$$.


 Now suppose that $k=r_i$.
We have $$\ui(b_k)=\ui(b_{r_i})=a_{1}-(c_{r_i}-1)n.$$ We claim
that $a_{1}<a$. Assume that $a\in\L_l$. If $l>1$ we of course have
$a_{1}<a$. Now suppose that $l=1$. To see that $a_1<a$ in this
case it suffices to show that $w(a)<w(a_1)$ since $a,a_1$ are in
$\L_1$ and $w\in\Ga_\l$.
 If $b_{r_i}=b_k$ is in
$\L_{r_i}$, by $(*)$ we have $c_1=c_{r_i}$. In this case we have
$w(a)=b'+cn<w(a_1)=b_1+c_1n$ since $c<c_k=c_{r_i}=c_1$. If
$b_{r_i}\in\L_h$ for some $1\le h\le r_i-1$, then
$b_1\in\L_{h+1}$, so $b_1>b_{r_i}>b'$. In this case we also have $
w(a)=b'+cn<w(a_1)=b_1+c_1n$ since $c\le c_{r_i}-1=c_k-1=c_1$. We
have seen that $a_1<a$ when $k=r_i$. Noting that $c_k-1\ge c$,
therefore if $k=r_i$ we also have
$$u^{-1}(b_k)=a_1-(c_k-1)n<a-cn=u^{-1}(b').$$

\medskip

Assume that $c_k=c$. Then $a_{k}<a$. Suppose $a\in\L_l$. Since
$b=b_k,b'$ are in the same $\L_q$, by $(*)$, we have
$$w(a_{k-1})=b_{k-1}+c_{k-1}n<w(a)=b'+cn<w(a_k)=b_k+c_kn.$$ By the
definition of $a_k$, we see that $a$ is not in $\L_k$. Thus $1\le
k<l$. When $l>k+1$, we have $a>a_{k+1}$, so in this case
$$\ui(b_k)=a_{k+1}-c_kn<a-cn=\wi(b')=\ui(b').$$ If $l=k+1$, since
$$w(a_{k+1})=b_{k+1}+c_{k+1}n>w(a_k)=b_k+c_kn>b'+cn=w(a)$$ and
$w\in\Ga_\l$, we have $a>a_{k+1}$. In this case we also have
$$\ui(b_k)=a_{k+1}-c_kn<a-cn=\wi(b')=\ui(b').$$

\medskip

We have showed  that $l(\ui w_\l)=l(\ui)-l(w_\l)$.

\medskip

(3) Now we prove  $\l(u)=\l$. By (1) we see $\l(u)\ge\l$. If we
can show  $\l(u)\le \l$, then we are forced to have $\l(u)=\l$.

\def\a{\alpha}
\def\b{\beta}
\def\d{\delta}
\def\e{\eta}
\def\t{\theta}

Let $\{A_1,A_2,...,A_j\}$ be an r-chain family of $u$ (see \S 2.2
for definition). If we can construct an r-chain family
$\{B_1,B_2,...,B_j\}$ of $w$ such that the cardinalities of
$A_1\cup\cdots\cup A_j$ and $B_1\cup\cdots\cup B_j$ are the same,
then we are done since it implies $\l(u)\le\l(w)=\l$. Note that
$u(a_1),...,u(a_m)$ form an r-antichain of $u$.  According to
Corollary 2.4.2, we have

\medskip

\no($\star 1$) each $A_k$ contains at most one element of the set
$H=\{w(a_m)+ln\ |\ 1\le m\le r_i,\ l\in\bbZ\}=\{u(a_m)+ln\ |\ 1\le
m\le r_i,\ l\in\bbZ\}.$

\medskip

Since $A_1,...,A_j$ form an r-chain family of $u$, we have

\medskip

 \no ($\star 2$) for any $1\le m\le r_i$, at most one
element of the set $H_m=\{w(a_m)+ln\ |\ l\in\bbZ\}$ is contained
in the union $A_1\cup\cdots\cup A_j$.

\medskip

Assume that $A_k$ contains some element of $H$ if $1\le k\le h$
and $A_k$ does not contain any element of $H$ if $h+1\le k\le j$.
Then $A_{h+1},...,A_{j}$ are also r-chains of $w$. Set $B_k=A_k$
if $h+1\le k\le j$.

Now we consider the r-chains $A_1,...,A_h$ of $u$. Assume that
$w(a_{m_k})+l_kn$ ($1\le k\le h$) is contained in $A_k$. It is no
harm to assume $l_k=0$ for $k=1,...,h$. Otherwise we may replace
$A_k$ by $\{u(a)-l_kn\ |\ u(a)\in A_k\}$. We may further assume
$m_h>m_{h-1}>\cdots>m_1$.

Let
$$u(b_{k1})>\cdots>u(b_{kq_k})>u(a_{m_k})>u(b_{k,q_{k}+2})>\cdots$$
be the elements in $A_k$. We have two cases.

\medskip

\no Case 1. For any $1\le k\le h$, if $b_{kq_k}$ exists, (that is,
$u(a_{m_k})$ is not the greatest number in $A_k$,) then $b_{kq_k}$
is in $\L_{m_k}$.

\medskip

\no Case 2.   We can find some $1\le k\le h$ such that $b_{kq_k}$
exists and $b_{kq_k}$ is not in $\L_{m_k}$.

\medskip

 Set $a_{r_i+1}=a_1+n$.

 \medskip

Assume that we are in case 1. By ($\star
1), \ (\star 2)$, the definition of $a_{m_k}$ and Lemma 5.1.10, we
see $u(b_{kq_k})>u(a_{m_k+1})>u(a_{m_k})$ if $u(b_{kq_k})$ exists.
Let $B_k$ ($1\le k\le h$) be the set consisting of the elements
$$u(b_{k1})>\cdots>u(b_{kq_k})>u(a_{m_k+1})>u(b_{k,q_{k}+2})>
\cdots.$$ Note that $\L_{m_k}$ is a d-chain of $u$.  By the
definition of $u$ we see that $B_1,...,B_h$ are r-chains of $w$.
Clearly $B_1,...,B_j$ form an r-chain family of $w$ and the
cardinalities of $A_1\cup\cdots\cup A_j$ and $B_1\cup\cdots\cup
B_j$ are the same.

\medskip

Now assume that we are in case 2.  Let $k$ ($1\le k\le h$) be the
maximal such that either $b_{kq_k}$ does not exist or  $b_{kq_k}$
exists and is not in $\L_{m_k}$. We claim that $A_k$ is an r-chain
of $w$. When $m_k>1$, this is clear since $\L_{m_k-1}$ is a
d-chain of $w$. When $m_k=1$, we necessarily have $k=1$ since
$m_1<\cdots <m_{h-1}<m_h$. Noting that $u(a_1)$ is positive and
$u(b)=w(b)=w_\l(b)$ if $a_{r_ij}<b\le n$ (cf. Lemma 5.1.7), we see
$b_{1q_1}<a_{r_ij}-n$ since $b_{1q_1}$ is not in $\L_1$ when
$b_{1q_1}$ exists. Therefore $A_k$ is also an r-chain of $w$ when
$k=1$.

\medskip

If $k>1$, we would like to construct two r-chains $A'_{k-1}, A'_k$
of $w$ from $A_{k-1},A_k$ such that $A'_{k-1}\cup A'_k$ and
$A_{k-1}\cup A_k$ have the same cardinality.

\medskip

If $b_{k-1,q_{k-1}}$ does not exist or $b_{k-1,q_{k-1}}$ is not in
$\L_{m_{k-1}}$ when it exists, we set $A'_{k-1}=A_{k-1}$ and
$A'_k=A_k$.

\medskip

 Now assume that $b_{k-1,q_{k-1}}$ exists and is in $\L_{m_{k-1}}$.

\medskip

  If
  $m_{k-1}<m_k-1$, we set
$A'_{k-1}=(A_{k-1}-\{u(a_{m_{k-1}})\})\cup\{u(a_{m_{k-1}+1})\}$
and $A'_k=A_k$.

 Suppose that $m_{k-1}=m_k-1$ and
$b_{k-1,q_{k-1}}$ is in $\L_{m_{k-1}}$. If    $b_{kq_k}$ exists
and is in $\L_{m_{k-1}}$, we choose $1\le p_k\le q_k$ and $ 1\le
p_{k-1}\le q_{k-1}$ such that both $b_{k-1,p_{k-1}}$ and $
b_{kp_k}$ are in $\L_{m_{k-1}}=\L_{m_k-1}$ but neither
$b_{k-1,p_{k-1}-1}$ nor $ b_{k,p_k-1}$ is in $\L_{m_{k-1}}$. We
can move all elements in the intersection of
$u(\L_{m_k-1})=u(\L_{m_{k-1}})$ and $A_{k-1}\cup A_k$ to $A_k$ and
form two sets $A'_k, A'_{k-1}$ as follows. If
$u(b_{kp_k})>u(b_{k-1,p_{k-1}})$, we set
$$A'_k=A_k\cup\{u(b_{k-1,p_{k-1}}),...,u(b_{k-1,q_{k-1}})\},$$
$$A'_{k-1}=A_{k-1}-\{u(b_{k-1,p_{k-1}}),...,u(b_{k-1,q_{k-1}})\}.$$
If $u(b_{kp_k})<u(b_{k-1,p_{k-1}})$, let $A'_k$ be the set
consisting of $$u(b_{k-1,1}),\ u(b_{k-1,2}),\ ...,\
u(b_{k-1,q_{k-1}}),$$ $$ u(b_{kp_k}),\ ... ,\ u(b_{kq_k}) ,\
u(a_{m_k}),\ u(b_{k,q_{k}+2}),\ ...,$$ and let $A'_{k-1}$ be the
set consisting of $$
 u(b_{k1}),\ ...,\ u(b_{k,p_k-1}),\ u(a_{m_{k-1}}),\ u(b_{k-1,q_{k-1}+2}),
 \  ....$$

 Suppose that $m_{k-1}=m_k-1$ and
 $b_{k-1,q_{k-1}}$ is in $\L_{m_{k-1}}$. If $b_{kq_k}$ does not
 exist or  $b_{kq_k}$
 exists
 but $b_{kq_k}$ is not in $\L_{m_{k-1}}$, we choose
 $ 1\le p_{k-1}\le q_{k-1}$ such that
$b_{k-1,p_{k-1}}$ is in $\L_{m_{k-1}}=\L_{m_k-1}$ but
$b_{k-1,p_{k-1}-1}$  is not in $\L_{m_{k-1}}$. We can move all
elements in the intersection of $u(\L_{m_k-1})=u(\L_{m_{k-1}})$
and $A_{k-1} $ to $A_k$ and form two sets $A'_k, A'_{k-1}$ as
follows. If   $b_{kq_k}$ does not
 exist or   $b_{kq_k}$
 exists and  $u(b_{kq_k})>u(b_{k-1,p_{k-1}})$, we set
$$A'_k=A_k\cup\{u(b_{k-1,p_{k-1}}),...,u(b_{k-1,q_{k-1}})\},$$
$$A'_{k-1}=A_{k-1}-\{u(b_{k-1,p_{k-1}}),...,u(b_{k-1,q_{k-1}})\}.$$
If $b_{kq_k}$
 exists and $u(b_{kq_k})<u(b_{k-1,p_{k-1}})$, let $A'_k$ be the set
consisting of $$u(b_{k-1,1}),\ u(b_{k-1,2}),\ ...,\
u(b_{k-1,q_{k-1}})  ,\ u(a_{m_k}),\ u(b_{k,q_{k}+2}),\ ...,$$ and
let $A'_{k-1}$ be the set consisting of $$
 u(b_{k1}),\ ...,\ u(b_{kq_k}),\ u(a_{m_{k-1}}),\ u(b_{k-1,q_{k-1}+2}),
 \  ....$$

It is easy to see that in any case $A'_{k-1}$ and $A'_k$ are
r-chains of $w$. Moreover  $A'_{k-1}\cup A'_k$ and $A_{k-1}\cup
A_k$ have the same cardinality. If $k-1>1$, for the pair
$A'_{k-1}, A_{k-2}$, we apply the same process, then we get two
r-chains $A''_{k-1},A'_{k-2}$ of $w$ such that $A''_{k-1}\cup
A'_{k-2}$ and $A'_{k-1}\cup A_{k-2}$ have the same cardinality.
Continuing this procedure, we get $k$ r-chains
$A'_k,A''_{k-1},...,A''_2,A'_1$ of $w$. From the construction we
see that each of the r-chains $A'_k,A''_{k-1},...,A''_2,A'_1$ of
$w$ contains some $u(a_m)$ with $1\le m\le m_k$.

When $h=k>1$, we set
$B_k=A'_k,B_{k-1}=A''_{k-1},...,B_2=A''_2,B_1=A'_1$.

When $h=k=1$, we set $B_1=A_1$.

\medskip

Now assume $1\le k<h$. Recall that we are in case 2 and $A_k$ is
also an r-chain of $w$. When $k>1$, we have constructed  $k$
r-chains $A'_k,A''_{k-1},...,A''_2,A'_1$ of $w$. If $k=1$ we set
$A'_1=A_1$. By our assumption on $k$,   for any $k+1\le l\le h$,
  $b_{lq_l}$     exists and $b_{lq_l}$ is in $\L_l$. If $m_h<r_i$ or $u(a_1)$ is not in  $A'_1$,
then let $B_l$ ($k+1\le l\le h$) be the set consisting of  the
elements
$$u(b_{l1})>\cdots>u(b_{lq_l})>u(a_{m_l+1})>u(b_{l,q_{l}+2})>
\cdots.$$ And we set
$B_k=A'_k,B_{k-1}=A''_{k-1},...,B_2=A''_2,B_1=A'_1$.

\medskip

Assume that $m_h=r_i$ and $u(a_1)$ is in $A'_1$. We construct two
r-chains $A''_1,A'_h$ of $w$ as follows.


\medskip

 Note that $b_{hq_h}$ exists and is in
 $\L_{m_{h}}=\L_{r_i}$.
Let $$u(c_1)>\cdots>u(c_q)>u(a_1)>u(c_{q+2})>\cdots$$ be the
elements in $A'_1,$ and
$$u(d_1)>\cdots>u(d_{q'})>u(a_{r_i})>u(d_{q'+2})>\cdots$$ be the
elements in $A_h.$ We have   $d_{q'}=b_{hq_h}$.

Note that $u(a_1)$ is positive. We have $c_q< a_{r_ij}-n$ if
$u(c_q)$ exists,  since  $A'_1$ is an r-chain of $w$ and
$u(b)=w(b)=w_\l(b)$ if $a_{r_ij}<b\le n$ (cf. Lemma 5.1.7). If
$c_q$ exists and $c_q+n$ is in $\L_{r_i}$, we choose $1\le p\le q$
and $ 1\le p'\le q'$ such that both $c_p+n$ and $ d_{p'}$ are in
$\L_{r_i}=\L_{m_h}$ but neither $c_{p-1}+n$ nor $ d_{p'-1}$ is in
$\L_{r_i}$.  If $u(c_p+n)>u(d_{p'})$, we set
$$A''_1=A'_1\cup\{u(d_{p'}-n),\ ...,\ u(d_{q'}-n)\},$$
$$A'_{h}=A_{h}-\{u(d_{p'}),\ ...,\ u(d_{q'})\}.$$ If
$u(c_{p}+n)<u(d_{p'})$, let $A''_1$ be the set consisting of
$$u(d_1-n),\ u(d_2-n),\ ...,\ u(d_{q'}-n),\ u(c_p),\ ..., \
u(c_q),\    u(a_1),\ u(c_{q+2}),\ ...,$$ and let $A'_{h}$ be the
set consisting of $$
 u(c_1+n),\ ...,\ u(c_{p-1}+n),\ u(a_{m_h}),\ u(d_{q'+2}),\  ....$$

 If $c_q$ does not exist or
$c_q$ exists but $c_q+n$ is not in $\L_{r_i}$, we choose  $ 1\le
p'\le q'$ such that   $ d_{p'}$ is in $\L_{r_i}=\L_{m_h}$ but
 $ d_{p'-1}$ is not in $\L_{r_i}$.  If $c_q$ does not exists or
 $c_q$ exists and
$u(c_q+n)>u(d_{p'})$, we set $$A''_1=A'_1\cup\{u(d_{p'}-n),\ ...,\
u(d_{q'}-n)\},$$ $$A'_{h}=A_{h}-\{u(d_{p'}),\ ...,\ u(d_{q'})\}.$$
If    $c_q$ exists and $u(c_{q}+n)<u(d_{p'})$, let $A''_1$ be the
set consisting of $$u(d_1-n),\ u(d_2-n),\ ...,\ u(d_{q'}-n),\
u(a_1),\ u(c_{q+2}),\ ...,$$ and let $A'_{h}$ be the set
consisting of $$
 u(c_1+n),\ ...,\ u(c_{q}+n),\ u(a_{m_h}),\ u(d_{q'+2}),\  ....$$

It is easy to see that in any case $A''_1$ and $A'_h$ are r-chains
of $w$. Moreover  $A''_1\cup A'_h$ and $A'_1\cup A_h$ have the
same cardinality.  If $k<h-1$, for the pair $A'_{h}, A_{h-1}$, we
apply the same process as that for $A_k, A_{k-1}$ (in the case
$k>1$), then we get two r-chains $A''_{h},A'_{h-1}$ of $w$ such
that $A''_{h}\cup A'_{h-1}$ and $A'_{h}\cup A_{h-1}$ have the same
cardinality. Continuing this procedure, we get $h-k$ r-chains
$A''_h,A''_{h-1},...,A''_{k+2},A'_{k+1}$ of $w$. From the
construction we see that each of the $h-k$ r-chains of $w$
contains some $u(a_m)$ with $m_k<m\le r_i$. We set
$B_h=A''_h,...,B_{k+2}=A''_{k+2},B_{k+1}=A'_{k+1},$ and
$B_k=A'_k,B_{k-1}=A''_{k-1},...,B_2=A''_2,B_1=A''_1$.

From the construction it is clear that in any case $B_1,...,B_j$
form an r-chain family of $w$ and the cardinalities of
$A_1\cup\cdots\cup A_j$ and $B_1\cup\cdots\cup B_j$ are the same.

\medskip

The lemma is proved.

\bigskip

Now we are in the position to state the main result of this
section.

\bigskip

\no{\bf Theorem 5.1.12.} {\sl Each element of
$\Ga_\l\cap\Ga_\l^{-1}$ has a complete r-antichain family.}

\bigskip

\itp Let $w$ be an element in $\Ga_\l\cap\Ga_\l^{-1}$. We need
show that $w$ has a complete r-antichain family. We show the
result first in the case when $w$ is positive (see 5.1.5 for
definition), and then in general.

Suppose that $w$ is positive. We use induction on the sum
$E(w)=w(1)+\cdots+w(n)$. When $E(w)=1+\cdots+n$, by Lemma 5.1.6,
$w=w_\l$. According to 2.4.6, the result is true

Now suppose that $E(w)>1+\cdots+n$ and the result is true if
$u\in\Ga_\l\cap\Ga_\l^{-1}$  is positive and $E(u)<E(w)$. We can
find $r_i$ and $j$ such that $w(a_{r_ij})>n$ and $w(a)\le n$
whenever $a_{r_ij}<a\le n$. Let $u$ be as in Lemma 5.1.11. Then
$u\in\Ga_\l\cap\Ga_\l^{-1}$  is positive and $E(u)=E(w)-n<E(w)$.
By induction hypothesis, $u$ has a complete r-antichain family.

Keep the notation in 5.1.9. Let $Z$ be a complete r-antichain
family of $u$. We may require that the r-antichain in $Z$
containing $u(a_{r_i})$ has length $r_i$ for the reason explained
below. If $u(a_{r_i})>n$, then the r-antichain in $Z$ containing
$u(a_{r_i})$ has length $r_i$, since $u(\L_q)=w(\L_q)=\L_q$ when
$r_i<q\le r$. If $u(a_{r_i})\le n$, by Lemma 5.1.7, we have
$u(a)=w_\l(a)$ if $a_{r_i}\le a\le n$. Since $w(a_{r_i})>n$, by
Lemma 5.1.1, we see $\l_{r_i}-j\ge \l_{r_i+1}$. Thus at least one
r-antichain $A$ in $Z$ that contains $u(a)$ for some $a_{r_i}\le
a\le a_{r_i,\l_{r_i}}$ and has length $r_i$. Let $B$ be the
r-antichain in $Z$ containing $u(a_{r_i})$. According to Lemma
5.1.4, $(A-\{u(a)\})\cup \{u(a_{r_i})\}$ and
$(B-\{u(a_{r_i})\})\cup \{u(a)\}$ are r-antichains of $u$. Thus it
is harmless to assume that the r-antichain $B$ in $Z$ containing
$u(a_{r_i})$ has length $r_i$.

Assume that $\xi_k\in\L_k\ (1\le k<r_i-1)$ and $u(\xi_k)$ is in
$B$ for all $k$. Note that $u(a_{r_1}),...,u(a_1)$ form an
r-antichain of $u$ and $u(\L_q)=\L_q$ if $r_i<q\le r$. Let $B_k$
be the r-antichain in $Z$ containing $u(a_k)$. According to Lemma
5.1.4 and the assertion ($*)$ in the proof of Lemma 5.1.1, the
following sets are also r-antichains:
$B'=\{u(a_{r_i}),...,u(a_1)\}$,
$B'_k=(B_k-\{u(a_k)\})\cup\{u(\xi_k)\}$ ($1\le k\le r_i-1$).
Replacing $B, B_k\ (1\le k\le r_i-1)$ in $Z$ by $B',  B'_k\ (1\le
k\le r_i-1)$, respectively, we get a complete r-antichain family
$Y$ of $u$ such that the r-antichain
$B'=\{u(a_{r_i}),...,u(a_1)\}$ of $u$ is in $Y$. Replacing $B'$ in
$Y$ by the r-antichain $\{w(a_{r_i}),...,w(a_1)\}$ of $w$, we then
get a complete r-antichain family of $w$.

We have showed the result for positive elements in
$\Ga_\l\cap\Ga_\l^{-1}$. For any element $w$ in
$\Ga_\l\cap\Ga_\l^{-1}$, we can find a positive integer $a$ such
that $\omega^{an}w\in\Ga_\l\cap\Ga_\l^{-1}$ is positive. Noting
that  any d-antichain of $\omega^{an}w$ is also a d-antichain of
$w$, thus $w$ has a complete r-antichain family since
$\omega^{an}w$ has one.

The theorem is proved.

\bigskip

In next section we use this result to establish a bijection
between $\Ga_\l\cap\Ga_\l^{-1}$ and Dom$(F_\l)$.

\medskip

\section{A map from $\Ga_\lambda\cap\Ga_\lambda^{-1}$ to Dom$(F_\lambda)$}

In this section we establish the bijection between
$\Ga_\lambda\cap\Ga_\lambda^{-1}$ and  Dom$(F_\lambda)$ by means
of complete r-antichain family.

\no{\bf 5.2.1.} Let $w$ be in $\Ga_\lambda\cap\Ga_\lambda^{-1}$.
By Theorem 5.1.12, $w$ has a complete r-antichain family $Z$. In
5.1.2 we have seen that  $Z$ contains
$n_i(=\l_{r_i}-\l_{r_{i+1}})$ r-antichains of $w$ of length $r_i$.
Let $B_{i1},...,B_{in_i}$ be the r-antichains in $Z$ of length
$r_i$. Let $$b_{r_i,j}+c_{r_i,j}n>b_{r_i-1,j}+c_{r_i-1,j}n>\cdots
>b_{1,j}+c_{1,j}n$$ be elements in $B_{ij}$, where $1\le b_{k,j}\le n$ and
$c_{k,j}\in\bbZ$ for all $1\le k\le r_i$ and $1\le j\le n_i$. It
is no harm to assume that
$$b_{r_i,1}+c_{r_i,1}n>b_{r_i,2}+c_{r_i,2}n>\cdots
>b_{r_i,n_i}+c_{r_i,n_i}n.$$

Define $\vare_{ij}(Z)=\vare(B_{ij})=\sum_{1\le k\le r_i}c_{k,j}$.
From Lemma 5.1.4  and our assumption on the arrangement
$B_{i1},B_{i2},$...,$B_{in_i}$, we see
 $\vare_{ij}(Z)=\vare(B_{ij})\ge\vare_{i,j+1}(Z)=\vare(B_{i,j+1})$
  if $j+1\le n_i$.
Thus we have defined an element
$$\vare(Z)=(\vare_{11}(Z),...,\vare_{1n_1}(Z),...,\vare_{p1}(Z),...,
\vare_{pn_p}(Z))\in\text{Dom}(F_\l).$$

\medskip

We will show that $\vare(Z)$ is independent of the choice of the
complete r-antichain family $Z$. To do this, we construct a
particular complete r-antichain family of $w$ in next subsection.

\bigskip

\no{\bf 5.2.2.} Assume that
$A=(x_{11},x_{12},...,x_{1\l_1},x_{21},x_{22},...,x_{2\l_2},
...,x_{r1},...,x_{r\l_r})$ is in $\bbZ^n$. We call that $A$ is
$\l$-{\bf admissible} if

\no(1) any two components of $A$ are different,

\no(2) $x_{i1}>x_{i2}>\cdots
>x_{i\l_i}$ for $i=1,...,r$,

\no(3) given $1\leq h< i\leq r$,
$1\leq j\leq \l_{i}$, we have $x_{ij}>x_{h,\l_h-\l_i+j}$, in
particular, we have $x_{ij}>x_{i-1,\l_{i-1}-\l_{i}+j}$.

\bigskip

Let $A =(x_{11},...,x_{1\l_1},...,x_{r1},...,x_{r\l_r})$ be a
$\l$-admissible element of $\bbZ^n$. We define $\vare_{k,i,j}$
inductively for all $1\le k\le {r_i},\ 1\le i\le p,$ and $ 0\le
j\le n_i$, (see the beginning of   section 5.1 for the definition
of $r_i,\ p, n_i$). First we define $\vare_{k,i,0}=\infty$ for all
$k,i$. Let $x$ be the unique greatest number among all components
of $A$. Then $x=x_{r_i,j}$ for some $1\le i\le p$ and $1\le j\le
\l_{r_i}$. Let $\vare_{r_i,i,1}=x_{r_i,j}$. Assume that we have
defined $\vare_{k,i,1}$ for some $k\le r_i$. We then define
$\vare_{k-1,i,1}$ to be the maximal number in $\{x_{k-1,q}\ |\
x_{k-1,q}<\vare_{k,i,1},\ 1\leq q\leq \l_{k-1}\}$.

Suppose that we have defined $\vare_{l,h,q}$ for all $h$ in
$\{1,2,...,p\}$ and $1\leq l\le r_h$, $0\le q\le b_h$. Here $b_h$
are some nonnegative integers, (note that we have defined
$\vare_{l,h,0}=\infty$ for all $l,h$). Let $B\in \bbZ^{n'}$ be
obtained from $A$ by removing all the components $\vare_{l,h,q}$
($1\le q\le b_h$)   for some $n'$. Let $g_{p'}=x_{p'q'}$ be the
greatest number among all components of $B$. Then it is easy to
see $p'=r_{h'}$ for some $h'$. Suppose that we have defined $g_k$
for $k\le p'=r_{h'}$. Then we define $g_{k-1}$ to be the maximal
number in $$\{x_{k-1,q'} \text{ is a component of $B$}\ |\
x_{k-1,q'}<g_k,\ 1\leq q'\leq \l_{k-1}\}.$$ Then we define
$\vare_{k,h',b_{h'}+1}=g_k$ for any $1\le k\le r_{h'}$. Continuing
this way we define all $\vare_{k,i,j}$ for $1\le k\le r_i,\ 1\le
i\le p$ and $ 1\le j\le n_i$.

\bigskip

\no{\bf Example:} Let $\l=(4,3,2,2)$. Then $r_1=1,\ r_2=2,\
r_3=4$. And $A=(11,7,4,3;12,6,5;$ $10,8;14,9)$ is $\l$-admissible.
We have $$\vare_{4,3,1}=14,\ \ \vare_{3,3,1}=10, \ \
\vare_{2,3,1}=6, \  \ \vare_{1,3,1}=4;$$ $$\vare_{4,3,2}=9,\ \
\vare_{3,3,2}=8,\ \ \vare_{2,3,2}=5,\ \ \vare_{1,3,2}=3;$$
$$\vare_{2,2,1}=12,\ \ \vare_{1,2,1}=11;\qquad \vare_{1,1,1}=7.$$

\bigskip

Let $w\in\Ga_\l\cap\Ga_\l^{-1}$. Set $x_{ij}=w(e_{i-1}+j)$ for any
$1\leq i\leq r$ and $1\leq j\leq \l_i$. According to  Lemma 5.1.1
we see that
$A=(x_{11},...,x_{1\l_1},...,x_{r1},...,x_{r\l_r})\in\bbZ^n$ is
$\l$-admissible. Thus we can define $\vare_{k,i,j}$. To indicate
its relation with $w$, we shall denote it by $\vare_{k,i,j}(w)$.
We would like to understand the properties of $\vare_{k,i,j}(w)$.

\bigskip

\noindent{\bf Lemma 5.2.3.} {\sl Let $w$ be in
$\Ga_\l\cap\Ga_\l^{-1}$.  For given integers $1\le i\le p$ and
$1\le j\le n_i$, the numbers $\vare_{k,i,j}(w)\ (1\le k\le r_i)$
form an r-antichain of $w$ of length $r_i$. Thus the r-antichains
$$\{\vare_{k,i,j}(w)\ |\ 1\le k\le r_i\}\quad (1\le i\le p,\ 1\le
j\le n_i)$$ of $w$ form a complete r-antichain family, we denote
it by $Z_w$.}

\bigskip

\itp For simplicity we write $\vare_{k,i,j}$ instead of
$\vare_{k,i,j}(w)$ in this proof. Let $Z$ be a complete
r-antichain of $w$. If $\vare_{r_i,i,j}$ is maximal among all
$w(1),...,w(n)$, using 5.1.2 (a), we see that
 the r-antichain $B$ in $Z$ containing $\vare_{r_i,i,j}$ has
 length $r_i$.
 Moreover by the definition of $\vare_{k,i,j}$ we see that the
 smallest number in $B$ is not greater than $\vare_{1,i,j}$. Thus
 $\vare_{r_i,i,j}-n<\vare_{1,i,j}.$ So all
 $\vare_{k,i,j}$ ($1\le k\le r_i$) form an r-antichain of $w$ of
 length $r_i$.

\medskip

 In general, we consider the r-antichains of $w$ in
  $Z$.

\medskip

Let $$I: w(\xi_s)>w(\xi_{s-1})>\cdots>w(\xi_1)$$ $$J:
w(\eta_t)>w(\eta_{t-1})>\cdots>w(\eta_1)$$ be two r-antichains of $w$
in $Z$. Then $\xi_h$ and $\eta_h$ are in $\L_h$ for all $h$. Suppose
$w(\xi_s)>w(\eta_t)$. If

$$\begin{array}{rl} w(\xi_{u+1})>w(\eta_u),&
w(\eta_u)>w(\xi_u),\ \ w(\eta_{u-1})>w(\xi_{u-1}),\\
...,&w(\eta_{u-v})>w(\xi_{u-v}),\ \ w(\eta_{u-v-1})<w(\xi_{u-v-1})
\end{array}$$

\no for some $0\leq v<u\le t,s-1$, then

$$\begin{array}{rl}
I':& w(\xi_s)>w(\xi_{s-1})>\cdots>w(\xi_{u+1})>w(\eta_u)\\
 &\qquad >w(\eta_{u-1})
>\cdots >w(\eta_{u-v})>w(\xi_{u-v-1})>\cdots >w(\xi_1),\\
& \\
J':& w(\eta_t)>w(\eta_{t-1})>\cdots>w(\eta_{u+1})>w(\xi_u)\\
 & >w(\xi_{u-1})>
\cdots >w(\xi_{u-v})>w(\eta_{u-v-1})>\cdots >w(\eta_1)\end{array}$$

\medskip

\no are two r-antichains of $w$. Replacing $I,J$ by $I',J'$
   in $Z$ respectively, we then get
 a
new complete r-antichain family $Z'$ of $w$. Continuing this process,
finally we get a complete r-antichain family $Z''$ of $w$ with the
following property. If $I: w(\xi_s)>w(\xi_{s-1})>\cdots>w(\xi_1)$ and
$J: w(\eta_t)>w(\eta_{t-1})>\cdots>w(\eta_1)$ are two r-antichains of
$w$ in $Z''$ and $w(\xi_s)>w(\eta_t)$, then $w(\xi_h)>w(\eta_h)$ for
all $1\le h\le s,t$ whenever $w(\xi_{h+1})>w(\eta_h)$. By the
definition of $\vare_{k,i,j}$ we see that
$$\vare_{r_i,i,j}>\vare_{r_i-1,i,j}>\cdots
>\vare_{1,i,j}$$ is exactly an r-antichain in $Z''$. The lemma is proved.

\bigskip

\noindent{\bf Lemma 5.2.4.} {\sl Let $w\in\Ga_\l\cap\Ga_\l^{-1}$
and $Z$ be a complete r-antichain family of $w$. Then the dominant
weight $\vare(Z)$ is independent of the choice of the complete
r-antichain family $Z$ of $w$ and only depends on $w$. Thus we can
denote the dominant weight by $\vare(w)$ and denote
$\vare_{ij}(Z)$ by $\vare_{ij}(w)$.}

\bigskip

\itp Let $I,J,I',J'$ be as in the proof of Lemma 5.2.3. We claim
that $\vare(I)=\vare(I')$ and $\vare(J)=\vare(J')$. By Lemma
5.1.4, $\vare(I)$ and $\vare(I')$ are completely determined by
$w(\xi_s$). So we always have $\vare(I)=\vare(I')$. If $u\ne t$,
then $w(\eta_t$) is the maximal number in $J$ and $J'$ as well. By
Lemma 5.1.4 we see in this case $\vare(J)=\vare(J')$. If $u=t$,
then $s>t$. Write $w(\xi_t)=b_1+c_1n$ and $w(\eta_t)=b_2+c_2n$,
where $1\le b_1,b_2\le n$ and $c_1,c_2\in \bbZ$. Using  Lemma
5.1.4 for the r-antichains $I$ and $I'$ we can find $k$ such that
both $b_1,b_2$ are in $\L_k$ and then $c_1=c_2$. Using Lemma 5.1.4
 we see $\vare(J)=\vare(J')$. Thus we have $\vare(Z)=\vare(Z')=\vare(Z_w)$. The
lemma is proved.

\bigskip

It seems that the number $\vare_{ij}(w)$ has a strange property
which now we are going to state. Define
$\vare'_{k,i,j}(w)=\vare_{k,i,j}(w)-w_\l(f)$ if
$\vare_{k,i,j}=w(f)$.

\bigskip

\noindent{\bf Lemma 5.2.5.} $\vare'_{k,i,j}(w)\geq
\vare'_{k,i,j+1}(w)$ {\sl if} $j+1\leq n_i.$

\bigskip

\itp Let $\vare_{k,i,j}(w)=w(e_{k-1}+f)$ and $\vare_{k,i,j+1}(w)
=w(e_{k-1}+g)$. Then $$\vare_{k,i,j}(w)-\vare_{k,i,j+1}(w)\geq
g-f= w_\l(e_{k-1}+f)-w_\l(e_{k-1}+g).$$ Therefore
$\vare'_{k,i,j}(w)\geq \vare'_{k,i,j+1}(w)$. The lemma is proved.

\bigskip

\no{\bf Remark:} It is likely that $\vare_{ij}(w)=\frac
1n\sum_{1\leq k\leq r_i}\vare'_{k,i,j}(w).$

\bigskip

The following is the main result of this chapter.

\medskip

\noindent{\bf Theorem 5.2.6.} {\sl Let $w$ be in
$\Ga_\l\cap\Ga_\l^{-1}$. We have}

\smallskip

\no(a)\qquad\qquad
$\vare(w^{-1})=(-\vare_{1n_1}(w),...,-\vare_{11}(w),...,-\vare_{pn_p}(w),...,-\vare_{p1}(w))$.

\smallskip

\no {\sl if }$\vare(w)=(\vare_{11}(w),...,
\vare_{1n_1}(w),...,\vare_{p1}(w),...,\vare_{pn_p}(w))$.

\smallskip

\no (b) {\sl The map $\vare:w\to \vare(w)$ defines a bijection between
$\Ga_\l\cap\Ga_\l^{-1}$ and} Dom$(F_\l)$.

\bigskip

\itp (a)   Let $B=B_{ij}\in Z_w$ be an r-antichain of $w$ of
length $r_i$ consisting of $$b_{r_i}+c_{r_i}n,\
b_{r_i-1}+c_{r_i-1}n,\ ...,\ b_1+c_1n,$$ where all $b_k$ are in
$[1,n]$ and $c_k$ are in $\bbZ$. Using Lemma 5.1.4 it is easy to
see that all $w^{-1}(b_k)=w^{-1}(b_k+c_kn)-c_kn$ form an
r-antichain $B'$ of $w^{-1}$ that is equivalent to $w^{-1}(B)$.
(See Lemma
 2.4.3 for the definition of equivalent r-antichains).
  All such r-antichains $B'$ of course form a complete r-antichain family
of $w^{-1}$. By definition we have
$$\vare(B')=-c_1-c_2-\cdots-c_{r_i}=-\vare(B).$$ By the definition of
$\vare(w)$ and of $\vare(w^{-1})$ we see that (a) is true.

\medskip

We will prove (b) in next section.

\bigskip

Before going further, we give two examples.

\medskip

\no(1) Assume that $\l=(n)$. Then the dual $\mu$ of $\l$ is
(1,...,1) and $w_\l=w_0$ is the longest element of
$W_0=<s_1,...,s_{n-1}>$. In this case we have
$$\Ga_\l\cap\Ga_\l^{-1}=\{w_0x\ |\ x\in X^+\}$$ and
$F_\l=GL_n(\bbC)$. Therefore Dom$(F_\l)=\bbZ^n_{\text{dom}}$. Let
$w$ be in $\Ga_\l\cap\Ga_\l^{-1}$. Then $w$ has a unique complete
r-antichain family, which consists of $\{w(1)\},...,\{w(n)\}$.
Assume that $w=w_0x^{a_1}_1\cdots x_{n-1}^{a_{n-1}}x_n^{a_n}$,
where $a_1,...,a_{n-1}\in\bbN$ and $a_n\in\bbZ$. Then we have
 $$\vare(w)= (a_1+a_2+a_3+\cdots+a_n,\ a_2 +a_3+\cdots+a_n,\ ...,
 \ a_{n-1}+a_n,\
 a_n),$$
 which is in $\text{Dom}(F_\l)=\bbZ^n_{\text{dom}}.$

\medskip

\no(2) Assume that $n\ge 3$ and $\l=(2,1,...,1)$. Then the dual $\mu$ of $\l$ is
($n-1$,1) and $w_\l=s_1$. In this case we have
$$\Ga_\l\cap\Ga_\l^{-1}=\{\om^{an}s_1(\om s_1)^b,\  \
\om^{an}s_1(\om^{-1} s_1)^b\ |\ a\in\bbZ,\ b\in\bbN\}$$
and $F_\l$ is isomorphic to $\bbC^*\times\bbC^*$. Therefore
Dom$(F_\l)=\bbZ^2$. We have
 $$\vare(\om^{an}s_1(\om s_1)^b)= (a,a(n-1)+b) \in\bbZ^2=\text{Dom}(F_\l)
 $$
 and
$$\vare(\om^{an}s_1(\om^{-1} s_1)^b)= (a,a(n-1)-b) \in\bbZ^2=\text{Dom}(F_\l)
 .$$

\def\medksip{\medskip}

\medskip

\section{Constructing elements of $\Ga_\lambda\cap\Ga^{-1}_\lambda$}

In this section we will give a proof of Theorem 5.2.6 (b). To do
this we need construct elements of $\Ga_\l\cap\Ga^{-1}_\l$. More
precisely for a given dominant weight
$\vare=(\vare_{11},...,\vare_{1n_1},...,\vare_{p1},$
$...,\vare_{pn_p})$ in Dom$(F_\l)$ we will construct an element
$w_\vare$ in $\Ga_\l\cap\Ga^{-1}_\l$ such that
$\vare(w_\vare)=\vare$. Then we show that the map $\vare$ is a
bijection from $\Ga_\l\cap\Ga^{-1}_\l$ to Dom$(F_\l)$.

\medskip

First we assume that all $\vare_{ij}$ are nonnegative.

\medskip

If $\vare=(1,0,...,0)\in$Dom$(F_\l)$, then $w_\vare$ is defined by
$$w_\vare(a)=\cases (i+1)\l_1 & \text{if } a=a_{i1},\ 1\le i\le
r_1-1,\\ \l_1+n& \text{if } a=a_{r_11},\\ w_\l(a)& \text{if }
a\not\equiv a_{i1}\text{(mod }n),\ \text{ for all }1\le i\le
r_1.\endcases $$ It is easy to see that $\l(w_\vare)\ge \l$  and
$\mu(w_\vare)\ge \mu$ (see 2.2 for the definition of $\l(w )$ and
$\mu(w)$).
 This forces that $\l(w_\vare)=\l$. It is easy to check that
$w_\vare$ is in $\Ga_\l\cap\Ga_\l^{-1}$ and
$\vare(w_\vare)=(1,0,...,0)$.

\medskip

Now Suppose that $\vare_{ij}\ge 1$ and we have defined $w=w_\vare\in
\Ga_\l\cap\Ga_\l^{-1}$ such that

 \no(1) $\vare(w_\vare)=\vare$, here
$\vare=(\vare_{11},...,\vare_{1n_1},...,\vare_{i1},...,
\vare_{i,j-1}, \vare_{ij}-1,0,...,0)$,

\no (2)
$w(a)=w_\l(a)$ if $a_{r_ij}<a\le n$, and $w(a_{r_ij})=w_\l(a_{r_ij})$
if $\vare_{ij}-1=0$.

\no(3) all $w(a)>0$ whenever $a>0$.

\no Then we define $u=w_{\vare'}$ for
$\vare'=(\vare_{11},...,\vare_{1n_1},...,\vare_{i1},...,
\vare_{i,j-1}, \vare_{ij},0,...,0)$ as follows. Set $j_{r_i}=j$
and choose $1\le j_{k-1}\le\l_{k-1}$ for $k=r_i,...,3,2$ such that
$$w(a_{k-1,j_{k-1}-1})>w(a_{k,j_k})>w(a_{k-1,j_{k-1}})$$ for
$k=r_i,...,3,2$, (we set $w(a_{k-1,0})=\infty$). Then define
$u=w_{\vare'}$ by $$w_{\vare'}(a)=\cases w(a_{k,j_k}) & \text{if }
a=a_{k-1,j_{k-1}},\ 2\le k\le r_i,\\ w(a_{1,j_1})+n& \text{if }
a=a_{r_i,j},\\ w(a)& \text{if } a\not\equiv a_{kj_k}\text{(mod
}n),\ 1\le k\le r_i.\endcases $$

\bigskip

\noindent{\bf Lemma 5.3.1.} {\sl $u=w_{\vare'}$ is in
 $\Ga_\l\cap\Ga_\l^{-1}$ and $\vare(u)=\vare'$. Moreover,
 $u(a)=w_\l(a)$ if $a_{r_ij}<a\le n$ and all $u(a)>0$ whenever $a>0$.}

\bigskip

\itp We shall prove (1) $l(u)=l(uw_\l)+l(w_\l),$ (2)
$l(u^{-1})=l(u^{-1}w_\l)+l(w_\l),$ (3) $\l(u)= \l,$ and (4)
$\vare(u)=\vare'$.

\medskip

Set $$a_{k}=a_{kj_k}\text{\ \ for $k=1,2,...,r_i$.}$$ Write
$$w(a_k)=b_k+c_kn, \ \ 1\le b_k\le n,\ c_k\in \bbZ.$$ Then
$b_k\in\L_{q_k}$ for some $1\le q_k\le r_i$, since by the assumption
we have $w(\L_q)=w_\l(\L_q)=\L_q$ if $r_i<q\le r$. We have
$$w(a_{r_i,j-1})=\xi+\eta n> b_{r_i}+c_{r_i}n=w(a_{r_ij})$$ since $w\in
\Ga_\l$, where $1\le
\xi\le n$ and $\eta\in\bbZ$.

It is easy to see that $w(a_1),...,w(a_{r_i})$ form an r-antichain
of $w$. Note that $w$ is in $\Ga_\l\cap\Ga_\l^{-1}$. Since
$w(\L_q)=w_\l(\L_q)=\L_q$ whenever $r\ge q>r_i$, using Corollary
2.4.2 we see that each $\L_q\ (1\le q\le r_i)$ contains exactly
one $b_k$ $(1\le k\le r_i)$. As in the proof of Lemma 5.1.4  we
get

\medskip

\no($*)$
If $b_{r_i}$ is in $ \L_k$, then (a) $b_{r_i-h}$ is in $\L_{k-h}$ for
all
 $1\le h< k$, and  $b_{r_i-k-h}$
 is in $\L_{r_i-h}$ if $0\le h<r_i-k$;
 (b)  $c_{r_i}=\cdots =c_{r_i-k+1}$, and $c_{r_i-k}=\cdots
=c_{1}=c_{r_i}-1$.

\medskip

(1) By the definition of $u$ and using 2.1.3 (f), to check
$l(u)=l(uw_\l)+l(w_\l),$ we only need to verify that $\xi+\eta
n>b_1+c_1n+n=w(a_1)+n$.

\medskip

When $j=1$ we need do nothing. Now assume that $j>1$, $\xi$ is in
$\L_k$ and $b_{r_i}$ is in $\L_{h}$. Note that $1\le k,h\le r_i$.

\medskip

If $k\le h$, then $\eta>c_{r_i}$. Otherwise, $\vare_{i,j-1}(w)\le
\vare_{ij}-1$. This contradicts that
$\vare_{i,j-1}(w)=\vare_{i,j-1}\ge\vare_{ij}$ since $\vare$ is in
Dom$(F_\l)$.

By $(*)$, we have $b_1\in\L_g$, where $g=h+1$ if $1\le h<r_i$ and
$g=1$ if $h=r_i$.

If $1\le h<r_i$, then $c_1=c_{r_i}-1$. So $c_1+1=c_{r_i}<\eta$. Thus
we have $\xi+\eta n>w(a_1)+n$ in this case.

If $h=r_i$ and $k>1$, then $g=1$ and $c_1=c_{r_i}$. In this case,
$\xi>b_1$ and $\eta>c_1$, we also have $\xi+\eta n>w(a_1)+n$.

When $h=r_i,\ k=1,$ we have $g=1$ and $c_1=c_{r_i}<\eta$. We then have
$$w^{-1}(b_1)=a_1-c_1n>w^{-1}(\xi)=a_{r_i,j-1}-\eta n.$$ Since
$w^{-1}\in\Ga_\l$, we see $\xi>b_1$. Therefore $\xi+\eta n>w(a_1)+n$.

We have showed  that if $k\le h$ then $\xi+\eta n>w(a_1)+n$.

\medskip

Now suppose that $k>h$. Then $\eta\ge c_{r_i}$ and $\xi>b_{r_i}$. We
have $g=h+1$ and $c_1=c_{r_i}-1$. If $k>h+1$, we have $\xi>b_1$, so
$\xi+\eta n>w(a_1)+n$. If $k=h+1$, we have
$$w^{-1}(b_1)=a_1-c_1n>w^{-1}(\xi)=a_{r_i,j-1}-\eta n.$$ Since
$w^{-1}$ is in $\Ga_\l$, thus $\xi>b_1$, so $\xi+\eta n>w(a_1)+n$. So
if $k>h$ we also have $\xi+\eta n>w(a_1)+n$.

\medskip

We have showed  that $\xi+\eta n>w(a_1)+n$ is true. So
$l(uw_\l)=l(u)-l(w_\l)$.

\medskip

\def\at{&\text}
(2) Now we show that $l(u^{-1}w_\l)=l(u^{-1})-l(w_\l)$. We have
$$u^{-1}(b)=\cases a_{r_i}-(c_1+1)n\at{if } b=b_1,
 \\a_{k-1}-c_kn\at{if }b=b_k,\ 2\le k\le r_i,\\w^{-1}(b)\at{if
 }b\not\equiv b_k\text{(mod }n), \ 1\le k\le r_i.\endcases$$

\def\ui{u^{-1}}
\def\wi{w^{-1}}
Let $b>b'$ be in the same $\L_q$ for some $q$. We need to show that
$u^{-1}(b)<\ui(b')$. Note that $w^{-1}(b)<w^{-1}(b')$ since $w^{-1}$
is in $\Ga_\l$. When both $b,b'$ are different from any $b_k$, we have
$$\ui(b)=w^{-1}(b)<w^{-1}(b')=\ui(b').$$ If $b=b_k>b'$ for some $k$,
then $b'\ne b_h$ for any $1\le h\le r_i$ since $b,b'$ are in the same
$\L_q$ and each $\L_q$ contains at most one element of
$\{b_1,...,b_{r_i}\}$. Then we have
$$\ui(b)=\ui(b_k)<w^{-1}(b_k)<w^{-1}(b')=\ui(b').$$

\medskip

Now suppose that $b>b_k=b'$. Then $b\ne b_h$ for any $1\le h\le
r_i$ since $b,b'$ are in the same $\L_q$ and $\L_q$ contains at
most one element of $\{b_1,...,b_{r_i}\}$. Since $\wi\in\Ga_\l$ we
have $$w^{-1}(b_k)=a_k-c_kn>\wi(b)=a-cn,$$ where $1\le a\le n$ and
$c\in\bbZ$. Thus $c_k\le c$.

\medskip

If $c_k<c$, we have $$\ui(b_k)=a_{k-1}-c_kn>a-cn=\wi(b)=\ui(b)$$ for
$2\le k\le r_i$. For $k=1$, we have
$$\ui(b_k)=\ui(b_1)=a_{r_i}-(c_1+1)n.$$ We claim that $a_{r_i}>a$.
Otherwise, $a>a_{r_i}$, then we must have $c=0$ since $w(a)=w_\l(a)$
for $a_{r_i}=a_{r_ij}<a\le n$ and then $c_1<0$ since $c_k=c_1<c$. By
assumptions on $w$, $c_1\ge 0$, a contradiction, so $a_{r_i}>a$. Thus
in this case we also have
$$\ui(b_k)=\ui(b_1)=a_{r_i}-c_1n-n>a-cn=\wi(b)=\ui(b).$$ We have
showed that $\ui(b_k)>\ui(b)$ if $c_k<c$.

\medskip

If $c_k=c$, then $a_k>a$. Note that $a_k\in\L_k$.

If $1\le k<r_i$, we have
$$w(a_{k+1})=b_{k+1}+c_{k+1}n>w(a)=b+cn>w(a_k)=b_k+c_kn,$$ since
$b,b_k$ are contained in the same $\L_q$ and $b_{k+1}$ is not in
$\L_q$ and $w(a_{k+1})>w(a_k)$. By the definition of $a_k$, we see
that $a$ is not in $\L_k$ in this case. If $k=r_i$, we also have that
$a$ is not in $\L_{k}$. Otherwise we have $a\in\L_{r_i}$. Then
$c>c_{r_i}$, since $b, b_k$ are in the same $\L_q$, $a_k>a$ and
$\vare_{i1}\ge\cdots\ge\vare_{i,j-1}>\vare_{ij}-1$. This contradicts
the assumption $c=c_k(=c_{r_i}$ in this case). Therefore in any case
$a$ is not in $\L_{k}$.

Since $a_k>a$, we see that $a$ is in $\L_h$ for some $1\le h<k$.
 If $h+1<k$, then
$$\ui(b_k)=a_{k-1}-c_kn>a-cn=\ui(b).$$ If $h+1=k$, then $a,a_{k-1}$ are
contained in $\L_h$. Since $b>b_k$ and $c=c_k$, we have
$$w(a)=b+cn>b_{k}+c_kn>b_{k-1}+c_{k-1}n=w(a_{k-1}).$$ Thus $a<a_{k-1}$
since $w\in\Ga_\l$. Therefore $$\ui(b_k)=a_{k-1}-c_kn>a-cn=\ui(b)$$ if
$k=h+1$.

\medskip

We have showed  that $l(\ui w_\l)=l(\ui)-l(w_\l)$.

\medskip

(3) Now we prove that $\l(u)=\l$. By (1) we see $\l(u)\ge\l$. If we
can show that $\mu(u)\ge \mu$, then we are forced to have $\l(u)=\l$.

\def\a{\alpha}
\def\b{\beta}
\def\d{\delta}
\def\e{\eta}
\def\t{\theta}

\medskip

Let $Z$ be a complete r-antichain family of $w$. If $w(a_{r_i})>n$, by
the assumptions on $w$, we see that the r-antichain in $Z$ that
contains $w(a_{r_i})$ has length $r_i$.

If $w(a_{r_i})<n$, then $\vare_{ij}(w)=\vare_{ij}-1=0$. By the
assumptions on $w$ we see that $w(a_{r_i})=w_\l(a_{r_i})$. Let $B$
be the r-antichain in $Z$ that provides $\vare_{ij}(w)$ (i.e.
$\vare(B)=\vare_{ij}(w)$ and $B$ has length $r_i$). Let $C$ be the
r-antichain in $Z$ that contains $w(a_{r_i})$. Assume that
$w(a_{r_il})$ is in $B$. Then by the assumptions on $w$, we have
$l\ge j$ and $w(a_{r_il})\in\L_{r_i}$. Replacing $w(a_{r_il})$ in
$B$ by $w(a_{r_i})$ and replacing $w(a_{r_i})$ in $C$ by
$w(a_{r_il})$, we get a complete r-antichain family $Z'$ of $w$.
The r-antichain in $Z'$ that contains $w(a_{r_i})$ has length
$r_i$.

Thus it is harmless to assume that the r-antichain in $Z$ that
contains $w(a_{r_i})$ has length $r_i$,

Let $$I: w(\xi_s)>w(\xi_{s-1})>\cdots>w(\xi_1),$$ $$J:
w(\eta_t)>w(\eta_{t-1})>\cdots>w(\eta_1)$$ be two r-antichains in $Z$.
Note that $\xi_k,\e_k$ are contained in $\L_k$. Assume that $I$
contains $w(a_{r_i})$. Then $s=r_i$ and $\xi_{r_i}=a_{r_i}$. Write
$$w(\xi_k)=\a_k+\b_kn,\ \ \ w(\e_k)=\d_k+\theta_kn,$$ where $1\le
\a_k,\d_k\le n$ and $\b_k,\t_k\in\bbZ$. Suppose that
 for some $1\le h\le k<r_i$ we
 have $$\begin{array}
 {rl}
 &w(\xi_{k+1})>w(\e_k)>w(\xi_k),\ w(\e_{k-1})>w(\xi_{k-1}),\\
 &..., \ w(\e_{h})>w(\xi_h),\ w(\e_{h-1})<w(\xi_{h-1}).\end{array}$$
We claim that

\medskip

\no(3a) $\a_k,\d_k$ are contained in the same $\L_q$ for
 some $q$ and $\b_k=\t_k$.

\medskip

Assume that $\a_k$ is in $\L_q$. Since
$w(\xi_{k+1})>w(\e_k)>w(\xi_k)$, we have $\b_{k+1}\ge\t_k\ge\b_k$.
If $1\le q\le r_i-1$, using Lemma 5.1.4, then $\a_{k+1}\in
\L_{q+1}$ and $\b_{k+1}=\b_k$. Thus $\a_{k+1}>\d_k$ and $\d_k$ is
in $\L_q$ or $\L_{q+1}$. In this case we have
$$w^{-1}(\d_k)=\eta_k-\t_kn=\eta_k-\b_kn<\xi_{k+1}-\b_kn=
w^{-1}(\a_{k+1}).$$ Since $w^{-1}\in\Ga_\l$, we must have
$\d_k\in\L_q$. Thus (3a) is true if $1\le q\le r_{i}-1$.

Now suppose that $q=r_i$. By Lemma 5.1.4, $\a_{k+1}\in \L_{1}$ and
$\b_{k+1}=\b_k+1$. Thus $\t_k=\b_k$ or $\b_{k}+1$. If
$\t_k=\b_k+1$, then we must have $\d_k\in\L_1$ and
$\a_{k+1}>\d_k$, since $w(\xi_{k+1})>w(\e_k)>w(\xi_k)$. But in
this case we have $w^{-1}(\d_k)<w^{-1}(\a_{k+1}).$ This
contradicts that $w^{-1}\in\Ga_\l$. Therefore $\t_k=\b_k$. Now
$1\le k\le r_i-1$ and by assumptions on $w$ we have
$w(\L_l)=w_\l(\L_l)=\L_l$ if $r_i<l\le r$. Thus $\d_k$ is not in
$\L_l$ for any $r_i<l\le r$. This forces that $\d_k\in
\L_{r_i}=\L_q$ since $w(\e_k)>w(\xi_k)$ and $\t_k=\b_k$. In
conclusion, (3a) is true if $q=r_i$.

We have seen that (3a) is always true.

\medskip

Thus the following two sequences $$\begin{array} {rl} &I':
w(\xi_s)>\cdots
>w(\xi_{k+1})>w(\e_k)>\cdots>w(\e_h)\\[4mm]
&\qquad\qquad>\cdots> w(\xi_{h-1})>\cdots >w(\xi_1),\\[4mm] &J':
w(\e_t)>\cdots> w(\e_{k+1})>w(\xi_k)>\cdots>w(\xi_h)\\[4mm]
&\qquad\qquad>\cdots> w(\e_{h-1})>\cdots
>w(\e_1)\end{array}$$
are two r-antichains of $w$. Replacing
 $I,J$ in $Z$ by $I',J'$ respectively,
 we get a new complete r-antichain family of $w$.
 Continuing this process, finally we can find a complete r-antichain
family $Y$ of $w$ with the following two properties

\medskip

\no(3b) The r-antichain $I$ in $Y$ that contains $w(a_{r_i})$ has
length $r_i$,

\medskip

\no(3c) Let $I:
w(\xi_s)>w(\xi_{s-1})>\cdots>w(\xi_1)$, $J:
w(\eta_t)>w(\eta_{t-1})>\cdots>w(\eta_1)$, be two r-antichains in $Y$
with $s=r_i$ and $\xi_{r_i}=a_{r_i}$. For any $1\le k\le s-1,t$, if
$w(\xi_{k+1})>w(\e_{k})$, then $w(\xi_k)>w(\e_k)$.

\medskip

\no By the definition of $a_k$, we see $\xi_k=a_k$.

\medskip

Replacing the r-antichain $w(a_{r_i})>\cdots>w(a_2)>w(a_1)$ in $Y$ by
$u(a_{r_i})>
\cdots>u(a_2)>u(a_1)$, we then get a set $Y'$ consisting of  $r'$ r-antichains
  of $u$. The lengths of the r-antichains of $u$ in $Y'$ are
  $\mu_1,...,\mu_{r'}$ respectively. Moreover the union of all
  elements in the r-antichains of $u$ in $Y'$ is
  $\{u(1),u(2),...,u(n)\}$. By definition we have
   $\mu(u)\ge\mu$. Therefore we have $\mu(u)=\mu$ and
$\l(u)=\l$.

\medskip

(4) From the proof above we see clearly $\vare(u)=\vare'$.

\medskip

From the definition of $u$ we see $u(a)=w_\l(a)$ if $a_{r_ij}<a\le n$
and all $u(a)>0$ whenever $a>0$.
\medskip

The lemma is proved.

\bigskip

\noindent{\bf Lemma 5.3.2.} {\sl Let $w$ be in
$\Ga_\l\cap\Ga_\l^{-1}$. If all $\vare_{ij}(w) $ are non-negative,
then $w(k)>0$ for all $1\le k\le n$.}

\bigskip

\itp Otherwise, we can find some $1\le k\le n$ such that $w(k)=a-bn$,
where $1\le a\le n$ and $b$ is a positive integer. Let $I$ be an
r-antichain of $w$ in a complete r-antichain family of $w$ that
contains $w(k)$. Let $c+dn$ ($1\le c\le n,\ d\in\bbZ$) be the largest
number in $I$. Since $\vare_{ij}(w)\ge 0$ for all $i,j$, we have $d\ge
1$. Then $c+dn>a-bn+n$. This contradicts that $I$ is an r-antichain of
$w$ containing $w(k)$. The lemma is proved.

\bigskip

\noindent{\bf Lemma 5.3.3.} {\sl Let $w\in\Ga_\l\cap\Ga_\l^{-1}$
be such that all components of $\vare_{ij}(w)$ are non-negative.
Let $1\le i\le p$ and $1\le j\le n_i$.
 If all
$\vare_{kl}(w) =0$ for $k>i$, or $k=i$ and $l>j$, then $w(a)=w_\l(a)$
whenever $a_{r_ij}<a\le n$.}

\bigskip

\itp By  Lemma 5.3.2, $w$ is positive. By the assumption on $w$ we
see that $w(a)\le n$ if $a_{r_ij}<a\le n$. Using Lemma 5.1.7 we
see that $w(a)=w_\l(a)$ if   $a_{r_ij}<a\le n$. The lemma is
proved.

\bigskip

\no{\bf 5.3.4.} Let $w$ be in $\Ga_\l\cap\Ga_{\l}^{-1}$. Suppose
that $\vare_{ij}\ge 1$
 and  all $\vare_{kl}$ are non-negative and
we have
$$\vare(w)=(\vare_{11},...,\vare_{1n_1},...,\vare_{i1},...,
\vare_{i,j-1}, \vare_{ij},0,...,0).$$ Then we define $u$ for
$$\vare'=(\vare_{11},...,\vare_{1n_1},...,\vare_{i1},...,
\vare_{i,j-1}, \vare_{ij}-1,0,...,0)$$ as follows. Choose $1\le
j_1\le\l_1$ such that $ w(a_{1j_1})>w(a_{r_ij})-n$ but
$w(a_{1,j_1+1})<w(a_{r_ij})-n$ if $j_1+1\le\l_1$. Then choose
$2\le j_k\le\l_k$ for $k=1,2,...,r_i-1$ such that
$$w(a_{k,j_k})>w(a_{k-1,j_{k-1}})>w(a_{k,j_{k}+1}) $$
 for $k=2,3,...,r_i-1$. (We  set
$w(a_{k,\l_k+1})=-\infty$ for all $k$). By Lemma 5.1.8, such
$j_1,...,j_{r_i-1}$ exist. Finally let $j_{r_i}=j$. For simplicity
we set $$a_k=a_{kj_k}\ \ \text{ for }k=1,...,r_i.$$

\bigskip

Now we define $u=u_{\vare'}$ by $$u (a)=\cases  w(a_{r_i})-n&
\text{if } a=a_1,\\w(a_{k-1}) & \text{if } a=a_k,\ 2\le k\le
r_i,\\ w(a)& \text{if } a\not\equiv a_{k}\text{(mod }n),\ \text{
for all }1\le k\le r_i.\endcases $$ Actually we have defined this
element in section 5.1 (see Lemma 5.1.11). According to Lemma
5.1.11 we get

\bigskip

\no{\bf 5.3.4 (a)} \ $u$ is in  $\Ga_\l\cap\Ga_{\l}^{-1}$.

\bigskip

From the proofs of Lemma 5.1.11 and of Theorem 5.1.2 we see
clearly

\bigskip

\no{\bf 5.3.4 (b)}\quad
$\vare(u)=\vare'=(\vare_{11},...,\vare_{1n_1},...,\vare_{i1},...,
\vare_{i,j-1}, \vare_{ij}-1,0,...,0)$.

\bigskip

\no{\bf 5.3.5.} {\it Proof of Theorem 5.2.6 (b):} Let
$\vare=(\vare_{11},...,\vare_{1n_1},...,\vare_{p1},...
,\vare_{pn_p})$ be in Dom$(F_\l)$.   Choose $k\in\bbN$ such that
$\vare_{ij}+kr_i\ge 0$ for all $1\le i\le p,\ 1\le j\le n_i$. By
Lemma 5.3.1 we can find $w\in\Ga_\l\cap\Ga_\l^{-1}$ such that
$$\vare(w)=(\vare_{11}+kr_1,...,\vare_{1n_1}+kr_1,...,\vare_{p1}+kr_p
,...,\vare_{pn_p}+kr_p) \in\text{Dom}(F_\l).$$ Then we have
$\om^{-kn}w\in\Ga_\l\cap\Ga_\l^{-1}$ and
$$\vare(\om^{-kn}w)=(\vare_{11},...,\vare_{1n_1},...,
\vare_{p1},...,\vare_{pn_p}).$$ Therefore the map $\vare$ is
surjective.

\medskip

Suppose that $\vare(w)=(0,...,0)\in$Dom$(F_\l)$. By Lemma 5.3.2,
$w$ is positive. By the definition of $\vare(w)$, we see $w(a)\le
n$ if $1\le a\le n$. Thus $w(1)+\cdots+w(n)=1+\cdots+n$. Using
Lemma 5.1.7 we get $w=w_\l$. Now suppose that
$w,w'\in\Ga_\l\cap\Ga_\l^{-1}$ and $$\vare(w)=\vare(w')=
\vare=(\vare_{11},...,\vare_{1n_1},...,\vare_{i1},...,\vare_{ij},0
...,0).$$ Suppose all $\vare_{kl}$ are nonnegative and
$\vare_{ij}\ge 1$. By subsection 5.3.4 we can construct $u,u'$
from $w,w'$ respectively such that
$$\vare(u)=\vare(u')=\vare=(\vare_{11},...,\vare_{1n_1},...,
\vare_{i,j-1},\vare_{ij}-1,0...,0).$$ We use induction on the sum
of all components of $\vare(w)$. By induction hypothesis we see
$u=u'$. Now we can recover $w,w'$ from $u,u'$ using the
construction at the beginning of this section, see Lemma 5.3.1.
Therefore $w=w'$ if all components of $\vare(w)=\vare(w')$ are
nonnegative. In general we can find $k$ such that all components
of $\vare(\om^{kn}w)=\vare(\om^{kn}w'$) are nonnegative. Thus
$\om^{kn}w=\om^{kn}w$. Hence $w=w'$ if $\vare(w)=\vare(w')$. We
proved that the map $\vare$ is injective.

\medskip

Theorem 5.2.6 (b) is proved.

\bigskip

 In the following section we give some   simple properties of elements
  in $\Ga_\l\cap\Ga_\l^{-1}$.

\medskip

\section{Some simple properties of elements in $\Ga_\lambda\cap\Ga_\lambda^{-1}$}

Recall that we defined $a_{ij}=\l_1+\cdots+\l_{i-1}+j$ for $1\le i\le
r$, $1\le j\le\l_i$, and $\L_k=\{a_{k1},...,a_{k\l_k}\}$.
  Fix
$w$ in $\Ga_\l\cap\Ga^{-1}_\lambda$. We write
$w(a_{ij})=b_{ij}+l_{ij}n$, where $1\leq b_{ij}\leq n$ and
$l_{ij}\in\bbZ$.

\bigskip

\noindent{\bf Lemma 5.4.1.} {\sl Assume that  $1\le j< k\le
\l_{i}$. If both $b_{ij}$ and $b_{ik}$ are contained in the same
$\L_q$ for some $q$, then $b_{ij}>b_{ik}$.}

\bigskip

\itp  Note that $l_{ij}\ge l_{ik}$ since  $w(a_{ij})>w(b_{ik})$.
Thus we have
$$w^{-1}(b_{ij})=a_{ij}-l_{ij}n<w^{-1}(b_{ik})=a_{ik}-l_{ik}n.$$Since
$w^{-1}\in\Ga_\l$ and both $b_{ij}$ and $b_{ik}$ are contained in
$\L_q$, we see $b_{ij}>b_{ik}$.

\bigskip

\noindent{\bf Lemma 5.4.2.} {\sl Assume that  $1\le a,b\le n$ and
$w(a)=a_{ij}+l_an$ and $w(b)=a_{ik}+l_bn$, where $1\le j<k\le
\l_{i+1}$. Then $l_a\le l_b$. Moreover, if $l_a=l_b$ then $a>b$.}

\bigskip

\itp Since $w^{-1}(a_{ij})>w^{-1}(a_{ik})$ we see $l_a\le l_b$.
If $l_a=l_b$, from $w^{-1}(a_{ij})>w^{-1}(a_{ik})$ we get $a>b$.

\bigskip

\noindent {\bf Lemma 5.4.3.} (a) $w(\l_1)=e_i+1\text{(mod }n)$
{\sl for some } $i$.

\smallskip

\no(b) $w(e_i+1)\equiv\l_1\text{(mod }n)$ {\sl for some} $i$.

\no(c) $w(1)\equiv e_i\text{(mod }n)$ {\sl for some} $i$.

\no(d) $w(e_i)\equiv 1\text{(mod }n)$ {\sl for some} $i$.

\no(e) {\sl If $w(a_{i1})$ is maximal among all $w(a_{jk})$, then}
$w(a_{i1})\equiv e_j\text{(mod }n)$ {\sl for some $j$.}

\bigskip

\itp (a) Assume that  $w(\l_1)=e_i+j+l_{\l_1}n,$
where $1\leq j\le \l_{i+1}$ and $l_{\l_1}\in\bbZ$. If $j\ne 1$, we can
find $1\le k\le n$ such that $w(k)=e_i+1+l_kn$. Then $l_k> l_{\l_1}$
since $w(k)>w(\l_1)$, see Lemma 5.1.1 (b). Thus
$w^{-1}(e_i+1)<w^{-1}(e_i+j)$. This contradicts $w^{-1}\in
\Ga_\l$. So we must have $j=1$.

\medskip

(b) Applying (a) to $w^{-1}$ we see that (b) is true.

\medskip

(c) Assume that $w(1)=e_{i-1}+j+l_1n,$ where $1\le j\le \l_i$ and
$l_1\in\bbZ$. If $j\ne\l_i$, we can find $1\le k\le n$ such that
$w(k)=e_i+l_kn$. Since $w^{-1}\in\Ga_\l$, we have $l_1<l_k$. Then
$w(k)>w(1+n)>\cdots> w(\l_1+n)$. This contradicts $\l(w)=\l$. So
$j=\l_i$.

\medskip

(d) It follows from $w^{-1}\in\Ga_\l\cap\Ga^{-1}_\l$ and (c).

\medskip

The proof of (e) is similar to that of (a).

\bigskip

\no{\bf Lemma 5.4.4.} {\sl Let $w\in\Ga_\l\cap\Ga_\l^{-1}$ be such
that all components of $\vare(w)$ are non-negative. If the
$(k,l)$-component of $\vare(w)$ is 0 whenever $k\le i-1$, then
$w(a_{\alpha\ga})>w(a_{\beta\ga})$ whenever $r_i\ge
\alpha>\beta\ge 1$ and $\l_\alpha,\l_\beta\ge \ga\ge 1$.}

\bigskip

\itp It follows from the construction of elements in
$\Ga_\l\cap\Ga_\l^{-1}$, see section 5.3. To see it more clearly
we
 use induction on the sum of all components of $\vare(w)$ as in
section 5.3. When the sum of all components of $\vare(w)$ is 0, we
have $w=w_\l$. In this case the lemma is true. Now assume that
$$\vare(w)=(0,...,0,\vare_{i1}(w),...,\vare_{in_i}(w),...,\vare_{j1}(w),...,
\vare_{jh}(w),0,...,0)\in\text{Dom}(F_\l),$$ where $i\le j$, $1\le
h\le n_j$ and $\vare_{jh}(w)\ge 1$.

Let $u$ be in $\Ga_\l\cap\Ga_\l^{-1}$ such that
$\vare_{\alpha\beta}(w)=\vare_{\alpha\beta}(w)$ whenever
$(\alpha,\beta)\ne(j,h)$ and $\vare_{jh}(u)=\vare_{jh}(w)-1$. By
induction hypothesis, we have

\medskip

\no($*)$\qquad\qquad\qquad$u(a_{\alpha\ga})>u(a_{\beta\ga})$

\medskip

\no whenever $r_i\ge \alpha>\beta\ge
1$ and $\l_\alpha,\l_\beta\ge \ga\ge 1$.

\medskip

According to Lemma 5.3.1 and Theorem 5.2.6, we have $$w (a)=\cases
u(a_{k,m_k}) & \text{if } a=a_{k-1,m_{k-1}},\ 2\le k\le r_j,\\
u(a_{1,m_1})+n& \text{if } a=a_{r_j,h},\\ w(a)& \text{if }
a\not\equiv a_{km_k}\text{(mod }n),\ 1\le k\le r_j,\endcases $$
where $m_{r_j}=h$ and $m_{r_j-1},...,m_2,m_1$ are inductively
defined by $$w(a_{k,m_{k}-1})>w(a_{k+1,m_{k+1}})>w(a_{k ,m_{k
}})$$ for $k=r_j-1,...,2,1$, (we set $w(a_{k,0})=\infty$). Using
$(*)$ we see that

\medskip

\no$(\star)$ $m_{r_j}\ge m_{r_j-1}\ge\cdots\ge m_2\ge m_1.$

\medskip

Let $r_i\ge \alpha>\beta\ge 1$ and $\l_\alpha,\l_\beta\ge \ga\ge 1$.

If both $(\alpha,\ga)$ and $(\beta,\ga)$ are different from any
($k,m_k)$ for $k=1,2,...,r_j$, we have
$$w(a_{\alpha\ga})=u(a_{\alpha\ga})>u(a_{\beta\ga})=w(a_{\beta\ga}).$$

If $(\alpha,\ga)=(q,m_q)$ for some $1\le q\le r_j$ and
$(\beta,\ga)$ is different from any ($k,m_k)$ for $k=1,2,...,r_j$,
noting that $u(a_{1m_1})$,...,$u(a_{r_jm_{r_j}})$ form an
r-antichain of $u$ (see the proof of Lemma 5.3.1), we have
$$w(a_{\alpha\ga})>u(a_{\alpha\ga})>u(a_{\beta\ga})=w(a_{\beta\ga}).$$

\def\a{\alpha}
\def\b{\beta}

Suppose that $(\alpha,\ga)$ is different from any ($k,m_k)$ for
$k=1,2,...,r_i$ and $(\beta,\ga)=(q,m_q)$ for some $1\le q\le r_i$. By
$(\star)$ we have $m_\alpha\ge\cdots\ge m_{q+1}\ge m_q=\gamma$. Now
$\ga\ne m_\alpha$, so $\gamma<m_{\alpha}$. When $\alpha=\beta+1$, we
then have
$$w(a_{\alpha\ga})=u(a_{\alpha\ga})>u(a_{\alpha,\ga+1})\ge
u(a_{\alpha m_{\alpha}}) =w(a_{\beta\ga}).$$ When $\a>\b+1$, we have
$$w(a_{\alpha\ga})=u(a_{\alpha\ga})>u(a_{\alpha,m_\a})>\cdots>
u(a_{\b+1, m_{\b+1}}) =w(a_{\beta\ga}).$$

Finally suppose that $(\a,\ga)=(\a,m_\a)$ and $(\b,\ga)=(\b,m_\b)$. By
the relation between $w$ and $u$ we have
$$w(a_{\alpha\ga})>w(a_{\beta\ga}).$$

\medskip

By induction we see that the lemma is true.

\medskip

\section{Some elements of $\Ga_\lambda\cap\Ga_\lambda^{-1}$}

In this section we write explicitly the elements in
$\Ga_\l\cap\Ga_\l^{-1}$ corresponding to fundamental weights of
$F_\l$.

\def\t{\theta}

Let $\theta_{ij}$ ($1\le i\le p,\ 1\le j\le n_i$) be in Dom$(F_\l)$
such that its $(k,l)$-component is 0 whenever $k\ne i$ or $l>j$ and
its $(i,l)$-component is 1 for $l=1,...,j$. Then $\theta_{ij}$ is the
$(i,j)$th fundamental weight of $F_\l$. Given $1\le i_1,...,i_k\le n$,
let $s(i_1,...,i_k) $ be the element of $W$ defined by $i_l\to
i_{l+1}$ for $l=1,...,k-1$, $i_l\to i_1$, and $m\to m$ if $m\not\equiv
i_l$(mod$n$) for all $l$.

\medskip

Let $1\le i\le p$ and $1\le j\le n_i$. For simplicity we write $h$ for
$ {r_i}$. Define
$$u_{ij}=\tau_{\l_1}s(e_1,...,e_h)\tau_{\l_1-1}s(e_1-1,...,e_h-1)
\cdots\tau_{\l_1-j+1}s(e_1-j+1,...,e_h-j+1).$$ Then
$$u_{ij}(a)=\cases e_{k+1}-l+1,&\text{ if } a=e_k-l+1,\ \ 1\le k\le
h-1,\ \ 1\le l\le j;\\ e_{1}-l+1+n,&\text{ if } a=e_h-l+1,\ \ 1\le
l\le j;\\ a ,&\text{ if } a\not= e_k-l+1,\ \ \text{for all }1\le k\le
h,\
\ 1\le l\le j,\endcases$$
where $1\le a\le n$. Note that $w_\l(e_{k-1}+l)=e_k-l+1$ for all $1\le
k\le r$ and $1\le l\le \l_k$. Thus we have
$$u_{ij}w_\l(a)=\cases e_{k+1}-l+1,&\text{ if } a=e_{k-1}+l,\ \
1\le k\le h-1,\ \ 1\le l\le j;\\ e_{1}-l+1+n,&\text{ if }
a=e_{h-1}+l,\
\ 1\le l\le j;\\ w_\l(a) ,&\text{ if } a\not= e_{k-1}+l,\ \ \text{for all
}1\le k\le h,\
\ 1\le l\le j,\endcases$$
where $1\le a\le n$.

\medskip

Recall that $a_{kl}=e_{k-1}+l$ for $1\le k\le r$ and $1\le
l\le\l_k$. According to Lemma 5.3.1 we see easily

\bigskip

\no(a) $u_{ij}w_\l$ is in $\Ga_\l\cap\Ga_\l^{-1}$.

\medskip

\no(b) $\vare(u_{ij}w_\l)=\t_{ij}$.

\bigskip

Using 2.1.3 (e) and noting that $u_{ij}s_k\ge u_{ij}$ if $e_{m-1}+1
\le
k\le e_{m}-1$ for some $1\le m\le r$ (cf. 2.1.3 (f)), we get

\medskip

\no(c) $l(u_{ij})=(n-r_ij)j$ and $l(u_{ij}w_\l)=l(u_{ij})+l(w_\l)$.

\bigskip

Let $j\t_{i1}\in$Dom$(F_\l)$ be such that its $(i,1)$-component is
$j$ and other components are 0. Using Lemma 5.3.1 we see

\bigskip

\no(d) $u_{i1}^jw_\l$ is in $\Ga_\l\cap\Ga_\l^{-1}$.

\medskip

\no(e) $\vare(u_{i1}^jw_\l)=j\t_{i1}$.

\bigskip

Using 2.1.3 (e) we see that
$$l(u_{i1}^{jr_i})=(n-r_i)jr_i=jr_il(u_{i1}).$$ Obviously we have
$u_{i1}^js_k\ge u_{i1}^j$ if $e_{m-1}+1
\le
k\le e_{m}-1$ for some $1\le m\le r$ (cf. 2.1.3 (f)). Therefore we
have

\medskip

\no(f) $l(u_{i1}^j)=jl(u_{i1})$ and $l(u_{i1}^jw_\l)=jl(u_{i1})+l(w_\l)$.

\bigskip

Using the reduced expressions in 2.1.3 (d) we get

\medskip

\no(g) A reduced expression of $u_{i1}$ is
$$\begin{array}
{rl} &\om s_{e_1-2}s_{e_1-3}\cdots s_1s_0s_{n-1}s_{n-2}\cdots
s_{e_h}\\
&\times \hat s_{e_h-1}s_{e_h-2}\cdots s_{e_{ h-1}}\hat
s_{e_{h-1}-1}s_{e_{h-1}-2}\cdots s_{e_{h-2}}\\
&\times\cdots\hat s_{e_2-1}s_{e_2-2}\cdots s_{e_{1}}
\end{array}$$ where $\hat s_k$ means that
$s_k$ is omitted.

\medskip

\chapter{A Factorization Formula in  $J_{\Ga_\lambda\cap\Ga_\lambda^{-1}}$}


\def\tt{\tilde T}
In this chapter we will establish a factorization formula in
$J_{\Ga_\l\cap\Ga_\l^{-1}}$, see Theorem 6.4.1. This formula is
important for our proof of Lusztig Conjecture on based ring for
type $\tilde A_{n-1}$ and is obviously motivated by the
corresponding formula in $R_{F_\l}$, which says that each
irreducible rational $F_\l$-module is a tensor product of some
irreducible modules of reductive components of $F_\l$. In section
6.1 we note that in some cases the integers $\ga_{u,v,w}$ can be
computed through the basis $\{\tilde T_x\ |\ x\in W\}$ instead of
the basis $\{C_x\ |\ x\in W\}$. In section 6.2 we compute the
product $\tilde T_u\tt_v$ for some special $u,v$. This computation
is a key to our factorization formula. In section 6.3 we give some
consequences of the computation in section 2. In section 6.4 we
prove the factorization formula.

\medskip

\section{The integers $\ga_{u,v,w}$}

\def\L{\Lambda}
\def\gg{\Ga_\l\cap\Ga_\l^{-1}}
\def\jgg{J_{\gg}}
\def\cwl{C_{w_\l}}
\def\tt{\tilde T}

In this section we show that in some cases the integer $\ga_{u,v,w}$
can be computed through the product $\tt_u\tt_v$ instead of $C_uC_v$,
the latter is usually much more difficult to compute. Let $\l$ and
$w_\l$ be as in section 2.2. Using induction on $l(w)$ we see

\smallskip

\no(a) If $w=uw_\l$ and $l(w)=l(u)+l(w_\l)$, then $C_w=\phi C_{w_\l}$ for
some $\phi\in H$ with $\bar \phi=\phi$.

\smallskip

\no(b) If $w=w_\l u$ and $l(w)=l(u)+l(w_\l)$, then $C_w= C_{w_\l}\phi$
for some $\phi\in H$ with $\bar \phi=\phi$.

\def\ca{\mathcal A}

\bigskip

Noting that if $u,v$ are in $\gg$, we can find $u_1,v_1$ in $W$ such
that $u=u_1w_\l$, $v=w_\l v_1$ and $l(u)=l(u_1)+l(w_\l)$,
$l(v)=l(v_1)+l(w_\l)$. Thus we can find bar invariant elements
$\phi_1,\phi_2\in H$ such that $$C_u=\phi_1C_{w_\l}\text{ and }\
\ C_v=C_{w_\l}\phi_2.$$ Note that
$$C_{w_\l}\cwl=q^{-l(w_\l)}\sum_{w\le
w_\l}q^{2l(w)}\cwl.$$ Let $$\phi_1\cwl\phi_2=\sum_{w\in W}\eta_wC_w,\
\eta_w\in
\mathcal A.$$ Then $\eta_w$ is bar invariant, i.e. $\bar\eta_w=\eta_w$.
 We have
$$h_{u,v,w}=q^{-l(w_\l)}\sum_{y\le w_\l}q^{2l(y)}\eta_w$$ and then
the degree of $h_{u,v,w}$ is either bigger than or equal to
$a(w_\l)$ whenever $\eta_w\ne 0$. If $w$ is in $\gg$, we must have
$\eta_w\in \bbN$ and $\eta_w=\ga_{u,v,w}$.

\bigskip

\no {\bf Lemma 6.1.1.} {\sl Let $u,v\in\gg$. Write (see 1.2 for the definition of $\tt_w$)
$$\tt_u\tt_v=\sum_{w\in
W}f_{u,v,w}\tt_w,\ \ f_{u,v,w}\in \ca.$$ If deg$f_{u,v,w}\le a(w_\l)=a
$ for all $w$, then $$f_{u,v,w}=\ga_{u,v,w}q^{a}+\text{lower degree
terms,}$$ for any $w\in \gg$.}

\bigskip

\itp Write $$\tt_u\tt_v=\sum_{w\in W}h'_{u,v,w}C_w,\ h'_{u,v,w}\in
\ca.$$

\no Note that
$$\tt_x\in C_x+q^{-1}\sum_{y\in W}\bbZ[q^{-1}] C_y\quad\text{ for any $x\in W$.}$$
 We see
deg$h'_{u,v,w}\le a(w_\l)=a$ for any $w$ since deg$f_{u,v,w}\le a $.
Moreover deg$h'_{u,v,w} =a$ if and only if deg$f_{u,v,w}=a $, and in
this case the leading coefficients of $h'_{u,v,w}$ and $f_{u,v,w}$
coincide.

For $w\in\gg$ we have $$h'_{u,v,w}=\ga_{u,v,w}q^a+\text{ lower degree
terms.}$$ Therefore $$f_{u,v,w}=\ga_{u,v,w}q^{a}+\text {lower degree
terms,}$$ for any $w\in \gg$.

The lemma is proved.

\bigskip

\def\mb{\mathbf}

\no{\bf Lemma 6.1.2.} {\sl Let $u,v\in W$ and $s_{i_1}\cdots s_{i_k}$ a
reduced expression of $v$. Then
$$\tt_u\tt_v=\sum_{\mathbf j}g_{u,v,\mathbf j}\tt_{uw_{\mb j}},\ \
g_{u,v,\mb j}\in\ca,$$ where $\mb j$ runs through all subsequences
$j_1,...,j_m$ (including empty subsequence) of $i_1,...,i_k$ and
$w_{\mb j}=s_{j_1}\cdots s_{j_m}$.}

\bigskip

\itp It follows from the following fact
$$\tt_{x}\tt_{s_i}=\cases (q-q^{-1})\tt_x+\tt_{xs_i}&\text{ if
$xs_i\le x$}\\ \tt_{xs_i}&\text{ if $xs_i\ge x.$}\endcases$$

\bigskip

Obviously when $i_1,...,i_k$ are pairwise different, then $w_{\mb
j}\ne w_{\mb j'}$ if $\mb j\ne \mb j'$, and $g_{u,v,\mb j}$ is either
0 or a power of $q-q^{-1}$.

\bigskip

\no{\bf Lemma 6.1.3.} {\sl Let $x\in W$ and $j,k$ be integers with
 $0\le k-j\le n-2$. Assume that
$xs_l\le x$ for $l=j,j+1,...,k,$ and let $y=s_ks_{k-1}\cdots s_j$.
Then $$\tt_x\tt_y=\sum_{k\ge k_1>k_2>\cdots>k_m\ge j}\xi^{k-j+1-m}
\tt_{xk_1\cdots k_m},$$ here we use $xk_1\cdots k_m$ for
$xs_{k_1}\cdots s_{k_m}$, and $\xi$ stands for $q-q^{-1}$.}

\bigskip

{\it Proof.} We use induction on $k-j$. When $k-j=0$, we have
$$\tt_x\tt_y=\xi\tt_x+\tt_{xs_k}.$$
The lemma is true in this case. Now suppose that the lemma is true
when $y$ is replaced by $ys_j$. Then we have
$$\tt_x\tt_{ys_j}=\sum_{k\ge k_1>k_2>\cdots>k_m\ge j+1}
\xi^{k-j-m}\tt_{xk_1\cdots k_m}.$$
Note that $xs_{k_1}\cdots s_{k_m}s_j\le xs_{k_1}\cdots s_{k_m}$ for
any sequence $k\ge k_1>k_2>\cdots>k_m\ge j+1$. Thus we have
$$\begin{array}{rl}  \tt_x\tt_y&=\tt_x\tt_{ys_j}\tt_{s_j}
   \\[4mm]
&=\sum_{k\ge k_1>k_2>\cdots>k_m\ge j+1}\xi^{k-j-m}
\tt_{xk_1\cdots k_m}\tt_{s_j}
\\[4mm]
&=\sum_{k\ge k_1>k_2>\cdots>k_m\ge j+1}\xi^{k-j-m}(\xi \tt_{xk_1\cdots k_m}+\tt_{xk_1\cdots k_m j})
  \\[4mm]
&=\sum_{k\ge k_1>k_2>\cdots>k_m\ge j}\xi^{k-j+1-m}\tt_{xk_1
\cdots k_m}.
\end{array}$$
The lemma is proved.

\bigskip

Let $W_\l$ be the subgroup of $W$ generated by all simple reflections
that appear in a reduced expression of $w_\l$.

\bigskip

\no{\bf Lemma 6.1.4.} {\sl Let $u\in W$. Write $u=u_1u_2$ such that
$u_1$ is the shortest element in the coset $uW_\l$ and $u_2$ is in
$W_\l$. Then in
$$\tt_u\tt_{w_\l}=\sum_{w\in W}f_{u,w_\l,w}\tt_w,\ \
f_{u,w_\l,w}\in\A,$$ we have deg$f_{u,w_\l,w}<l(u_2)$ if $w\ne
u_1w_\l$ and $f_{u,w_\l,u_1w_\l}=q^{l(u_2)}+$ lower degree terms.}

\bigskip

\itp We have $l(u_1)+l(u_2)=l(u)$, $l(u_1w_\l)=l(u_1)+l(w_\l)$ and
$l(u_2w_\l)=l(w_\l)-l(u_2)$. Using induction on $l(u_2)$ we can prove
the lemma.

\bigskip

 In next section we compute some $\tt_u\tt_v$.

\medskip

\section{A computation for some $\tt_u\tt_v$}

In this section we compute $\tt_u\tt_v$ for some special $u,v$ in
$\gg$. Suppose that $w\in\gg$ and
$\vare(w)=(\vare_{11}(w),...,\vare_{pn_p}(w))$. Denote by $\vare_i(w)$
the element in Dom($F_\l)$ whose $(i,j)$-component is $\vare_{ij}(w)$
for $j=1,2,..., n_i$ and other components are 0. Then
$\vare(w)=\vare_1(w)+\cdots+\vare_p(w)$. The main result of this
section is

\bigskip

\no{\bf Proposition 6.2.1.} {\sl Let $u,v\in \gg$ be such that all
components of $\vare(u)$ and $\vare(v)$ are nonnegative. Assume
that $\vare_1(u)=\vare_2(u)=\cdots=\vare_{i-1}(u)=0$ and
$\vare_{i1}(v)=1$ and $\vare_{kl}(v)=0$ if $(k,l)\ne(i,1)$. Then
deg$f_{u,v,w}\le a(w_\l)$ for all $w\in W$, see Lemma 6.1.1 for
the definition of $f_{u,v,w}$. }

\bigskip

\no{\bf 6.2.2.} {\it Proof of Proposition 6.2.1.} When $r=1$,
the proposition is proved in [L1]. Now we suppose that $r>1$. Recall
that $e_0=0,\ e_k=\l_1+\l_2+\cdots+\l_k$, $a_{kl}=e_{k-1}+l$ for
$k=1,2,...,r$, $1\le l\le \l_{k}$. Set $h=r_i$, see section 5.1 for
definition of $r_i$. Let
$$\begin{array}
{rl} v_1=& s_{e_1-1}s_{e_1-2}\cdots s_2s_1\om s_{n-1}s_{n-2}\cdots
s_{e_h}\\
&\times \hat s_{e_h-1}s_{e_h-2}\cdots s_{e_{ h-1}}\hat
s_{e_{h-1}-1}s_{e_{h-1}-2}\cdots s_{e_{h-2}}\\
&\times\cdots\hat s_{e_2-1}s_{e_2-2}\cdots s_{e_{1}}
\end{array}$$ where $\hat s_k$ means
$s_k$ is omitted. According to 5.5 (a)-(c) and 5.5 (g) we have
$v=v_1w_\l$ and $l(v)=l(v_1)+l(w_\l)$.

\medskip

Now we compute $\tt_u\tt_{v_1}$. To avoid complicated subscripts we
also use $\tt(w)$ for $\tt_w$ for any $w\in W$. Set
$$\begin{array} {rl} v_j=
s_{e_j-1}s_{e_j-2}\cdots s_{e_h}&\hat s_{e_h-1} s_{e_h-2}
\cdots s_{e_{h-1}}\\
 \times &\hat s_{e_{h-1}-1}s_{e_{h-1}-2}
\cdots s_{e_{h-2}}\\
 &\cdots\\
\times&\hat s_{e_2-1}s_{e_2-2}
\cdots s_{e_1},
\end{array}$$
for $j=r,r-1,..., h+1$, note that $r_i=h$;
$$\begin{array} {rl} v_{j}=
&\hat s_{e_j-1}s_{e_j-2}
\cdots s_{e_{j-1}}\\
 \times&\hat s_{e_{j-1}-1}s_{e_{j-1}-2}
\cdots s_{e_{j-2}}\\
 &\cdots\\
\times&\hat s_{e_2-1}s_{e_2-2}
\cdots s_{e_1}.
\end{array}$$
for $j=2,3,..., h$.

\medskip

For simplicity we often use $i_1i_2\cdots i_k$ for
$s_{i_1}s_{i_2}\cdots s_{i_k}$. Since $us_q\le u$ for all $1\le q\le
e_1-1=\l_1-1$, according to Lemma 6.1.3 we have
$$\tt_u\tt_v=\sum_{e_1-1\ge j_{1,k_1}>j_{1,k_1-1}>
\cdots>j_{1,1}\ge1}\xi^{\l_1-1-k_1}\tt(uj_{1,k_1}j_{1,k_1-1}\cdots
j_{1,1}\om)\tt( v_r),$$ where $\xi$ stands for $q-q^{-1}$.

\medskip

Let $$u_1=uj_{1,k_1}j_{1,k_1-1}\cdots j_{1,1}\om.$$ Recall that
$a_{jk}=e_{j-1}+k$ for all $1\le j\le r$ and $1\le k\le
\l_j$. We have

\medskip

\no(a) $u_1(n)=u(a_{1l_1})+n $ for some $1\le l_1\le \l_1$.

\medskip

It is easy to see

\medskip

\no (b) $l_1\le k_1+1$. Moreover $l_1=k_1+1$ if and only if
$j_{1,q}=q$ for $q=1,...,k_1$.

\medskip

\no(c) $u_1(a)=u(a+1)$ for all $e_1\le a\le n-1$.

\medskip

According to Lemma 5.1.1 (c) we have $$u_1
(n)=u(a_{1l_1})+n>u(a_{ml_1})=u_1(a_{ml_1}-1)$$ for $m=2,3,...,r$. For
convenience we set $$u(a_{m0})=\infty\text{\ \ and \ \
}u(a_{m,\l_m+1})=-\infty$$ for all $m$. Choose $0\le m_r\le \l_r$ such
that
$$u(a_{rm_r})>u(a_{1l_1})+n>u(a_{r,m_r+1}).$$
Then $m_r\le l_1-1$. Set $$u'_r=u_1s_{e_{r}-1}s_{e_r-2}\cdots
s_{e_{r-1}+m_r}.$$ Then $l(u'_r)=l(u_1)+\l_r-m_r$. Moreover we have
$u'_rs_q\le u'_r$ for
$q=e_{r-1},e_{r-1}+1,...,e_{r-1}+m_r-1=a_{rm_r}-1$. Thus, using Lemma
6.1.3 we get
$$\begin{array} {rl}
\tt_{u_1}\tt_{
v_{r}}&=\tt(u'_r )
\tt(s_{e_{r-1}+m_r-1}\cdots s_{e_{r-1}})\tt(v_{r-1})\\[3mm]
&=\displaystyle{\sum\xi^{m_r-k_r}\tt(u'_r j_{r,k_r}
j_{r, k_r-1}\cdots j_{r,1} )\tt( v_{r-1})},
\end{array}$$
where the sum is for all sequences $e_{r-1}+m_r-1\ge
j_{r,k_r}>j_{r,k_r-1}>$ $ \cdots>j_{r,1}\ge e_{r-1}$.

\smallskip

Set $$u_{r}=u'_r j_{r,k_r} j_{r,k_r-1}\cdots j_{r,1}.$$ We have

\medskip

\no(d) If $m_r=k_r$ then $u_{r}(e_{r-1})=u(a_{1l_1})+n$, in this case
we set $l_r=m_r+1$. If $0\le k_r<m_r$ then $u_r(e_{r-1})=u(a_{rl_r})$
for some $1\le l_r\le k_r+1$. Note that in any case the $m_r-l_r+1$
elements $$u(a_{rm_r}),\ u(a_{r,m_r-1}),\
...,\ u(a_{r,l_r+1}),\ u_r(e_{r-1})
$$
  are bigger than $u(a_{1l_1})+n$, and the $m_r-l_r+1$ elements
$$u(a_{rm_r}),\ u(a_{r,m_r-1}),\ ...,\ u(a_{r,l_r+1}),\
u_1(e_r) $$ are less than $u_r(e_{r-1})$.

\medskip
\no(e) $u_r(a)=u(a+1)$ for all $e_1\le a\le e_{r-1}-1.$

\medskip

Suppose that we have defined $u'_{{c'} }, u_{{c'} }, m_{{c'} }\ge
k_{{c'}}\ge 0, l_{{c'}}$ for all $r\ge {c'}>c\ge r_i+1=h+1$ such that

\medskip

\no(f) $u_{c'}(e_{c'-1})=u(a_{{c'},l_{{c'}}})$ for some $1\le
l_{c'}\le k_{c'}+1$ if $k_{c'}<m_{c'}$; and $u_{c'}(
e_{{c'}-1})=u_{{c'}+1}(e_{c'})$ if $m_{c'}=k_{c'}$, in this case we
set $l_{c'}=m_{c'}+1$, (we understand that $u_{r+1}=u_1$;)

\medskip

\no(g) $u_{c'}(a)=u(a+1)$ for all $e_1\le a\le e_{c'-1}-1$.

\bigskip
  Now we define $u_c,u'_c,k_c,l_c,m_c$ as follows.

\bigskip

Choose $0\le m_c\le \l_c$ such that $$u(a_{cm_c})>u_{c+1}(e_c)
>u(a_{c,m_c+1})=u_{c+1}(a_{c,m_c}).$$ Set
$$u'_c=u_{c+1} s_{e_c-1}s_{e_c-2}\cdots s_{e_{c-1}+m_c}.$$ Then
$l(u'_c)=l(u_{c+1})+\l_c-m_c$ and $u'_cs_q\le u'_c$ for
$q=e_{c-1},e_{c-1}+1,...,e_{c-1}+m_c-1$. Using Lemma 6.1.3 we see
$$\begin{array} {rl}
\tt_{u_{c+1}}\tt_{
v_{c}}&=\tt(u'_{c} )
\tt(s_{e_{c-1}+m_c-1}\cdots s_{e_{c-1}})\tt(v_{c-1})\\[3mm]
&=\displaystyle{\sum\xi^{m_c-k_c}\tt(u'_{c} j_{c,k_c}
j_{c,k_c-1}\cdots j_{c,1})\tt( v_{c-1}),}
\end{array}$$
where the sum is for all sequences $e_{c-1}+m_c-1\ge
j_{c,k_c}>j_{c,k_c-1}> \cdots$ $>j_{c,1}\ge  e_{c-1}.$ \no Set
$$u_c=u'_{c} j_{c,k_c} j_{c,k_c-1 }\cdots j_{c,1 }.$$

\bigskip

If $k_c<m_c$ we can find $1\le l_c\le k_c+1$ such that
$u_c(e_{c-1})=u(a_{cl_c})$; if $m_c=k_c$ we set $l_c=m_c+1$, and in
this case we have $u_c(e_{c-1})=u_{c+1}(e_c)$. Obviously we have

\bigskip

\no(h) The  $m_c-l_c+1$ elements
$$u(a_{cm_c}),\ u(a_{c,m_c-1}),\ ...,\ u(a_{c,l_c+1}),\ u_c(e_{c-1})$$ are
bigger than $u_{c+1}(e_{c})$, and the $m_c-l_c+1$ elements
$$u(a_{cm_c}),\ u(a_{c,m_c-1}),\ ...,\ u(a_{c,l_c+1}),\
u_{c+1}(e_{c})$$ are less than $u_{c}(e_{c-1})$.

\medskip

\no(i) $u_c(a)=u(a+1)$ for all $e_1\le a\le e_{c-1}-1$.

\bigskip

In this way we defined $u_c,u'_c,l_c, k_c,m_c $ for all $r\ge c\ge
r_i+1=h+1$. From (i) we see

\bigskip

\no(j)   $u_{h+1}s_j\le u_{h+1}$ if $e_1\le j\le
e_h-2$ and $j\ne e_2-1,e_3-1,...,e_{h-1}-1$.

\bigskip

By Lemma 6.1.3 we have

$$\begin{array} {rl}
\tt_{u_{h+1}}\tt_{
v_{h}}&=\displaystyle{\sum_{\st{\sc h\ge c\ge 2}{e_c -2\ge
j_{c,k_c}>j_{c,k_c-1}>
\cdots>j_{c,1}\ge  e_{c-1}}}\xi^{\sum_{h\ge c\ge 2}(\l_c-1-k_c)}}\\[4mm]
&\qquad \displaystyle{\times \tt(u_{h+1} \prod_{h\ge c\ge 2}(j_{c,k_c}
j_{c,k_c-1}\cdots j_{c,1}) ).}
\end{array}$$

\bigskip

Write $$w'=u_{h+1}\displaystyle{ \prod_{h\ge c\ge 2}j_{c,k_c}
j_{c,k_c-1}\cdots j_{c,1}}.$$ We are concerned with the degree of
$f_{u,v_1,w'}$. From the construction above we see

\bigskip

\no(k) The degree of $f_{u,v_1,w'}$ is $$\l_1-1-k_1+\sum_{r\ge c\ge
h+1}(m_r-k_r)+\sum_{h\ge c\ge 2}(\l_c-1-k_c).$$

\bigskip

We have
$$\begin{array}{rl}
\tilde T_u\tilde T_v & =\tilde T_u\tilde T_{v_1}\tilde T_{w_\l}\\
 & =\sum_{w'\in W}f_{u,v_1,w'}\tilde T_{w'}\tilde T_{w_\l}\\
& =\sum_{w',w\in W}f_{u,v_1,w'}f_{w',w_\l,w}\tilde T_w.
\end{array}$$

\medskip

\no Now we consider the degree of $f_{w',w_\l,w}$. Note that
$w'(e_{c})=u(a_{c+1,l_{c+1}})$ for some $1\le l_{c+1}\le
\l_{c+1}$ whenever $1\le  c\le  h-1$. Obviously we have

\bigskip

\no(l)   $l_{c+1}\le  k_{c+1}+1$ for $1\le  c\le  h-1$.

\bigskip

According to Lemma 5.4.4 we have

\bigskip

\no(m) $u(a_{c+1,l_{c+1}})>u(a_{c,l_{c+1}})$ if $1\le  c\le  h-1$.

\bigskip

Note that if $1\le c\le h-1$ then
$w'(a_{c1}),...,w'(a_{c,\l_c-1})$ are in
$\{u(a_{c1}),...,u(a_{c\l_c})\}$. Thus,

\bigskip

\no(n) among $w'(a_{c1}),...,w'(a_{c,\l_c-1})$ at
least $\l_c-(l_{c+1}-1)-1=\l_c-l_{c+1}$ of them are less than
$w'(e_c)=u(a_{c+1,l_{c+1}})$.

\medskip

\no(o)  Suppose that among $u(a_{h1}),...,u(a_{h\l_h})$
we have that $\beta$ of them are less than
$w'(a_{h\l_h})=u_{h+1}(e_h)$. Since
$w'(a_{h1}),...,w'(a_{h,\l_{h}-1})$ are in $\{u(a_{h1}),...,
u(a_{h\l_{h}})\}$ we see at least $\beta-1 $ of
$w'(a_{h1}),...,w'(a_{h,\l_{h}-1})$ are less than $w'(a_{h\l_h})$.

\bigskip

 Using
(h) and (o) we get an r-chain of $u$ consisting of
$$\begin{array}{rl}
&u(a_{h1}),\ ...,\ u(a_{h,\l_{h}-\beta}),\  u_{h+1}(e_h),\
u(a_{h+1,l_{h+1}+1}),\ ...,\ u(a_{h+1,m_{h+1}}),\\[2mm]
&\ \ ...,\  u_r(e_{r-1}),\
u(a_{r,l_r+1}),\ ...,\ u(a_{rm_r}),\ u(a_{1l_1}+n),\
...,\ u(a_{1\l_1}+n),
\end{array}$$ whose length is
$$L=\l_{h}-\beta+m_{h+1}-l_{h+1}+1+\cdots+m_r-l_r+1+\l_1-l_1+1.$$
If $$-(\beta-1)+m_r-k_r+\cdots+m_{h+1}-k_{h+1}+\l_{h}-1-k_1>0,$$ then
we have $$L\ge
\l_{h}-\beta+m_{h+1}-k_{h+1}+\cdots+m_r-k_r+\l_1-(k_1+1)+1>\l_1.$$
This is impossible since the partition of $u$ is $\l$. So we have

\bigskip

\no(p)\ \  $-(\beta-1)+m_r-k_r+\cdots+m_{h+1}-k_{h+1}+\l_{h}-1-k_1\le 0$

\bigskip

Write $w'=w'_1w'_2$ such that $w'_1$ is the shortest element in
$w'W_\l$ and $w'_2$ is in $W_\l$. According to (n) and (o), we have

\bigskip

\no(q) $l(w'_2)\le l(w_\l)-\sum_{1\le c\le h-1}(\l_c-l_{c+1})-(\beta-1).$

\bigskip

Using Lemma 6.1.4 we see that the degree of $f_{w',w_\l,w}$ is less
than or equal to $l(w_\l)-\sum_{1\le c\le
h-1}(\l_c-l_{c+1})-(\beta-1).$

\bigskip

\no(r) Since the degree of $f_{u,v,w}$ is the maximal number in
 $$\{\text{deg}f_{u,v_1,w'}+\text{deg}f_{w',w_\l,w}\ |\ w'\in W\},$$ we
 see that the degree of
$f_{u,v,w}$ is less than or equal to
$$\begin{array}{rl}\ &\l_1-k_1-1+m_r-k_r+\cdots+m_{h+1}-k_{h+1}\\
\ +&\l_{h}-1-k_{h}+\cdots+\l_2-1-k_2\\
\ +&l(w_\l)-(\l_{h-1}-l_{h})-\cdots-(\l_1-l_2)-(\beta-1)\\
\le&l(w_\l)=a(w_\l).
\end{array}$$

The Proposition is proved.

\bigskip

\no{\bf Corollary 6.2.3.} {\sl Keep the notation in 6.2.2. Suppose that
  $f_{u,v,w}$ has degree $a(w_\l)$. Then we   have}

\smallskip

\no(a) {\sl $l_c$ is equal to $k_c+1$ for $1\le c\le h$.}

\smallskip

\no(b) $u(a_{c-1,l_c-1})>u(a_{cl_c})>u(a_{c-1,l_c})$ {\sl for $2\le
c\le h$. (Note that $h=r_i$.)}

\smallskip

\no(c) $l_c\le l_{c-1}$ {\sl for} $2\le c\le h $.

\smallskip

\no(d)\ \ $-(\beta-1)+m_r-k_r+\cdots+m_{h+1}-k_{h+1}+\l_{h}-1-k_1=0$.

\smallskip

\no(e)\ \ $u_{h+1}(e_h)>u(a_{hl_h})$ and
$u(a_{cl_c})>u(a_{c-1,l_{c-1}})$ for $c=2,3,...,h$.

\bigskip

\itp (a) If $l_c<k_c+1$ for some $c$, then $w'(a)<w'(b)$ for some
$$e_{c-1}+l_c+1\le a<b\le e_{c-1}+k_c+1<e_c.$$ Thus $l(w'_2)$ would be
less than $$l(w_\l)-\sum_{1\le c\le h-1}(\l_c-l_{c+1})-(\beta-1).$$
From the proof in 6.2.2 we see that the degree of $f_{u,v,w}$ then
would be less than $a(w_\l)$. (When $2\le c\le h$ we can see that
$l_c=k_c+1$ using 6.2.2 (l) and 6.2.2 (r).)

\medskip

(b) By Lemma 5.4.4 we have $u(a_{cl_c})>u(a_{c-1,l_c})$. Suppose
that $u(a_{cl_c})>u(a_{c-1,l_c-1})$ for some $c$ with $2\le c\le
h$. Then among $w'(a_{c-1,1}),...,w'(a_{c-1,\l_{c-1}-1})$ at least
$\l_{c-1}-(l_{c}-2)-1=\l_{c-1}-l_{c}+1$ of them are less than
$w'(a_{c-1,l_{c-1}})$. Then $l(w'_2)$ would be less than
$l(w_\l)-\sum_{1\le c\le h-1}(\l_c-l_{c+1})-(\beta-1).$ From 6.2.2
(l) and 6.2.2 (p-r) we see that the degree of $f_{u,v,w}$ then
would be less than $a(w_\l)$. Therefore (b) is true.

\medskip

 (c) The proof is similar. If $l_c>l_{c-1}$
 for some $2\le c\le h$, then
 among $w'(a_{c-1,1}),$...,$w'(a_{c-1,\l_{c-1}-1})$ at least
$\l_{c-1}-(l_{c}-2)-1=\l_{c-1}-l_{c}+1$ of them are less than
$w'(a_{c-1,l_{c-1}})$. As the proof of (b), this is impossible.
Therefore (c) is true.

\medskip

(d) It is clear from the proof of Prop. 6.2.1.

\medskip

(e) When $c=2,3,...,h$, using (c) and Lemma 5.4.4, we see that
$u(a_{cl_c})>u(a_{c-1,l_{c-1}})$. If $u_{h+1}(e_h)<u(a_{hl_h})$,
then $\beta $ of $w'(a_{h1}),...,w'(a_{h,\l_{h}-1})$ are less than
$w'(a_{h\l_h})$. Thus $l(w'_2)$ would be less than
$l(w_\l)-\sum_{1\le c\le h-1}(\l_c-l_{c+1})-(\beta-1).$ By 6.2.2
(q-r), this contradicts that $f_{u,v,w}$ has degree $a(w_\l)$.
Therefore $u_{h+1}(e_h)>u(a_{hl_h})$. (e) is proved.

\medskip

The proof is completed.

\bigskip

\no{\bf Corollary 6.2.4.} {\sl Keep the notation in 6.2.2. Suppose that
deg$f_{u,v,w}=a(w_\l)$. }

\smallskip

\no(a) {\sl We have $w(a_{cl})>w(a_{cl'})$ for any $1\le c\le r$ and $1\le
l<l'\le\l_c$.}

\smallskip

\no(b) {\sl Suppose that $1\le c<h$. Then
$$w(a_{cl})=\cases u(a_{cl}) & \text{ if } 1\le l\le l_{c+1}-1,\\
u(a_{c+1,l_{c+1}}) & \text{ if } l=l_{c+1},\\
 u(a_{c,l-1}) & \text{ if
} l_{c+1}<l\le l_c\\
u(a_{cl}) & \text{ if } l_c <l\le
\l_c.\endcases$$}

\smallskip

\no(c) {\sl $$w(a_{hl})=\cases u(a_{hl}) & \text{ if }
 1\le l\le \l_h-\beta, \\
 w'(e_h) & \text{ if } l=\l_h-\beta+1,\\
u(a_{h,l-1}) & \text{ if } \l_h-\beta+1<l\le l_h,\\
 u(a_{hl}) & \text{ if } l_h<l\le\l_h.\endcases$$}

\smallskip

\no(d) {\sl Suppose that $h+1\le c\le r$. We have
$$w(a_{cl})=\cases u(a_{cl}) & \text{ if } 1\le l<l_c ,\\
u(a_{c,l+1}) & \text{ if } l_c\le l< m_c,\\
  u_{c+1}(e_c) & \text{ if } l=m_c,\\
  u(a_{cl}) & \text{ if } m_c<l\le\l_c . \endcases$$}

\bigskip

\itp (a) From the proof of Prop. 6.2.1 we see
$l(w)=l(ww_\l)+l(w_\l)$. Using Lemma 2.5.1 we get (a).

\medskip

(b) Using Corollary 6.2.3 (b)-(c) we see that (b) is true.

\medskip

(c) and (d) are clear from the proof of Prop. 6.2.1.

\medskip

The corollary is proved.

\medskip

\section{Some consequences}

Let $u,v$ be as in section 6.1. In this section we will figure out the
$w$ in $\gg$ with $f_{u,v,w}$ having degree $a(w_\l)$. We keep the
notation in 6.2.2.

\bigskip

\no{\bf Lemma 6.3.1.} {\sl If $\beta=0$ then deg$f_{u,v,w}$ is less
than
$a(w_\l)$.}

\bigskip

\itp In 6.2.2 we have showed  that
$$-(\beta-1)+m_r-k_r+\cdots+m_{h+1}-k_{h+1}+\l_{h}-1-k_1\le0.$$
When $\beta=0$, we have
$$m_r-k_r+\cdots+m_{h+1}-k_{h+1}+\l_{h}-k_1\le0.$$ Therefore the
degree of $f_{u,v,w}$ is either less than or equal to
$$\begin{array}{rl}\
&\l_1-k_1-1+m_r-k_r+\cdots+m_{h+1}-k_{h+1}\\[3mm] \
&+\l_{h}-1-k_{h}+\cdots+\l_2-1-k_2\\[3mm] \
&+l(w_\l)-(\l_{h-1}-l_{h})-\cdots-(\l_1-l_2) \\[3mm]
<&l(w_\l)=a(w_\l).
\end{array}$$

The lemma is proved.

\bigskip

\no{\bf Lemma 6.3.2.} {\sl Assume that deg$f_{u,v,w}=a(w_\l)$.}

\smallskip

\no(a) {\sl If $l_1>\l_{h}$  then $w$ is not in
$\gg$.}

\smallskip

\no(b) If $w$ is in $\gg$, then $\vare_k(u)$ is equal to $\vare_k(w)$ if $1\le k\le i-1$.

\bigskip

\itp (a) According to Corollary 6.2.3 (d) we must have
$$-(\beta-1)+m_r-k_r+\cdots+m_{h+1}-k_{h+1}+\l_{h}-1-k_1=0.$$ Thus
the r-chain of $w$, (here we need Corollary 6.2.4 and 6.2.2 (h))
$$\begin{array} {rl} &w(a_{h1}),\ \ ...,\ \
w(a_{h,\l_{h}-\beta}),\ \ w(a_{h,\l_h-\beta+1}),\\[2mm] &w(a_{h+1
,l_{h+1}}),\ \ ...,\ \ w(a_{h+1,m_{h+1}}),\\[2mm] &\qquad\quad
...,\\[2mm] &w(a_{r,l_r}),\ \ ...,\ \ w(a_{r,m_r}), \ \
\end{array}$$ has length $$L=\l_{h}-\beta+1+m_{h+1
}-l_{h+1}+1+\cdots+m_r-l_r+1 .$$ We have
$$\begin{array}{rl}
L&\ge\l_{h}-\beta+1+m_{h+1}-k_{h+1}+\cdots+m_r-k_r\\
& =k_1+1\\
& \ge
l_1\\
&>\l_{h}.
\end{array}$$
By Corollary 6.2.4 (a), the sequence
$w(a_{c1})>\cdots>w(a_{c\l_c})$ is an r-chain of $w$ of length
$\l_c$ if $1\le c\le h-1$. Thus $w$ has an r-chain family set of
index $h$ and the cardinality of the r-chain family set is bigger
than $e_h$. Therefore $w\not\in\gg$.

\medskip

(b) Since $\vare_k(u)=0$ for $k=1,...,i-1$ and all components of
$\vare(u)$ are non-negative, using Lemma 5.3.2, Theorem 5.1.12 and
Lemma 5.1.4 we see that $$ 1\le u(a_{pq})\le n\text{\qquad
if\quad} 1\le p\le r_{i-1}\text{ and } q>\l_h.\leqno\hbox{(1)}$$
By (a) and Corollary 6.2.3 (c) we get $$\l_h\ge l_1\ge
l_2\ge\cdots\ge l_h.\leqno\hbox{(2)}$$ Using Corollary 6.2.4 we
then get $$w(a_{pq})=u(a_{pq}) \text{\qquad if\quad} 1\le p\le
r_{i-1}\text{ and } q>\l_h.\leqno\hbox{(3)}$$

Using Lemma 5.4.4 for $u$ and Corollary 6.2.4 we see clearly that
$$w(a_{pq})>w(a_{p-1,q})\text{\qquad if\quad} 2\le p\le h-1\text{
and } 1\le q\le \l_p,\l_{p-1}.\leqno\hbox{(4)}$$

Now we show that $$w(a_{hq})>w(a_{h-1,q})\text{\qquad if\quad}
1\le q\le \l_h,\l_{h-1}.\leqno\hbox{(5)}$$ We must have
$l_h\ge\l_h-\beta+1$ since deg$f_{u,v,w}=a(w_\l)$. Using Lemma
5.4.4 and Corollary 6.2.4 we see that (5) is true.

Using (1), (3-5), Lemma 5.1.4 and Lemma 5.2.3, we see that
$$\vare_{k}(w)=(0,...,0)=\vare_k(u)$$ if $1\le k\le i-1$. (b) is
proved.

\medskip

The lemma is proved.

\bigskip

\no{\bf Definition 6.3.3.} Let $a\in\bbZ$. Write $a=a_{ij}+kn,\
1\le i\le r, \ 1\le j\le\l_i$. The {\bf level} of $a$ is defined
to be the pair $(i,k)$. We say that $(i,k)>(i',k')$ if $k>k'$ or
$i>i'$ and $k=k'$.

\bigskip

Let $w\in W$. An r-antichain of $w$ is called {\bf saturated} if the
r-antichain is contained in some complete r-antichain family of $w$.
Analogously we define {\bf saturated d-antichains} of $w$.

\bigskip

\no{\bf Lemma 6.3.4.} {\sl If $w$ is in $\gg$ and deg$f_{u,v,w}=a(w_\l)$,
then all $$u(a_{cl_c}),\ u(a_{c,l_c+1}),\ ...,\ u(a_{cm_c}),\
u_{c+1}(e_c)$$
  have the same level whenever $h+1\le c\le r$ and $l_c\le m_c$. (Recall that
  $u_{r+1}=u_1.$)}

\bigskip

\itp Suppose that the conclusion is not true. Then there exists
some $c$ with $h+1\le c\le r$ and $l_c\le m_c$ such that some pair
of neighboring terms in  the sequence $$u(a_{cl_c}),\
u(a_{c,l_c+1}),\ ...,\ u(a_{cm_c}),\ u_{c+1}(e_c)$$ have distinct
levels but all $$u(a_{c'l_{c'}}),\ u(a_{c',l_c'+1}),\ ...,\
u(a_{c'm_{c'}}),\ u_{c'+1}(e_{c'})$$ have the same level for all
$c<c'\le r$ with $l_{c'}\le m_{c'}$.

\medskip

Let $Z$ be a complete r-antichain family of $u$. For
$\theta=l_c,l_{c+1},...,m_c$ we choose $\beta_\t$ such that
$u(a_{c-1,\beta_\t})$ and $u(a_{c\t})$ are in an r-antichain of $u$
that is contained in $Z$.

Suppose that all $$u(a_{cl_c}),u(a_{c,l_c+1}),...,u(a_{c,\t-1}) $$
have the same level but $u(a_{c,\t-1})$ and $u(a_{c\t})$ have distinct
levels for some $l_c< \t\le m_c$. Obviously we have

\medskip

\no(a) $u(a_{c,\t-1})> u(a_{c\t}).$

\medskip

Suppose that $c-1>h$. We claim that

\medskip

\no(b) $u(a_{c-1,\beta_{\t-1}})>u(a_{c\t})$.

\medskip

Write $$u(a_{cl_c})=\xi+\eta n,\ \ 1\le \xi\le n, \ \eta\in\bbZ.$$
Then $\xi\in\L_q$ for some $1\le q\le r$. By assumption,
$w(a_{kl})$ and $u(a_{kl})$ have the same level for all $c<k\le r$
and $1\le l\le \l_k$. By Lemma 5.1.4, then  we may assume that the
r-antichain $C$ in $Z$ containing $u(a_{cl_c})$ has length $c$ and
there is a saturated r-antichain of $w$ that contains
$u_{c+1}(e_c)$ and has length $c$.

 Note that $u(a_{c-1,\beta_{l_{c}}})$ is in $C$ and
  $u_{c-1}(e_{c-2})\ge
 u(a_{cl_c})$ (since $c-1>h$). Using Lemma
5.1.4, then we can find a saturated r-antichain of $w$ that
contains $u(a_{c-1,\b_{l_{c}}})$ and has length $c-1$ since
$u(a_{cl_c})$ and $u_{c+1}(e_c)$ have distinct levels.

We claim that $q>1$. Otherwise, $q=1$. Write
$$u(a_{c-1,\b_{l_c}})=x+y n,\ \ \ 1\le x\le n,\ \ y\in\bbZ.$$ By
Lemma 5.1.4, we have $x\in\L_{c}$. But there is a saturated
r-antichain of $w$ that contains $u(a_{c-1,\b_{l_c}})$ and has
length $c-1$, by Lemma 5.1.4, this is impossible. Therefore $q>1$
and $x\in\L_{q-1}$ and $y=\eta$. Using Lemma 5.1.4 we see that
$u(a_{c-1,\b_{l_c}})$ and $u(a_{c-1,\b_{\t-1}})$ have the same
level. Write $$u(a_{c,\t})=\xi'+\eta'n,\ \ \ 1\le\xi'\le n,\ \
\eta'\in\bbZ,$$ and $$u(a_{c-1,\b_{\t-1}})=x'+y'n,\ \ \ 1\le x'\le
n,\ \ y'\in\bbZ.$$ then $x'\in\L_{q-1}$ and $y=y'$.

Since $u(a_{c,\t-1})>u(a_{c\t})$, $y'=\eta$, and $u(a_{cl_c})$,
$u(a_{c,\t-1})$ have the same level, we have $y'\ge\eta'$ and
$\xi'\in\L_{q'}$ for some $1\le q'<q$ if $y'=\eta'$.

When $y'>\eta'$, clearly (b) is true. If $y'=\eta'$ and $1\le
q'<q-1$, (b) is also obviously true. Now suppose that $y'=\eta'$
and $q'=q-1$. Using Lemma 5.4.2 we see that (b) is true. Thus (b)
is always true.

\medskip

 Similarly if $u(a_{cl_c}),...,u(a_{cm_c})$ have the same level and
$u(a_{cm_c}), u_{c+1}(e_c)$ have distinct levels and
$u_{c+1}(e_c)=u(a_{kl})$ for some $c<k\le r$ and $1\le l\le \l_k$, we
have $u(a_{c-1,\beta_{m_c}})>u_{c+1}(e_c)$. If
$u(a_{cl_c}),...,u(a_{cm_c})$ have the same level but $u(a_{cm_c})$
and $u_{c+1}(e_c)$ have distinct levels and
$u_{c+1}(e_c)=u(a_{1l_1}+n)$, using $u^{-1}\in\gg$ we can see that
$u(a_{c-1,\beta_{m_c}})>u_{c+1}(e_c)$. Thus in this case (b) is true.

\medskip

Note that $u(a_{cl_c})>u(a_{c-1,\beta_\t})$ for all $\beta_\t$. Thus
$m_{c-1}<\beta_\t$ for all $\beta_\t$ since $c-1\ge h+1$. Using
Corollary 6.2.4 we see that the following elements
$$\begin{array} {rl}
&w(a_{h1}),\ ...,\ w(a_{h,\l_{h}-\beta}),\
 w(a_{h,\l_h-\beta+1}) ,\\[3mm]
 &w(a_{h+1,l_{h+1}}),\ ...,\
 w(a_{h+1,m_{h+1}}),\ \   ...,\ \ w(a_{c-2,l_{c-2}}),\ ...,\
  w(a_{c-2,m_{c-2}}),\\[3mm]
& w(a_{c-1,l_{c-1} }),\ ...,\ w(a_{c-1,m_{c-1} }),\
w(a_{c-1,\beta_{l_c}}),\ w(a_{c-1,\beta_{l_c+1}}),\ ...,\
w(a_{c-1,\beta_{\theta-1}}),\\[3mm]
&w(a_{c,\theta-1}),\ ..., \
w(a_{cm_c}), \\[3mm]
&w(a_{c+1,l_{c+1}}),\ ...,\ w(a_{c+1,m_{c+1}}),\ \ ...,\ \
w(a_{r-1,l_{r-1}},\ ...,\ w(a_{r-1,m_{r-1}}),\\[3mm]
&w(a_{rl_r}),\
w(a_{r,l_r+1}),\ ...,\ w(a_{r,m_r}),\ w(a_{1,l_1+1}+n),\
...,\ w(a_{1\l_1}+n),
\end{array}$$
form an r-chain of $w$.

\medskip

Since
$$-(\beta-1)+m_r-k_r+\cdots+m_{h+1}-k_{h+1}+\l_{h}-1-k_1=0,$$
the above r-chain of $w$ has length

$$\begin{array}
{rl} L=&\l_{h}-\beta+1+m_{h+1}-l_{h+1}+1+\cdots+m_{c-1}-l_{c-1}+1\\
&+
m_c-l_c+2+m_{c+1}-l_{c+1}+1+\cdots+ m_r-l_r+1+\l_1-l_1.
\end{array}$$

\no  We have
$$\begin{array}
{rl}L&\ge
\l_{h}-\beta+1+m_{h+1}-k_{h+1}
+\cdots+m_{c-1}-k_{c-1}\\
&\qquad+m_c-k_c+1+m_{c+1}-k_{c+1}+\cdots+ m_r-k_r+1+\l_1-k_1 -1\\
&=\l_1+1\\
&>\l_1.
\end{array}$$

\no This is impossible since $w\in\gg$, so the lemma is true in
this case.

\medskip

Suppose that $c=h+1$. By assumption, we have

\medskip

\no(c) $u(a_{c'l_{c'}})$ and $u_{c'+1}(e_{c'})$ have the same level for all
$c=h+1<c'\le r$ with $l_{c'}\le m_{c'}$.

\medskip

Using Corollary 6.2.4 we get (note that $l_{h+1}\le m_{h+1}$ in our
case)

\medskip

\no(d)
 $u(a_{h+1,l_{h+1}})=w'(e_h)=w(a_{h,\l_h-\beta+1})$, and $w(a_{kl})$ and
$u(a_{kl})$ have the same level for all $r\ge k\ge h+2$ and
$\l_k\ge
l\ge 1$. And
$$w(\L_{h+1})-\{w(a_{h+1,m_{h+1}})\}=u(\L_{h+1})-
\{u(a_{h+1,l_{h+1}})\}.
$$

\medskip

Using Lemma 5.1.4, thus we have

\medskip

\no(e) there is  a complete r-antichain family $Y_u$ of $u$ such that
the r-antichain $B$ in $Y_u$ containing $u(a_{h+1,l_{h+1}})$ has
length $h+1$;

\medskip

\no(f) there is a complete r-antichain family $Y_w$ of $w$ such that
the r-antichain $C$ in $Y_w$ containing
$w(a_{h,\l_h-\b+1})=u(a_{h+1,l_{h+1}})$ has length $h$; and the
r-antichain $D$ in $Y_w$ containing
$w(a_{h+1,m_{h+1}})=u_{h+2}(e_{h+1})$ has length $h+1$.

\medskip

By Corollary 6.2.3 (e), $u(a_{kl_k})<u(a_{k+1,l_{k+1}})$ for
$k=1,2,...,h$. Using Lemma 5.4.2 we see

\medskip

\no(g) $u(a_{jl_j})$ and $u(a_{kl_k})$ have distinct levels whenever
$1\le j<k\le h+1$.

\medskip

Let $B'$ be an r-antichain in $Y_u$. Suppose that $B'$ has length
greater than $h+1$ or of length $h+1$ but does not contain
$u(a_{h+1,l_{h+1}})$. Assume that $u(a_{kl_k})$ is in $B'$ for
some $1\le k\le h$. By (g) and (d), we can find some
$u(a_{k\t_k})$ that has the same level as $u(a_{kl_k})$ and in an
r-antichain $B''$ in $Y_u$ of length less than the length of $B'$.
Replacing $u(a_{kl_k})$ by $u(a_{k\t_k})$ in $B'$ and replacing
$u(a_{k\t_k})$ by $u(a_{kl_k})$ in $B''$, and continuing this
process, finally we get a complete r-antichain family $Y'_u$ of
$u$ with the following property.

\medskip

\no(h)  Any r-antichain in $Y'_u$ of length greater than $h+1$ or
of length $h+1$ but not containing $u(a_{h+1,l_{h+1}})$, does not
contain any $u(a_{kl_k})$ for $k=1,2,...,h$. Moreover the r-antichain
in $Y'_u$ containing $u(a_{h+1,l_{h+1}})$ has length $h+1$.

\medskip

\no(i) It is harmless to require that $Y_u=Y'_u$. Then, for any
r-antichain $B'$ in $Y_u$ that has length greater than $h+1$ or
has length $h+1$ but does not contain $u(a_{h+1,l_{h+1}})$,
replacing $u(a_{c'l_{c'}})\ (h+1<c'\le r)$ by $w(a_{c'm_{c'}})$ if
$u(a_{c'l_{c'}})$ is $B'$, we get an r-antichain $C'$ of $w$. It
is harmless to assume all such $C'$ are in $Y_w$.

\medskip

Suppose that $u(a_{kq_k})$ are in $B$ for $k=1,2,...,h$. Write
$$u(a_{h+1,l_{h+1}})=x_{h+1}+y_{h+1}n,\ \ \ 1\le x_{h+1}\le n,\ \
y_{h+1}\in\bbZ,$$ $$u(a_{kq_k})=x_k+y_kn,\ \ \ 1\le x_{k}\le n,\ \
y_{k}\in\bbZ,\ \ k=1,2,...,h.$$ By (e), (f) and Lemma 5.1.4 we get

\medskip

\no(j) $x_{h+1}$ is in $\L_g$ for some $1\le g\le h$ and
  $x_{h+1-g}$ is in $\L_{h+1}$.

\medskip

Since $$u(a_{h+1,l_{h+1}})>u_{h+2}(e_{h+1})=w(a_{h+1,m_{h+1}})$$ and
they have distinct levels, by Corollary 6.2.4, if $q_{h+1-g}\ne
l_{h+1-g}$ then \break $u(a_{h+1-g,q_{h+1-g}})$ is in $w(\L_{h+1-g})$,
if $h+1-g\ge 2$ then $u(a_{h+1-g,q_{h+1-g}})$ is in $w(\L_{h+1-g})$ or
in $w(\L_{h+1-g-1})$. Then in both cases, by (i) and (f),
   $u(a_{h+1-g,q_{h+1-g}})$ is not in any r-antichain in $Y_w$ of length
greater than $h$. By Lemma 5.1.4, this is impossible. Hence
$h+1-g=1$ and $q_1=l_1$. But we have
$u(a_{h+1,l_{h+1}})>u_{h+2}(e_{h+1})\ge u(a_{1l_1})+n$, this
contradicts that $B$ is an r-antichain of $u$. Therefore
$u(a_{h+1,l_{h+1}})$ and $u_{h+2}(e_{h+1})=w(a_{h+1,m_{h+1}})$
have the same level.

\medskip

The lemma is proved

\bigskip

\no{\bf Corollary 6.3.5.} {\sl  $\vare_k(u)=\vare_k(w)$ whenever $k\ne i$.
}

\bigskip

\itp When $k<i$, this is Lemma 6.3.2 (b).  Using Lemma 6.3.4 and
Corollary 6.2.4 we see that $w(a_{cl})$ and $u(a_{cl})$ have the
same level if $c>r_i$ and $\l_c\ge l\ge 1$. Therefore
$\vare_k(u)=\vare_k(w)$ if $k>i$. The corollary is proved.

\bigskip

\no{\bf Lemma  6.3.6.} {\sl Let $u,v\in \gg$ be as in Prop. 6.2.1.
Keep the notation in 6.2.2. Assume that $w\in\gg$ and
deg$f_{u,v,w}=a(w_\l)$. Then}

\smallskip

\no(a) {\sl $u(a_{1l_1}) $ and $u(a_{1,l_1-1})$ have distinct levels and
$1\le l_1\le\l_h$.}

\no(b) {\sl $u(a_{cl_c})$ and $u(a_{c,l_c-1})$ have distinct levels for
$c=2,3,...,r$.}

\no(c) {\sl There is a complete r-antichain family $Z$ of $u$ such that
the r-antichain $B$ in $Z$ containing $u(a_{1l_1}) $ has length $h$.
(Recall that $r_i=h$.)}

\bigskip

\itp (a) By Lemma 6.3.2 (a) we have $1\le l_1\le\l_h$. If $u(a_{1,l_1-1})$ and
$u(a_{1l_1})$ have the same level then $$u(a_{1,l_1-1})=a_{kl}+qn\
\ \ \text{ and }u(a_{1l_1})=a_{kl'}+qn$$ for some $r\ge k\ge 1$,
$\l_k\ge l>l'\ge 1 $ and some $q\in\bbZ$. Then
$w^{-1}(a_{kl'})=\xi'-(q+1)n$ for some $1\le
\xi'\le n$, and $w^{-1}(a_{kl})=\xi-qn$ for some $1\le\xi\le\l_1$.
Thus $w^{-1}(a_{kl'})<w^{-1}(a_{kl})$. This contradicts that
$w\in\gg$. Therefore $u(a_{1,l_1-1})$ and $u(a_{1l_1})$ have distinct
levels.

\medskip

(b) Similarly we see that (b) is true.

\medskip

(c) Let $A_u$ (resp. $A_w$) be the set of elements in $u(\L_1)$
(resp. $w(\L_1)$) that have the same level as $u(a_{1l_1})$. Let
$Z$ be a complete r-antichain family of $u$. Assume that (c) is
not true. By Lemma 5.1.4, then any r-antichain in $Z$ that
contains one element in $A_u$ has length different from $h$. By
Corollary 6.2.3 (e) and Lemma 5.4.2, $u(a_{2l_2})$ and
$u(a_{1l_1})$ have distinct levels. Now $w(\L_1)$ is the union of
$\{u(L_1)-\{u(a_{1l_1})\}$ and $u(a_{2l_2})$. Thus
$|A_u|=|A_w|+1$. By Lemma 5.1.4, this forces that
$\vare_k(w)\ne\vare_k(u)$ for some $k\ne i$. This is impossible by
Corollary 6.3.5. So there must be some r-antichain in a complete
r-antichain family of $u$ that hash length $h$ and contains
$u(a_{1l_1})$.

\bigskip

\def\ci{\mathcal I}

\no{\bf Proposition 6.3.7.} {\sl Let $u,v\in \gg$ be such that all
components of $\vare(u)$ and $\vare(v)$ are nonnegative. Assume
that $$\vare_1(u)=\vare_2(u)=\cdots=\vare_{i-1}(u)=0$$ and
$$\vare_{i1}(v)=1\text{ \ and \ }\vare_{kl}(v)=0\ \ \text{if
$(k,l)\ne( i,1)$ }$$ Then $$t_ut_v=\sum t_w,$$ where $w$ runs
through the set $$\begin{array} {rl} \ci=\{w\in\gg \ |\ &
\vare_{ij}(w)=\vare_{ij}(u)+1\ \text{for some $1\le j\le n_i$,}\\
&\qquad\text{ and } \vare_{kl}(w)=\vare_{kl}(u)\ \text{if
$(k,l)\ne(i,j) $}\}.\end{array}$$}

\bigskip

\itp Suppose that $w\in\ci$ and $\vare_{ij}(w)=\vare_{ij}(u)+1$.
Then $ \vare_{i,j-1}(u)>\vare_{ij}(u)$. (We understand that
$\vare_{i0}(u)=\infty.$) Let $Z$ be a complete r-antichain family
of $u$ and $B$ an r-antichain in $Z$ that has length $r_i$ and
provides $\vare_{ij}(u)$. Let $u(a_{kl_k}$) be in $B$ for
$k=1,...,h$. We may choose $l_k$ as big as possible. By Lemma
5.1.4, then $u(a_{k,l_k-1})$ and $u(a_{kl_k})$ have distinct
levels. According to Lemmas 5.4.4 and 5.2.3, we see that $1\le
l_1\le\l_h$.

Now suppose that $1\le l_1\le \l_{h}$ satisfies (a) and (c) in
 Lemma 6.3.6. Assume that $B$ is an r-antichain in a complete
  r-antichain family
 $Z$ that has length $r_i$ and provides
$\vare_{ij}(u)$. Then obviously we have $
\vare_{i,j-1}(u)>\vare_{ij}(u)$. We need show that there exist unique
$k_1,\ k_c, \ 1\le l_c\le
\l_c+1$ for $c=2,3,...,r$ such that $k_c+1=l_c$ for all $1\le c\le r$ and
the corresponding $w$ in 6.2.2
satisfies that (1) $w\in\gg$ and (2) deg$f_{u,v,w}=a(w_\l)$.

\medskip

Keep the notation in 6.2.2.

\medskip

If $c\ge h+1$, we choose $ \ 1\le l_c\le m_c+1$ inductively so that
(1) $u(a_{cm_c})>u_{c+1}(e_c)>u(a_{c,m_c+1}),$ (2) $u(a_{cl_c})$ and
$u_{c+1}(e_c)$ have the same level but $u(a_{c,l_c-1})$ and
$u(a_{cl_c})$ have distinct levels, noting that if no such $l_c$ we
set $l_c=m_c+1$. Then set $k_c=l_c-1$.

\medskip

Let $u(a_{h\a_h}),...,u(a_{2\a_2}), u(a_{1l_1})$ be the r-antichain
$B$ in Lemma 6.3.6 (c). We may require that $u(a_{c\a_c})$ and
$u(a_{c,\a_c-1})$ have distinct levels for $c=2,...,h$. Then we set
$l_c=\a_c$ for $c=2,...,h$.

\medskip

According to the arguments in 6.2.2, it suffices to prove the
following three assertions.

\bigskip

\no(a) $l_h\le l_{h-1}\le\cdots\le l_1$.

\medskip

\no(b) $u(a_{c-1,l_c-1})>u(a_{cl_c})$ for $c=2,3,...,h$.

\medskip

\def\uh{u_{h+1}(e_h)}

\no(c) Let $1\le \beta\le \l_h$ be such that
$u(a_{h,\l_h-\beta})>\uh
>u(a_{h,\l_h-\beta+1})$. Then $$\l_h-\beta+m_{h+1}-l_{h+1}+1+
\cdots+m_r-l_r+1=l_1-1.$$
(Note that $u(a_{hl_h})<u(a_{1l_1})+n\le \uh$. So $\l_h-\beta+1\le
l_h$.)

\bigskip

For all $1\le c\le h$ we write
$$u(a_{cl_c})=\xi_{c}+\eta_{c}n,\qquad\text{ where }1\le\xi_{c}\le
n\text{ and }\eta_{c}\in\bbZ.$$ Let $Z$ be a complete r-antichain
family of $u$ containing $B$.

\bigskip

Now we argue for (a). If $l_c>l_{c-1}$ for some $2\le c\le h$, then
there is an r-antichain $C$ in $Z$ such that $C$ contains
$u(a_{cp_c})$ and $u(a_{c-1,p_{c-1}})$ for some $1\le p_c<l_c$ and
$l_{c-1}<p_{c-1}\le\l_{c-1}$. Write
$$u(a_{cp_c})=x+yn\text{ \  and\ \ }u(a_{c-1,p_{c-1}})=x'+y'n,$$
 where $1\le x,x'\le n$ and $y,y'\in\bbZ$.
Assume that $$x\in\L_j,\ \ x'\in\L_{k}\text{ \ and\ \ }\xi_c\in \L_q
.$$

\medskip

If $q >1$, by Lemma 5.1.4, then $\xi_{c-1}\in\L_{q-1}$ and
$\eta_c=\eta_{c-1}$. Since $u(a_{cp_c})\ge u(a_{c,l_c-1})$,\ \
$u(a_{c,l_c-1})$ and $u(a_{cl_c})$ have distinct levels, and
$u(a_{c,l_{c-1}})>u(a_{c-1,p_{c-1}})$ (see Lemma 5.4.4), we see

\medskip

\no(d) either $j>q$ and $ y\ge\eta_c$, or $y>\eta_c$; and
either $k\le q-1$ and $ y'\le
\eta_{c}$, or $y'<\eta_c$.

\medskip

If $j>1$, using Lemma 5.1.4 we see that $k=j-1$ and $y=y'$. This
contradicts  (d). If $j=1$, then $y>\eta_c$. By Lemma 5.3.2,
$\eta_c\ge 0$, thus $\vare(C)>0$, see section 5.2 for the
definition of $\vare(C)$. Since $\vare_m(u)=0$ for $1\le m<i$, the
length of $C$ is not less than $h$. By Lemma 5.1.4, then $k\ge h$.
But by Lemma 5.1.4, $q\le h$. Hence by (d) we see $y'<\eta_c$.
Therefore $u(a_{cp_c})\ge u(a_{c,l_c-1})+n$. This contradicts that
$C$ is an r-antichain of $u$.

\medksip

If $q=1$, then $\xi_{c-1}\in\L_h$ and $\eta_c=\eta_{c-1}+1$. Thus
$y\ge \eta_c$ and $y'\le\eta_c-1$. Since $C$ is an r-antichain in
$Z$, by Lemma 5.1.4 we must have $y=\eta_c$, $y'=\eta_c-1$ and
$j=1$. This contradicts that $u(a_{cl_c})$ and $u(a_{cp_c})$ have
distinct levels.

\medskip

Thus if $l_c>l_{c-1}$ for some $2\le c\le h$ we would be led to a
contradiction. Therefore (a) is true.

\medskip

Now we show (b). If $u(a_{c-1,l_c-1})<u(a_{cl_c})$ for some $2\le
c\le h$, using Lemma 5.4.2, we see that $u(a_{c-1,l_c-1})$ and
$u(a_{cl_c})$ have distinct levels. Write
$$u(a_{c-1,l_c-1})=\xi'+\eta'n,\ \ \ \ \xi'\in\L_{k'}\text{ and }
\eta'\in\bbZ.$$
 Recall that  $\xi_c\in\L_q$. Then

 \medskip

 \no(e) either $k'<q$ and $
\eta'\le\eta_c$,  or $\eta'<\eta_c$.

\medskip

Since $l_c\le l_{c-1}$ and $u(a_{c-1,l_{c-1}-1})$,
$u(a_{c-1,l_{c-1}})$ have distinct levels, if $q\ge 2$, by Lemma
5.1.4, we have that $k'>q-1$ and $ \eta'\ge\eta_{c-1}=\eta_c$,  or
$\eta'>\eta_{c-1}=\eta_c$. This contradicts (e). Therefore $q=1$
and $\eta'<\eta_c$. Then $\eta'\ge\eta_{c-1}$ implies that
$\eta'=\eta_{c-1}=\eta_c-1$. Moreover we have $k'>h$ since
$u(a_{c-1,l_{c-1}-1})$, $u(a_{c-1,l_{c-1}})$ have distinct levels
and $\xi_{c-1}\in\L_{h}$. Then it is easy to see that the
r-antichain $D$ in $Z$ containing $u(a_{c-1,l_c-1})$ does not
contain any of $u(a_{c1})$,...,$u(a_{c,l_c-1})$ since none of them
is in $\L_1+ (\eta'+1) n$ or $\L_{k'+1}+\eta' n$.

Thus there is an r-antichain $E$ in $Z$ such that $E$ contains
$u(a_{cp_c})$ and $u(a_{c-1,p_{c-1}})$ for some $1\le p_c<l_c$ and
$l_{c}-1< p_{c-1}\le\l_{c-1}$. As before, write
$u(a_{cp_c})=x+yn$, $u(a_{c-1,p_{c-1}})=x'+y'n$, where $1\le
x,x'\le n$ and $y,y'\in\bbZ$. Assume that $x\in\L_j,\
x'\in\L_{k}$. From the arguments above we see that $y'\le
\eta'=\eta_c-1$ and $y\ge\eta_c$. Thus by Lemma 5.1.4 we must have
$y=\eta_c$ and $j=1$. This is impossible since $u(a_{c,l_c-1})$
and $u(a_{cl_c})$ have distinct levels. We proved (b).

\medskip

\def\a{\alpha}
Now we prove (c). Let $U$ be the set consisting of $u(a_{11}),
u(a_{12}),...,$ $u(a_{1,l_1-1})$, and let $V$ be the union of the
two sets  $\{ u(a_{h1}),...,$ $u(a_{h,\l_h-\beta})\}$ and
$\{u(a_{c\delta})\ |\ h+1\le c\le r,\ 1\le\delta\le\l_c,\ \text{
and $u(a_{c\delta})$}$\text{ has the same level}{ as}{
$u(a_{1l_1})+n$} \}. Let $C$ be an r-antichain of $u$.  We shall
show that if $C$ contains some element in $U$, then $C$ contains
exactly one element of $V$, and vice versa.

\medskip

Let $\a\le\l_h-\beta$. Then $u(a_{h\a})>\uh$. Let $C$ be the
r-antichain in $Z$ containing $u(a_{h\a})$. Assume that
$u(a_{1\gamma})\in C$. We claim that $\gamma<l_1$. Write
$$u(a_{h\a})=\xi+\eta n,\ \ u(a_{1\gamma})=\xi'+\eta'n,$$ where
$1\le\xi,\xi'\le n,\
\eta,\eta'\in\bbZ$. Assume that  $\xi\in\L_m$ and
$\xi_{h+1}\in\L_{q_{h+1}}$. Then either $m>q_{h+1}$ and
$
\eta=\eta_{h+1}$, or $\eta>\eta_{h+1}.$

\medskip

Assume that $m>q_{h+1},\ \eta_{h+1}=\eta$. By Lemma 5.1.4, then
$\eta'=\eta-1$ and $\xi'\in\L_{m'}$ for some $m'>m$ whenever
$m<h$. If $m\ge h$, then $\eta'=\eta$. In any case we have
$u(a_{1\gamma})>u(a_{1l_1})$ since $\uh\ge u(a_{1l_1})+n$. So we
have $\gamma<l_1$.

If $\eta>\eta_{h+1}$, then $\eta'\ge \eta_{h+1}$, so we have
$u(a_{1\gamma})>u(a_{1l_1})$ since $\uh\ge u(a_{1l_1})+n$. Thus we
also have $\gamma<l_1$ in this case.

\medskip

Now suppose that $c>h,\ 1\le \delta\le \l_c$ and $u(a_{c\delta})$ has
the same level as $u(a_{1l_1})+n$ and $D$ is the r-antichain in $Z$
containing $u(a_{c\delta})$. Assume that $u(a_{1\gamma})$ is in $D$.
We claim that $\gamma<l_1$. We may write
$$u(a_{c\delta})=\zeta+(\eta_1+1) n,\ \ \ \ u(a_{1\gamma})=\xi'+\eta'n,$$
where $1\le\zeta,\xi'\le n,\
 \eta'\in\bbZ$.  Assume that  $\zeta\in\L_m$. Then either $m\ge c$,
 $\xi'\in\L_{m-c+1}$
and $\eta'=\eta_1+1$; or $m<c$, $\eta'=\eta_1$, and $\xi'\in\L_{m'}$
for some $m'>m$. In any case we have $u(a_{1\gamma})>u(a_{1l_1})$.
Hence $\gamma<l_1$.

\medskip

Obviously if an r-antichain in $Z$ contains some $u(a_{c\delta})$ for
some $c>h$ and $u(a_{c\delta})$ has the same level as $u(a_{1l_1})+n$,
then $C$ does not contains any $u(a_{h\a'})$ with $\a'\le\l_h-\beta$,
since $u(a_{h\a'})>\uh\ge u(a_{c\delta})$.

\medskip

Now suppose that $\gamma<l_1$ and $C$ is the r-antichain in $Z$
containing $u(a_{1\gamma})$. We claim that either $C$ contains some
$u(a_{h\a})$ for some $\a\le\l_h-\beta$ or $C$ contains some
$u(a_{c\a})$ for some $c>h$ and $u(a_{c\a})$ has the same level as
$u(a_{1l_1})+n$.

\medskip

As before write $u(a_{1\gamma})=\xi'+\eta'n$. Assume that
$\xi'\in\L_{m'}$ and $\xi_1\in\L_{q_1}$. Then either $m'>q_1$ and
$ \eta'=\eta_1$, or $\eta'>\eta_1$, since $u(a_{1\gamma})$ and
$u(a_{1l_1})$ have distinct levels. Since $\vare_j(u)=0$ for
$j<i$, by Lemma 5.1.4 we see that the length of $C$ is either $h$
or greater than $h$. Let $u(a_{h\a})=\xi+\eta n$ be in $C$, where
$\xi\in\L_q$ for some $1\le q\le r$ and $\eta\in\bbZ$.

\medskip

Assume that $\eta'>\eta_1$. Suppose that $u(a_{h\a})<\uh$. Then
$\eta=\eta_1+1$ since $\eta_1+1\ge\eta\ge\eta'>\eta_1$ and
$q<q_1$. By Lemma 5.1.4, $q\ge h$ and $q_1\le h$. A contradiction.
Therefore we have $u(a_{h\a})>\uh$ and we are done in this case.

Assume that $m'>q_1$ and $\eta'=\eta_1$. If $u(a_{h\a})>\uh$ then
we are done. Now suppose that $u(a_{h\a})<\uh$. Then we have two
cases, (1) $\eta=\eta_1+1$ and $q<q_1$, (2) $\eta=\eta_1$ and
$q>m'>q_1$. When $\eta=\eta_1+1$, since $q<q_1<m'$, by Lemma 5.1.4
   we see that there
is $r\ge c>h$ such that $u(a_{c\a'})=\xi''+\eta''n\in C$ for some
$1\le\a'\le\l_c$ with $\xi''\in\L_{q_1}$ and $\eta''=\eta$. Thus
$u(a_{c\a'})$ has the same level as $u(a_{1l_1})+n$. When
$\eta=\eta_1$ and $q>m'>q_1$, using Lemma 5.1.4 we see there is
some $r\ge c>h$ such that $u(a_{c\a'})=\xi''+\eta''n\in C$ for
some $1\le \a'\le\l_c$ with $\xi''\in\L_{q_1}$ and
$\eta''=\eta+1$.

Therefore (c) is true.

\medskip

The proposition is proved.

\bigskip

\no{\bf Proposition 6.3.8.} {\sl Let $u,v\in \gg$ be such that all
components of $\vare(u)$ and $\vare(v)$ are non-positive. Assume
that $$\vare_1(u)=\vare_2(u)=\cdots=\vare_{i-1}(u)=0$$ and
$$\vare_{in_i}(v)=-1\text{\ \  and \ }\vare_{kl}(v)=0$$ if
$(k,l)\ne(i,n_i)$. Then $$t_ut_{v }=\sum t_w,$$ where $w$ runs
through the set $$\begin{array} {rl} \{w\in\gg \ |\ &
\vare_{ij}(w)=\vare_{ij}(u)-1\ \text{for some $1\le j\le n_i$,}\\
&\qquad\text{ and } \vare_{kl}(w)=\vare_{kl}(u)\ \text {if
$(k,l)\ne(i,j) $}\}.\end{array}$$}

\bigskip

\itp Apply Prop. 6.3.7 and use 1.3 (b) and note that $J_{\gg}$ is
commutative (see Theorem 2.3.2 (b)).

\medskip

\section{The factorization  formula}

Now we can prove the factorization formula. We need a notation.
For $x,y\in\gg$ we define $x*y\in\gg$ to be the unique element $z$
in $\gg$ determined by $\vare(z)=\vare(x)+\vare(y)$.

\bigskip

\no{\bf Theorem 6.4.1.} {\sl Let $w\in\gg$ and    $w_i\in\gg$ (
$1\le i\le p$) be such that $\vare(w_i)=\vare_i(w)$ (see the
beginning of section 6.2 for the definition of $\vare_i(w)$). Then
$$t_w=t_{w_1}t_{w_2}\cdots t_{w_p}.$$}

\bigskip

\itp We show the result first in the case when $w$ satisfies
$\vare_{kl}(w)\ge 0$ for all $k,l$, and then in general.

\medskip

{\it Step 1.} Assume that all $\vare_{kl}(w)\ge 0$. It is sufficient
to prove that if $\vare_k(w)=0$ for $k=1,...,i$ and
$\vare(u)=\vare_i(u)$, then $t_ut_w=t_{u*w}$

Let $v\in\gg$ be such that $\vare_{kl}(v)=0$ for all pairs $(k,l)$
except $\vare_{i1}(v)=1$. Let $m\in\bbN$. According to Prop.
6.3.7, we have

\bigskip

\no(a)\qquad\qquad $t_v^mt_w=\sum_{u}\t_ut_ut_w=\sum_{u}\t_ut_{u*w},$

\medskip

\noindent where $u$ is in $\gg$ such that (1) all components of
$\vare(u)$ are nonnegative, (2) $\vare_i(u)=\vare(u)$, (3)
$\sum_{1\le j\le n_i}\vare_{ij}(u)=m$, and (4) $\t_u$ is given by
the following formula
$$V(\vare(v))^m=\sum_{u}\t_uV(\vare(u))\qquad \t_u\in\bbN.$$
Recall that $V(\vare(u))$ is an irreducible $F_\l$-module of
highest weight $\vare(u)$.

\bigskip

By 1.3 (d), the positivity of the structural coefficients
$\gamma_{x,y,z}$'s (see 1.3 (f)) and (a) we get

\bigskip

\no(b) If $\t_u\ne 0$, then $$t_ut_w=t_{x*w}$$ for some $x\in\gg$
with (1) $\vare_{kl}(x)\ge 0$ for all $k,l$, (2) $\sum_{1\le j\le
n_i}\vare_{ij}(x)=\sum_{1\le j\le n_i}\vare_{ij}(u)=m$ and (3)
$\vare_i(x)=\vare(x)$. Moreover if $u\ne u'$ and $t_{u'}t_w=t_{x'*w}$,
then $x\ne x'$.

\bigskip

Clearly $\theta_u$ is not zero if and only if (1) all components
of $\vare(u)$ are nonnegative, (2) $\vare_i(u)=\vare(u)$, (3)
$\sum_{1\le j\le n_i}\vare_{ij}(u)=m$.

\medskip

 Let $y\in\gg$ be such that
$\vare_{ij}(y)=k$ for all $j$ and other components of $\vare(y)$ are
0. Then $$t_yt_w=t_{z*w}$$ for some $z\in \gg$ with
$\vare_i(z)=\vare(z)$ and with $\sum_{1\le j\le
n_i}\vare_{ij}(z)=kn_i$.

By Prop. 6.3.7, we have $t_vt_y=t_{v*y}$. Hence
$$t_vt_yt_w=t_{v*y}t_w=t_{y'*w}$$ for some $y'\in
\gg$ with $\vare_i(y')=\vare(y')$ and with $\sum_{1\le j\le
n_i}\vare_{ij}(y')=kn_i+1$. Thus $$t_vt_{z*w}=t_{y'*w}.$$ By
Prop.6.3.7, this forces that $y=z$. Thus we have

\bigskip

\no(c) $t_yt_w=t_{y*w}$.

\bigskip

Now we define a total order on the subset $M_k$ of $\gg$ consisting of
all elements $u$ in $\gg$ with (1) all components of $\vare(u)$ are
nonnegative and are less than or equal to $k$ and (2)
$\vare_i(u)=\vare(u)$. Let $z,z'$ be in $M_k$. If
$\vare_{ij}(z)=\vare_{ij}(z')$ for $j=l+1,l+2,...,n_i$ but
$\vare_{il}(z)>\vare_{il}(z')$ then we define $z>z'$. This of course
introduces a total order on $M_k$. We have

\bigskip

\no(d) If $u\in M_k$ and $t_ut_w=t_{x*w}$, then $x\in M_k$.

\bigskip

Otherwise, by (b) we see that $\vare_{i1}(x)>k$. Let
$q=kn_i-\vare_{i1}(u)-\cdots-\vare_{in_i}(u)$. Write
$$t_v^qt_u=\sum_{u'}\xi_{u'}t_{u'},\qquad
\xi_{u'}\in\bbN.$$ By Prop. 6.3.7,  $\xi_y\ne 0$. By (c) and Prop. 6.3.7,
$t_{y*w}$ appears in $t_v^qt_ut_w$ with nonzero coefficient $\xi_y$
but $t_{y*w}$ appears in $t_v^qt_{x*w}$ with zero coefficient. Thus
 $t_v^qt_ut_w\ne t_v^qt_{x*w}$. This
contradiction shows that $x$ must be in $M_k$. (For a given
element $t=\sum a_zt_z\in J_{\gg}$, we say that $t_z$ appears in
$t$ with coefficient $a_z$.)

\bigskip

Now we claim that

\bigskip

\no(e) $t_ut_w=t_{u*w}$ for all $u$ in $M_k$.

\bigskip

   By (c) we see that (e) is true for the maximal
element in $M_k$. Assume that for any $z\in M_k$ with $z>u$ we have
$t_zt_w=t_{z*w}$. By (c) we may assume that $u$ is not maximal and we
can find $1\le j<n_i$ such that
$\vare_{ij}(u)>\vare_{i,j+1}(u)=\cdots=\vare_{in_i}(u)$.

We have $t_ut_w=t_{x*w}$ for some $x\in M_k$. Suppose that $u\ne x$.
Since for $x>u$ we have $t_xt_w=t_{x*w}$, using (b) we get $x<u$.

\medskip

Denote by $\tau_{il}$ the element in
$\bbZ^{n_1}\times\cdots\times\bbZ^{n_p}$ whose $(i,l)$-component is 1
and other components are 0. Consider the product
$$t_vt_ut_w=t_vt_{x*w}=\sum_{z}\eta_zt_{z*w},\ \  \eta_z\in\bbN.$$ Let
$u'\in M_k$ be determined by $\vare(u')=\vare(u)+\tau_{i,j+1}$, then
$t_{u'}$ appears in $t_vt_u$ with coefficient 1. By induction
hypothesis, $t_{u'*w}$ appears in $t_vt_{x*w}$ with nonzero
coefficient. According to Prop. 6.3.7, this forces that $\vare(x)=
\vare(u)+\tau_{i,j+1}-\tau_{i,n_i}$
since $x<u$.

If $\vare_{i1}(u)<k$ or there is $j'>j$ such that
$\vare_{ij'}(u)>\vare_{i,j'+1}(u)$, then $t_{v'}$ appears in $t_vt_u$
with coefficient 1, where $v'\in\gg$ with
$\vare(v')=\vare(u)+\tau_{i1}$ or $\vare(v')=\vare(u)+\tau_{i,j'+1}$.
Obviously $v'\in M_k$ and $v'>u$, by induction hypothesis, then
$t_{v'*w}$ appears in $t_vt_{x*w}$ with coefficient 1. But this is
impossible since by Prop. 6.3.7, if $x'\ne u'*w$ and $t_{x'}$ appears
in $t_vt_{x*w}$ with non zero coefficient, then
$\vare_{in_i}(x')=\vare_{in_i}(u)-1$.

\medskip

Thus we must have $$\vare_{i1}(u)=\cdots=\vare_{ij}(u)=k>
\vare_{i,j+1}(u)=\cdots=\vare_{in_i}(u).$$

\medskip

Define $\psi(\sum k_zt_z)=\sum k_z$.

\medskip

If $n_i-j>2$ or $n_i-j= 2$ but $\vare_{ij}(u)>\vare_{i,j+1}+1$, then
$\psi(t_vt_ut_w)=2$ but $\psi (t_vt_{x*w})\ge 3$. This is impossible.
If $n_i=j+1$ we have $x=u$ since
$\vare(x)=\vare(u)+\tau_{i,j+1}-\tau_{i,n_i}$.

 If  $n_i=j+2$ and $\vare_{ij}(u)=k=\vare_{i,j+1}(u)+1$, we can prove
directly that the assertion (e) is true in this case. In fact, from
$t_yt_w=t_{y*w}$ and by considering $t_vt_{y'}t_w$ we see that
$t_{y'}t_w=t_{y'*w}$ if $\vare(y')=\vare(y)-\tau_{in_i}$. Then by
considering $t_vt_{y''}t_w$ we see that $t_{y''}t_w=t_{y''*w}$ if
$\vare(y')=\vare(y)-2\tau_{in_i}$ and $k\ge 2$. Thus by considering
$t_vt_{u}t_w$ we see that $t_ut_w=t_{u*w}$ in this case.

\medskip

The theorem is proved when all $\vare_{kl}(w)\ge 0$.

\medskip

{\it Step 2.} Now $w\in\gg$ is arbitrary. We can find $q\in\bbN$ such
that $\vare_{kl}(w)+qr_k\ge 0$ for all $k,l$. Let $u=\om^{qn}w$ and
$u_i\in\gg$ with $\vare(u_i)=\vare_i(u)$. Then
$\vare_{kl}(u)=\vare_{kl}(w)+qr_k\ge 0$ for all $k,l$. By step 1 we
have

\bigskip

\no(f) $t_u=t_{u_1}t_{u_2}\cdots t_{u_p}$.

\bigskip

Obviously we have

\bigskip

\no(g) $t_w=t_ut_{\om^{-qn}w_\l}=t_{u_1}t_{u_2}\cdots t_{u_p}
t_{\om^{-qn}w_\l}.$

\bigskip

Since $\om^{qn}w_\l\in\gg$ and $\vare_{kl}(\om^{qn}w_\l)=qr_k\ge 0$,
by step 1, we have

\def\th{\theta}

\bigskip

\no(h) $t_{\om^{qn}w_\l}=t_{\th_1}t_{\th_2}\cdots t_{\th_p}$, where
$\th_k\in\gg$ is determined by $\vare(\th_k)=\vare_k(\om^{qn}w_\l)$.

\bigskip

Using 1.3 (b) and Theorem 2.3.2 (b) we get

\bigskip

\no(i) $t_{\om^{-qn}w_\l}=t_{\th_1^{-1}}t_{\th_2^{-1}}\cdots
t_{\th_p^{-1}}$.

\bigskip

Obviously we have $t_{\om^{qn}w_\l}t_{\om^{-qn}w_\l}=t_{w_\l}$. Using
1.3 (c) and 1.3 (f) we see that $t_{\th_k}t_{\th_k^{-1}}=t_{w_\l}$ for
all $k$. Thus we have

\bigskip

\no(j) $t_{\th_k}^{-1}=t_{\th_k^{-1}}$.

\bigskip

According to step 1 we have

\bigskip

\no(k) $t_{\th_1}t_{\th_2}\cdots t_{\th_{k-1}}t_{u_k}t_{\th_{k+1}}\cdots
t_{\th_p}=t_{v_k}$ for all $k$, where
$v_k=\th_1*\th_2*\cdots*\th_{k-1}*u_k*\th_{k+1}*\cdots*\th_p$.

\bigskip

Note that $\vare(v_k\om^{-qn})=\vare(w_k)=\vare(u_k)-\vare(\t_k)$
for all $k$. Using Theorem 5.2.6 (b), we get
 $v_k\om^{-qn}=w_k$ for all $k$. Then obviously we have
$t_{v_k}t_{\om^{-qn}w_\l}=t_{w_k}$. Using (f)-(k) we get
$$\begin{array} {rl} t_w&=t_ut_{\om^{-qn}w_\l}\\[3mm]
&=\prod_{k=1}^p(t_{\th_1}t_{\th_2}\cdots t_{\th_{k-1}}t_{u_k}t_{\th_{k+1}}\cdots
t_{\th_p}t_{\th_1^{-1}}t_{\th_2^{-1}}\cdots t_{\th_p^{-1}})\\[3mm]
&=\prod_{k=1}^pt_{v_k}t_{\om^{-qn}w_\l}\\[3mm]
&=t_{w_1}t_{w_2}\cdots t_{w_p}.
\end{array}$$

The theorem is proved.

\chapter{A Multiplication Formula in  $J_{\Ga_\lambda\cap\Ga_\lambda^{-1}}$}


\def\tt{\tilde T}
Let $\l$ be as in section 2.2 and $n_i$ be as in section 4.1. Then
$GL_{n_i}(\bbC)$ is the $i$th reductive component of $F_\l$. Let
$\vare=(\vare_{11},...,\vare_{pn_p})$ and
$\vare'=(\vare'_{11},...,\vare'_{pn_p})$ be two elements in
Dom$(F_\l)$ such that $\vare_{kl}=\vare'_{kl}=0$ whenever $k\ne i$
and $\vare_{i1}=2,\ \vare_{i2}=\cdots=\vare_{in_i}=0$. Let
$V(\vare)$ and $V(\vare')$ be two irreducible $F_\l$-modules of
highest weight $\vare$ and $\vare'$ respectively. Then we have a
simple formula for the product $V(\vare)V(\vare')$ in $R_{F_\l}$
(see 4.2 (g)). In this chapter we will establish a multiplication
formula in $J_{\Ga_\l\cap\Ga_\l^{-1}}$ that is corresponding to
the formula for $V(\vare)V(\vare')$, see Theorem 7.2.2. To do this
we compute some product $\tilde T_u\tilde T_v$ in section 7.1.
Then in section 7.2 we prove our formula. In Chapter 8, using this
formula and the factorization formula (Theorem 6.4.1) we show that
 Conjecture 2.3.3 is true.

\medskip

\section{A computation for some $\tt_u\tt_v$}

\def\ti{v_{i1}}

Let $u,v$ be in $\gg$ such that $\vare(u)=\vare_i(u)$ and
$\vare_{i1}(v)=2$ and $\vare_{kl}(v)=0$ whenever $(k,l)\ne (i,1)$.
In this section we will compute the product $\tt_u\tt_v$ provided
that all components of $\vare(u)$ are non-negative. We show that
$f_{u,v,w}$ have degree less than or equal to $a(w_\l)$ for the
$u,v$ and all $w\in W$. Keep the notation in Chapter 5. Let
$v_{i1}$ be the $u_{i1}$ in section 5.5. We have

\bigskip

\noindent{\bf Lemma 7.1.1.} {\sl Let $1\le k\le r$. Suppose
that $e_{k-1}+1\le q\le e_k-2$. Then }

\smallskip

\noindent (a) $\tt_{s_q}\tt_{\ti}=\tt_{\ti}\tt_{s_q}$,

\no (b) $\tt_{s_q}^{-1}\tt_{\ti}=\tt_{\ti}\tt_{s_q}^{-1}$.

\bigskip

\itp (a) According to section 5.5 we have
$v_{i1}=\tau_{\l_1}s(e_1,e_2,...,e_h).$ Using this we see easily
  that $s_q\ti=\ti s_q$ and
$l(s_q\ti)=1+l(\ti)$. Thus
$$\tt_{s_q}\tt_{\ti}=\tt_{s_q\ti}=\tt_ {\ti s_q}=\tt_{\ti}\tt_{s_q}.$$
We proved (a). (b) follows from (a).

\bigskip

\def\a{2}

\noindent{\bf 7.1.2.}
Let $u$ be in $\gg$ such that (1) $\vare_i(u)=\vare(u)$ and (2) all
components of $\vare (u) $ are nonnegative. Let $v=\ti^2 w_\l$.
 By 5.5 (d-e), we have  $v\in\gg$, $ \vare_{kl}(v)=0$ if $(k,l)\ne (i,1)$ and
 $\vare_{i1}(v)=2$.
 Now we compute  $\tt_u\tt_v$. Recall that we have set $e_k=\l_1+\cdots+\l_k$,
 $e_0=0$ and $r_i=h$.
 By the construction in section 5.3 we see

\bigskip

 \no ($*$) $u(a)=v(a)=w_\l(a)$ if either $e_h+1\le a\le n$ or
 $a=a_{kl}$ for some $k$ with $l>\l_h$.

\bigskip

 According to 5.5 (g) we get
$$\begin{array}
{rl} \ti=&\om s_{e_1-2}s_{e_1-3}\cdots s_1s_0\hat
s_{e_2-1}s_{e_2-2}\cdots s_{e_{1}}\\
&\times \hat
s_{e_{h-1}-1}s_{e_{h-1}-2}\cdots s_{e_{h-2}}\hat
s_{e_h-1}s_{e_h-2}\cdots s_{e_{ h-1}}\\
&\times s_{n-1}s_{n-2}\cdots s_{e_h}
\end{array}$$

 Set
$$\delta=s_{n-1}\cdots s_{e_h}.$$ Let $1\le j\le h-1$. Suppose that $z_{k}\le s_{e_k-2}\cdots
s_{e_{ k-1}}$ for $k=1,2,...,j$. Using 2.1.3 (f) we see

\bigskip

\no(a) $u\om z_{1}z_2\cdots z_js_q\le u\om z_{1}z_2\cdots z_j$ for all
$e_{j}\le q\le e_{j+1}-2$.

\bigskip

According to $(*)$ and Lemma 6.1.3, we have

$$
\tt_u\tt_{\ti}=\displaystyle\sum_{\st{\sc z_{1,1}\le s_{e_1-2}\cdots
s_{e_{0}}} {\st{\sc ...} {z_{h,1}\le s_{e_h-2}\cdots
s_{e_{h-1}}}}}(\prod_{k=1}^{h}\xi^{\l_k-1-l(z_{k,1})})\tt_{u\om
z_{1,1}\cdots z_{h,1}}\tt_\delta.\leqno{\hbox{(b)}}
$$

Using Lemma 5.3.2 and 2.1.3 (f) we see that $$l(u\om z_{1,1}\cdots
z_{h,1}\delta)=l(u\om z_{1,1}\cdots z_{h,1})+l(\delta).$$
 Thus we have
$$\tt_{u\om z_{1,1}\cdots z_{h,1}}\tt_\delta=\tt_{ u\om z_{1,1}\cdots
z_{h,1}
\delta}.$$

\smallskip

Let $$u'_1=u\om z_{1,1}\cdots z_{h,1}
\delta.$$ Then for $1\le k\le h-1$ we have
 $u'_1(e_k)=u(a_{k+1,p_{k+1,1}})$
for some $1\le p_{k+1,1}\le\l_{k+1}$ and $u'_1(e_h)=u(a_{1,p_{
1,1}})+n$ for some $1\le p_{1,1}\le\l_1$. Thus
$$z_{k,1}=s_{e_{k-1}+p_{k,1}-2}\cdots
 s_{e_{k-1}} y_{k,1}$$ for some $y_{k,1}\le s_{e_k-2}s_{e_k-3}
 \cdots s_{e_{k-1}+p_{k,1}}$.
 We have
 $$l(u'_1 y_{1,1}^{-1}\cdots y_{h,1}^{-1})=
 l(u'_1 )+l( y_{1,1}^{-1})+\cdots +l(y_{h,1}^{-1}).$$
 Let $u_1=u'_1 y_{1,1}^{-1}\cdots y_{h,1}^{-1}$. Then
 $$\tt_{ u'_1 }=\tt_{u_1}\tt_{y_{1,1}^{-1}}^{-1}
 \cdots\tt_{y_{h,1}^{-1}}^{-1}.$$

 Note that  $l(z_{k,1})=l(y_{k,1})+p_{k,1}-1$. Thus we have

$$
\tt_u\tt_{\ti}=\displaystyle\sum_{\st{\sc z_{1,1}\le s_{e_1-2}\cdots
s_{e_{0}}} {\st{\sc ...} {z_{h,1}\le s_{e_h-2}\cdots
s_{e_{h-1}}}}}(\prod_{k=1}^{h}\xi^{\l_k-p_{k,1}})\tt_{u_1}
(\xi^{-l(y_{1,1})}\tt_{y_{1,1}^{-1}}^{-1})\cdots(\xi^{-l(y_{h,1})}
\tt_{y_{h,1}^{-1}}^{-1}).\leqno{\hbox{(c)}}
$$

\smallskip

Moreover we have

\smallskip

$$u_1(a)=\cases u(a)&\text{ if } e_{k-1}+1\le a\le
e_{k-1}+p_{k,1}-1,\ \ k=1,2,...,h,\\ u(a+1)&\text{ if }
e_{k-1}+p_{k,1}\le a\le e_{k-1}-1,\ \ k=1,2,...,h,\\
u(a_{k+1,p_{k+1,1}})&\text{ if } a=e_k,\ \ k=1,2,...,h-1,\\
u(a_{1,p_{1,1}})+n&\text{ if } a=e_h,\\ u(a)&\text{ if } e_h<a\le
n.\endcases\leqno{\hbox{(d)}}$$

\smallskip

In particular we have

\no(e) $u_1(a)>u_1(b)$ if $e_k+1\le a<b\le e_{k+1}-1$ for some $0\le k\le
r-1$.

\bigskip

For any $w\in W$ we set $w(a_{k0})=\infty$ and
$w(a_{k,\l_k+1})=-\infty$ for all $k$.

\medskip

For $1\le k\le h$ choose $0\le q_{k,2}\le\l_k-1$ such that $$
u_1(a_{k,q_{k,2}})> u_1(e_k) \ge u_1(a_{k,q_{k,2}+1}).$$
 Using 2.1.3 (f) we get
  $$l(u_1\om \prod_{k=1}^h(s_{e_k-2}\cdots
s_{e_{k-1}+q_{k,2}}))=l(u_1)+\sum_{k=1}^h(\l_k-q_{k,2}-1).
\leqno{\hbox{(f)}}$$
Let $u'_2=u_1\om
\prod_{k=1}^h (s_{e_k-2}\cdots s_{e_{k-1}+q_{k,2}})$.
Then we have


\bigskip

$$\begin{array}
{rl}
 u'_2(a)&=\cases u_1(a+1)&\text{ if } e_{k-1} \le a\le
e_{k-1}+q_{k,2}-1,\ \ k=1,2,...,h,\\ u_1(e_k)&\text{ if }
a=e_{k-1}+q_{k,2}, \ \ k=1,2,...,h,\\ u_1(a)
&\text{ if
} e_{k-1}+q_{k,2}<a<e_k,\ \ k=1,2,...,h,\\ u_1(a+1)&\text{ if }
 e_h\le a\le n-1.\endcases\\
 u'_2\om^{-1}&=\cases u_1(a)&\text{ if } e_{k-1}+1\le a\le
e_{k-1}+q_{k,2},\ \ k=1,2,...,h,\\ u_1(e_k)&\text{ if }
a=e_{k-1}+q_{k,2}+1,\ \ k=1,2,...,h,\\ u_1(a-1)&\text{ if }
  e_{k-1}+q_{k,2}+1< a\le
e_{k},\
\ k=1,2,...,h,\\ u_1(a)
 &\text{ if }
 e_h+1\le a\le n.\endcases\end{array}\leqno{\hbox{(g)}}$$

\smallskip

As a consequence we have

\no(h) $u'_2(a)=u(a+1)$ if $e_h\le a\le n-1$ and $u'_2s_q\le u'_2$ if
$e_{k-1}\le q\le e_{k-1}+q_{k,2}-1$ for some $1\le k\le h$. Moreover,
$u'_2(a)>u'_2(b)$ if $e_k\le a<b\le e_{k+1}-1$ for some $0\le k\le
r-1$.

\bigskip

Thus we have

$$\tt_{u_1}\tt_{\ti}=\tt_{u'_2}\tt((\prod_{k=1}^hs_{e_{k-1}+q_{k,2}-1}\cdots
s_{e_{k-1}})\delta),\leqno{\hbox{(i)}}$$ and (here we need Lemma 6.1.3
and recall that $\tt(w)$ is set to be $\tt_w$)
$$\tt_{u_1}\tt_{\ti}=\displaystyle\sum_{\st{\sc z_{1,2}\le s_{e_{0}+q_{1,2}-1}\cdots
s_{e_{0}}} {\st{\sc ...} {z_{h,2}\le s_{e_{h-1}+q_{h,2}-1}\cdots
s_{e_{h-1}}}}}(\prod_{k=1}^{h}
\xi^{q_{k,2}-l(z_{k,2})})\tt(u'_2z_{1,2}\cdots z_{h,2})\tt_\delta.
\leqno\hbox{(j)}$$

\smallskip

We may check that $$l(u'_2 z_{1,2}\cdots z_{h,2}\delta)=l(u'_2
z_{1,2}\cdots z_{h,2})+l(\delta).$$ Therefore we have $$\tt_{u'_2
z_{1,2}\cdots z_{h,2}}\tt_\delta=\tt_{ u'_2 z_{1,2}\cdots z_{h,2}
\delta}.$$

\medskip

Assume that $$u'_2 z_{1,2}\cdots z_{h,2}
\delta(e_k)=u'_2\om^{-1}(a_{k+1,p_{k+1,2}})$$ for $k=1,...,h-1$ and
  $$u'_2 z_{1,2}\cdots z_{h,2}
\delta(e_{h})=u'_2\om^{-1}(a_{1,p_{ 1,2}})+n.$$ Then
$$z_{k,2}=s_{e_{k-1}+p_{k,2}-2}\cdots
 s_{e_{k-1}} y_{k,2}$$ for some $y_{k,2}\le
 s_{e_{k-1}+q_{k,2}-1}s_{e_{k-1}+q_{k,2}-2}\cdots s_{e_{k-1}+p_{k,2}}
 $.
  We have $$l(u'_2 z_{1,2}\cdots z_{h,2}
\delta y_{1,2}^{-1}\cdots y_{h,2}^{-1})=l(u'_2 z_{1,2}\cdots z_{h,2}
\delta)+l( y_{1,2}^{-1})+\cdots +l(y_{h,2}^{-1}).$$ Let $u_2=
u'_2 z_{1,2}\cdots z_{h,2}
\delta y_{1,2}^{-1}\cdots y_{h,2}^{-1}$. Then $$\tt_{
u'_2 z_{1,2}\cdots z_{h,2}
\delta}=\tt_{u_2}\tt_{y_{1,2}^{-1}}^{-1}
\cdots\tt_{y_{h,2}^{-1}}^{-1}.$$

Note that  $l(z_{k,2})=l(y_{k,2})
 +p_{k,2}-1$. Thus we have

$$
\tt_{u_1}\tt_{\ti}=\displaystyle\sum_{\st{\sc z_{1,2}\le s_{e_0+q_{1,2}-1}\cdots
s_{e_{0}}} {\st{\sc ...} {z_{h,2}\le s_{e_{h-1}+q_{h,2}-1}\cdots
s_{e_{h-1}}}}}(\prod_{k=1}^{h}
\xi^{q_{k,2}-p_{k,2}+1})\tt_{u_2}(\xi^{-l(y_{1,2})}\tt_{y_{1,2}^{-1}}^{-1})
\cdots(\xi^{-l(y_{h,2})}
\tt_{y_{h,2}^{-1}}^{-1}).
$$

Moreover we have

\bigskip

\no(k) $u_2(e_k)=u'_2\om^{-1}(a_{k+1,p_{k+1,2}})$ for
$k=1,...,h-1$ and $u_2(e_h)=u'_2\om^{-1}(a_{1p_{1,2}})+n$.

\bigskip

\no(l) $u_2(a)=u(a)$ for all $e_h<a\le n$, and $u_2(a)>u_2(b)$ if $e_k+1\le
a<b\le e_{k+1}-1$ for some $k=0,1,2,...,h-1$.

\bigskip

From the above discussion and using Lemma 7.1.1 we see that

$$\tt_u\tt_v=\displaystyle\sum_{u_\a}f'_{u,v,u_\a}\tt_{u_\a}
\prod_{m=2}^1((\xi^{-l(y_{1,m})}\tt_{y_{1,m}^{-1}}^{-1}
)\cdots(\xi^{-l(y_{h,m})}
\tt_{y_{h,m}^{-1}}^{-1}))\tt_{w_\l},\qquad f'_{u,v,u_\a}\in\A,
\leqno\hbox{(m)}$$
where
$f'_{u,v,u_\a}=\prod_{k=1}^h\xi^{\l_k-p_{k,1}+q_{k,2}-p_{k,2}+1}$ has
degree $$D=\sum_{k=1}^h(\l_k-p_{k,1}+q_{k,2}-p_{k,2}+1 ).$$

\bigskip

Using Lemma 5.4.4, (d)  and Lemma 5.1.1 (c)  we see that

\medskip

\no (n) $q_{k,2}\le p_{k+1,1}-1$ for $k=1,2,...,h-1$, and
$q_{h,2}\le p_{1,1}-1$.

\medskip

Since $u(a_{k+1,p_{k+1,1}})\le u'_2\om^{-1}(a_{k,p_{k,2}})$ for
$k=1,...,h-1$ and $u(a_{1,p_{ 1,1}})+n\le
u'_2\om^{-1}(a_{h,p_{h,2}})$, we have

\medskip

\no(o) Let $1\le k\le h-2$. We have
$$u_2(\L_{k})=(u(\L_k)-\{u(a_{k,p_{k,1}}),\
u'_2\om^{-1}(a_{k,p_{k,2}})\})\bigcup\{ u(a_{k+1,p_{k+1,1}}),\
u'_2\om^{-1}(a_{k+1,p_{k+1,2}})\}.$$ If
$u(a_{k+2,p_{k+2,1}})>u(a_{k+1,p_{k+1,1}})$, then $$1\le
p_{k+1,2}\le q_{k+1,2}+1\le p_{k+1,1}.$$ In this case we have
$$u'_2\om^{-1}(a_{k+1,p_{k+1,2}})\ge
u(a_{k+2,p_{k+2,1}})>u(a_{k+1,p_{k+1,1}}),$$ and $1\le
p_{k+1,2}\le q_{k+1,2}$ or $p_{k+1,2}=q_{k+1,2}+1$. Using (d), (g)
and Lemma 5.4.4 we see that

\no(o1)  among $u_2(a_{k,1}),...,u_2(a_{k,\l_k-1})$, at least
$\l_k-p_{k+1,2}$ of them are smaller than $u_2(a_{k\l_k})$.

Similarly, if $u(a_{k+2,p_{k+2,1}})<u(a_{k+1,p_{k+1,1}})$, then we
have

\no(o2)  among $u_2(a_{k,1}),...,u_2(a_{k,\l_k-1})$, at least
$\l_k-p_{k+1,2}-1$ of them are smaller than $u_2(a_{k\l_k})$, and
$q_{k+1,2}\le p_{k+2,1}-2$.

Similarly for $k=h-1,h$ we have

\no(o3) either  among $u_2(a_{k,1}),...,u_2(a_{k,\l_k-1})$, at least
$\l_k-p_{k+1,2}$ of them are smaller than $u_2(a_{k\l_k})$, or
 among $u_2(a_{k,1}),...,u_2(a_{k,\l_k-1})$, at least
$\l_k-p_{k+1,2}-1$ of them are smaller than $u_2(a_{k\l_k})$ and
$q_{k+1,2}\le p_{k+2,1}-2$, where we set $p_{h+k,2}=p_{k,2}$ and
$q_{h+k,2}=q_{k,2}$.

\medskip

Let $x$ be the element in $u_2W_\l$ of minimal length (see 6.1 for the
definition of $W_\l$). Then $u_2=xy$ for some $y\in W_\l$. Using (n)
and (o1-o3) we get

\medskip

\no(p) $D+l(y)\le l(w_\l)$.

\medskip

We also have

\medskip

\no(q)   $\tt_{y_{k,m}^{-1}}^{-1}\tt_{w_\l}=\tt_{y_{k,m}^{-1} w_\l}$ and
$l(y_{k,m}^{-1} w_\l)=l(w_\l)-l(y_{k,m}^{-1})$ for all $k,m$.

\medskip

\no(r)  For
any $w\in W, \ s\in S$ we have
$$\xi^{-1}\tt_s^{-1}\tt_w=\cases \xi^{-1}\tt_{sw}&\quad\text{ if } ws\le
w\\ \xi^{-1}\tt_{sw}-\tt_w&\quad\text{ if } ws\ge w.\endcases$$

\medskip

 Thus we see in the expression
$\tt_u\tt_v=\sum_wf_{u,v,w}\tt_w$, the degree $L$ of $f_{u,v,w}$ is at
most $D+l(y) $. Using (p) we see that
$$L\le a(w_\l)=l(w_\l).$$ Using (q) and (r) we see that

\medskip

\no(s) If $L=a(w_\l)$ then $y_{k,m}$ are the neutral element of $W$
 for all $k,m$, here we
set $p_{h+1,1}=p_{1,1}$.

\medskip

\section{A multiplication formula}

In this section we give a multiplication formula in $\jg$, based on
the computation in section 7.1.

\bigskip

\no{\bf Theorem 7.2.1.} {\sl Let $u,v$ be in $\gg$ such that all
components of $\vare(u)$ are non-negative, $\vare_i(u)=\vare(u)$,
$\vare_{i1}(v)=2$ and other components of $\vare(v)$ are 0. Then
$$t_ut_v=\sum_wt_w,$$ where $w$ runs through the set
$$\begin{array} {rl} \{w\in\gg \ &|\ \vare
(w)=\vare(u)+\tau_{ij_1}+\tau_{ij_2} \text{ for some $1\le j_2 \le
j_1\le n_i$ }\\ &\qquad\text{and }
\vare(u)+\tau_{ij_1}\in\text{Dom}(F_\l)
 \},\end{array}$$ see the proof of Theorem 6.4.1 for the definition
of $\tau_{ij}$.}

\bigskip

\itp Keep the notation in 7.1.2.

\medskip

Assume that $\vare(u)+ \tau_{ij_1}\in\text{Dom}(F_\l) $ and
$w\in\gg,\ \ \vare(w)=\vare(u)+\tau_{ij_1} +\tau_{ij_2}$. For all
$1\le k\le h$, set $$z_{k,1}=s_{e_{k-1}+j_1-2}\cdots
s_{e_{k-1}+1}s_{e_{k-1}}$$ and $$z_{k,2}=s_{e_{k-1}+j_2-2}\cdots
s_{e_{k-1}+1}s_{e_{k-1}}.$$ From the construction in section 5.3
and the arguments in 7.1.2 we see that $f_{u,v,w}$ has degree
$a(w_\l)$ and its leading coefficient is 1.

\medskip

Now suppose that $w\in\gg$ and $f_{u,v,w}$ has degree $a(w_\l$).
Let $x,y$ be in $\gg$ such that $\vare_{i1}(x)=1$,
$\vare_{i1}(y)=\vare_{i2}(y)=1$ and all other components of
$\vare(x) $ and $\vare(y)$ are 0. According to Prop. 6.3.7 we have
$t_x^2=t_y+t_v.$ Using Prop. 6.3.7 and the positivity 1.3 (f) we
see that $$\vare(w)=\vare_{i}(w)\text{ and }\sum_{1\le j\le
n_i}\vare_{ij}(w)=\sum_{1\le j\le n_i}\vare_{ij}(u)+2.\leqno{
{(\star)}}$$ Moreover, all $p_{k,m}$ ($1\le k\le h,\ m=1,2)$ are
not greater than $\l_h$. By Lemma 5.4.4, 7.1.2 (n-p),
$u(a_{k,p_{k+1,1}-1})>u(a_{k+1,p_{k+1,1}})$ for $k=1,2,...,h-1$
and $u(a_{h,p_{h+1,1}-1})>u(a_{ 1,p_{ 1,1}}+n)$. We claim that
$p_{1,1}=p_{2,1}=\cdots=p_{h,1}$.

Otherwise, $p_{k,1}\ne p_{k+1,1}$ for some $1\le k\le h$ (recall
that $p_{h+k,1}=p_{k,1}$). By $(\star)$ and the construction in
section 5.3, $w(\L_k)$ can not contain both $u(a_{k,p_{k+1,1}})$
and $u(a_{k+1,p_{k+1,1}})$ (we set
$u(a_{h+1,p_{h+1,1}})=u(a_{1,p_{1,1}})+n$). By the arguments in
7.1.2 (o) we see that

\no($*)$ $u(a_{k+1,p_{k+1,1}})$ is in $w(\L_{k-1})$ (we set $\L_0=\L_h$).

Using the construction in section 5.3, we must have
$p_{k,1}=p_{k+1,1}$. A contradiction, so the assumption
$p_{k,1}\ne p_{k+1,1}$ is not true.

In a complete similar way we see that
 $p_{1,2}=p_{2,2}=\cdots=p_{h,2}$.
 Since $w\in\gg$ and $\vare(w)=\vare_i(w)$,
we
must have $1\le p_{1,1},p_{1,2}\le n_i$. From the arguments in 7.1.2
we have $p_{1,2}\le p_{1,1}$. Thus $\vare
(w)=\vare(u)+\tau_{ij_1}+\tau_{ij_2}$ for $1\le j_2=p_{1,2}
\le j_1=p_{1,1}\le n_i$.

\medskip

The theorem is proved.

\bigskip

\no{\bf Theorem 7.2.2.} {\sl Let $u,v$ be in $\gg$ such that
$\vare_i(u)=\vare(u)$, all components of $\vare(u)$ are
non-negative, $\vare_{i1}(v)=\vare_{i2}(v)=1$ and other components
of $\vare(v)$ are 0. Then $$t_ut_v=\sum_wt_w,$$ where $w$ runs
through the set $$\{w\in\gg \ |\ \vare
(w)=\vare(u)+\tau_{ij_1}+\tau_{ij_2} \text{ for some $1\le j_2 <
j_1\le n_i$, }\}.$$}

\bigskip

\itp Let $x,y$ be in $\gg$ such that $\vare_{i1}(x)=1$,
$\vare_{i1}(y)=2$ and all other components of $\vare(x), \vare(y)$ are
0. According to Prop. 6.3.7 we have $t_x^2=t_y+t_v.$ Now considering
$t_ut_x^2$ and using Prop.6.3.7 and Theorem 7.2.1 we see that the
required result is true.

\chapter{The Based Rings $J_{\Gamma_\lambda\cap\Gamma_\lambda^{-1}}$ and $J_{\mathbf c}$}


\def\tt{\tilde T}
In this chapter we will prove Conjecture 2.3.3 using the formulas
in Chapters 6 and 7. This completes our proof of Lusztig
Conjecture for type $\tilde A_{n-1}$. In section 8.1 we give some
lemmas about multiplication in $J_{\gg}$. In section 8.2 we give a
proof for Conjecture 2.3.3 and give a summary of our proof of
Lusztig Conjecture on based ring for type $\tilde A_{n-1}$. Also
we give a few comments about a possible geometric realization of
$J_{\gg}$. In section 8.3 we consider the affine Weyl group
associated with  the projective linear group $PGL_n(\bbC)$. For
the affine Weyl group we show that Lusztig Conjecture on based
ring should have a weak form since the algebraic group
 $PGL_n(\bbC)$ is not simply connected. In section 8.4 we show that
 Lusztig Conjecture on based ring is true for the extended affine
  Weyl group associated with  the special linear group
$SL_n(\bbC)$.

\medskip

\section{Some lemmas}

\def\bzi{\bbZ_{\text{dom}}^{n_i}}

In this section we establish some lemmas about multiplication in
$J_{\gg}$. Let $N$ be the subset of $\gg$ consisting of all
elements $u$ in $\gg$ with $\vare(w)=\vare_i(w)$. Let $w_j$ (1$\le
j\le n_i$) be in $ N$ such that

\no(1)$\vare_{i1}(w_j)=\vare_{i2}(w_j)=\cdots=\vare_{ij}(w_j)=1,$

\no(2) other
components of $\vare(w_j)$ are 0.

\no Then $w_j$ is corresponding to the $(i,j)$-th fundamental weight of
$F_\l$. Using the bijection $\vare:
\gg\to$Dom$(F_\l)$ we shall identify $N$ with the set $\bzi$ (see 4.1
for the definition of $\bzi$). Under the identification, we have
$$w=(\vare_{i1}(w),\vare_{i2}(w),...,\vare_{in_i}(w))$$ if $w\in N$.

\bigskip

\noindent{\bf Lemma 8.1.1.} {\sl Let $w\in N$.
Then we have $$t_{w_{n_i}}t_w=t_{w_{n_i}*w}.$$ See section 6.4 for the
definition of $w_{n_i}*w$.}

\bigskip

\itp When all components of $\vare(w)$ are nonnegative,
it is entirely similar to the proof for Theorem 6.4.1. Here we need
consider $t_{w_{n_i}}t_{w_1}^m$ first. Thus the lemma is true if all
components of $\vare(w)$ are nonnegative.

Note that $\om^{kn}w$ is in $\gg$. In general we can find a
positive integer $k$ such that all components of
$\vare(\om^{kn}w)$ are nonnegative. Let $y_j\ (1\le j\le p)$ be in
$\gg$ such that $\vare(y_j)=\vare_j(\om^{kn}w)$. Note that all
components of $\vare(y_i)$ are nonnegative and $y_i\in N$.
According to Theorem 6.4.1 and its proof, and Prop. 6.3.7, we have
$$\begin{array} {rl}
t_{w_{n_i}}t_w&=t_{w_{n_i}}t_{\om^{kn}w}t_{\om^{-kn}w_\l}\\
&=t_{y_1}\cdots t_{y_{i-1}}t_{w_{n_i}}t_{y_i}t_{y_{i+1}}\cdots
t_{y_p}t_{\om^{-kn}w_\l}\\ &=t_{y_1}\cdots
t_{y_{i-1}}t_{w_{n_i}*y_i}t_{y_{i+1}}\cdots
t_{y_p}t_{\om^{-kn}w_\l}\\
&=t_{y_1*\cdots*y_{i-1}*(w_{n_i}*y_i)*y_{i+1} \cdots
*{y_p}}t_{\om^{-kn}w_\l}\\ &= t_{w_{n_i}*w}.
\end{array}$$

The lemma is proved.
\bigskip

\noindent{\bf Lemma 8.1.2.} {\sl Let $w\in N$. Recall that
$t_{w_j}t_w=\sum_{u\in\gg}\ga_{w_j,w,u}t_u$. We have $u\in N$ and
$\vare_{ik}(w)\le \vare_{ik}(u)\le\vare_{ik}(w)+1$ for all $k$ if
$\ga_{w_j,w,u}\ne 0$.}

\bigskip

\itp First we show that $\vare_{ik}(w)\le \vare_{ik}(u)$. Using
Lemma 8.1.1 we see that it is harmless to assume that all
components of $\vare(w)$ are non-negative. Using Prop. 6.3.7 we
see that $u$ is in $N$ if $\ga_{w_j,w,u}\ne 0$. Since $t_{w_j}$
appears in $t_{w_1}^j$ with non-zero coefficient, considering
$t_{w_1}^jt_w$, by the positivity (see 1.3 (f)) and Prop. 6.3.7 we
see that $\vare_{ik}(w)\le\vare_{ik}(u) $ if $\ga_{w_j,w,u}\ne 0$.

\medskip

Now we show that $\vare_{ik}(u)\le\vare_{ik}(w)+1$. Let
$w'_j=w_{n_i}^{-1}*w_j$. Using Prop. 6.3.8 we see that $t_{w'_j}$
appears in $t_{w_1^{-1}}^{n_i-j}$ with non-zero coefficient. Now
consider $t_{w_1^{-1}}^{n_i-j}t_w$. By Lemma 8.1.1 we may assume
that all components of $\vare(w)$ are non-positive. By the
positivity and Prop. 6.3.8 we see that
$
\vare_{ik}(u)\le\vare_{ik}(w)$ and $u$ is in $N$ if $\ga_{w'_j,w,u}\ne 0$. Multiplying
both sides of $t_{w'_j}t_w=\sum_{u\in N}\ga_{w'_j,w,u}t_u$ by
$t_{w_{n_i}}$ and using Lemma 8.1.1 we see that $
\vare_{ik}(u)\le\vare_{ik}(w)+1$ for all $k$ if $\ga_{w_j,w,u}\ne 0$.

The lemma is proved.

\bigskip

\noindent{\bf Proposition 8.1.3.} {\sl Let $w\in N$. Then
$t_{w_j}t_w=\sum_{u}t_u$, where $u$ runs through the set consisting of
all $x\in N$ such that
$\vare(x)=\vare(w)+\tau_{ik_1}+\tau_{ik_2}+\cdots+\tau_{ik_j}$ for
some sequence $1\le k_1<k_2<\cdots<k_j\le n_i$. See the proof of
Theorem 6.4.1 for the definition of $\tau_{ik}$.}

\bigskip

\itp Using Lemma 8.1.1 we may and  will assume that all components
of $\vare(w)$
 are non-negative.

\medksip

We use induction on the sum
$$E(w)=\vare_{i1}(w)+\cdots+\vare_{in_i}(w).$$ When $E(w)=0$, the lemma
is trivial. When $E(w)=1$, the lemma follows from Prop. 6.3.7. Now
suppose that the lemma is true when $E(w)\le a$. We need show that the
lemma then is true for $E(w)=a+1$.

\medskip

We define two bilinear forms, one for $J_N$, the subring of $J_{\gg}$
spanned by all $t_u\ (u\in N)$, the other for $R_{GL_{n_i}(\bbC)}$. We
define $$(\sum a_ut_u,\sum b_ut_u)=\sum a_ub_u$$ and
$$(\sum a_uV_u,\sum
b_uV_u)=\sum a_ub_u,$$ where $V_u\in R_{GL_{n_i}(\bbC)}$ stands for an
irreducible module of $GL_{n_i}(\bbC)$ of highest weight $u\in N$,
recall that we identify $N$ with $\bzi$.

\medskip

Let $u,y\in N$. We say that $t_u\sum_{x\in N}a_xt_x$ ($a_x\in\bbZ$)
has the right $y$-component if
$$(t_y,t_u\sum_{x\in N}a_xt_x)=(V_y,V_u\sum_{x\in N}a_xV_x).$$
We say that $t_u\sum_{x\in N}a_xt_x$ has the right form if it has the
right $y$-component for all $y\in N$.

By Prop. 6.3.7 and Theorem 7.2.2, we have

\bigskip

\no(a) if $t_u\sum_{x\in N}a_xt_x$ has the right $y$-component
(resp. right form), then both $(t_u\sum_{x\in N}a_xt_x )t_{w_1}$
and $(t_u\sum_{x\in N}a_xt_x )t_{w_2}$ have the right
$y$-component (resp. right form).

\bigskip

Obviously we have

\bigskip

\no(b) if $t_u\sum_{x\in N}a_xt_x$ and $t_u\sum b_xt_x$ have the right
$y$-component (resp. right form), then $t_u\sum_{x\in N}(a_x\pm
b_x)t_x $ has the right $y$-component (resp. right form).

\bigskip

We need to show that $t_{w_j}t_w$ has the right from.

\medskip

Let $\delta_k\in\bbZ^{n_i}$ be such that its $k$th component is 1
and other components are 0. By induction hypothesis, we know that
$t_{w_j}t_{a\delta_1}$ has the right form, recall that we identify
 $N$ with $\bzi$. Thus

 $$t_{w_j}t_{a\delta_1}=t_{(a+1)\delta_1+\delta_2+\cdots+\delta_j}+
 t_{a\delta_1+\delta_2+\cdots+\delta_j+\delta_{j+1}}.\leqno{\hbox{(c)}}$$

\no(Convention: If $k>n_i$, then $t_{x+\delta_k}=0$ for any $x\in
\bzi$;
 and if $x$ is in $\bbZ^{n_i}$ but not in $\bbZ^{n_i}_{\text{dom}}$, then
$t_x=0$.)

\medskip

Multiplying $t_{w_1}$ to both sides of (c) and using Prop. 6.3.7,
we
get

\bigskip

$$\begin{array}{rl}
t_{w_j}t_{a\delta_1}t_{w_1}=&t_{w_j}(t_{(a+1)\delta_1}+
t_{a\delta_1+\delta_2})\\
=&t_{(a+2)\delta_1+\delta_2+\cdots+\delta_j}+
 t_{(a+1)\delta_1+2\delta_2+\cdots+\delta_j}
 +t_{(a+1)\delta_1+\delta_2+\cdots+\delta_{j+1}}\\
 &+t_{(a+1)\delta_1+\delta_2+\cdots+\delta_j+\delta_{j+1}}
 +t_{a\delta_1+2\delta_2+\cdots+\delta_j+\delta_{j+1}}\\
 &+
 t_{a\delta_1+\delta_2+\cdots+\delta_j+\delta_{j+1}+\delta_{j+2}}.
\end{array}\leqno{\hbox{(d)}}$$

\bigskip

By (d) and Lemma 8.1.2, we see

\bigskip

\no(e) $t_{(a+2)\delta_1+\delta_2+\cdots+\delta_j}$ appears in
$t_{w_j}t_{(a+1)\delta_1}$ with coefficient 1; and all the three
elements $
 t_{(a+1)\delta_1+2\delta_2+\cdots+\delta_j},
 $
 $t_{a\delta_1+2\delta_2+\cdots+\delta_j+\delta_{j+1}}$,
$ t_{a\delta_1+\delta_2+\cdots+\delta_j+\delta_{j+1}+\delta_{j+2}}$
appear in $t_{w_j}t_{a \delta_1+\delta_2}$ with coefficient 1.

\bigskip

We wish to show that $t_{w_j}t_{(a+1)\delta_1}$ has the right form. By
induction hypothesis, $t_{w_j}t_{(a-1)\delta_1+\delta_2}$ has the
right form, so $t_{w_j}t_{(a-1)\delta_1+\delta_2}t_{w_1}$ has the
right form. We may use Prop. 6.3.7 to expand the expression
$t_{w_j}t_{(a-1)\delta_1+\delta_2}t_{w_1}$. Using Lemma 8.1.2 we see
that $ t_{(a+1)\delta_1+\delta_2+\cdots+\delta_j+\delta_{j+1} }$
appears in $t_{w_j}t_{a \delta_1+\delta_2}$ with coefficient 1. Thus
$t_{w_j}t_{a\delta_1+\delta_2}$ has the right form. By (d),
$t_{w_j}t_{(a+1)\delta_1}$ has the right form.

\medskip

We shall use the lexicographical order on $\bbZ^{n_i}$ with
$(1,0,...,0)>(0,1,...,0)$. Now suppose that if all components of
$w'$ are non-negative and $E(w')=E(w)=a+1,\ w'>w$, then
$t_{w_j}t_{w'} $ has the right form. We need show that
$t_{w_j}t_w$ has the right form. Choose $k$ such that
$\vare_{ik}(w)>0$ but $\vare_{il}(w)=0$ for all $l>k$. If $k=1$,
we have showed  that $t_{w_j}t_w$ has the right form. So we may
assume that $k>1$.

\medskip

 By induction hypothesis and (a), $t_{w_j}t_{w-\delta_{k-1}-\delta_k}t_{w_2}$ has the
right form. By induction hypothesis again, we see

\medskip

\no(f) $t_{w_j}(t_{w
}+t_{w-\delta_{k-1}+\delta_{k+1}}+t_{w-\delta_k+\delta_{k+1}}+
t_{w-\delta_{k-1}-\delta_{k}+\delta_{k+1}+\delta_{k+2}})$ has the
right form.

\medskip

Since $t_{w_j}t_{w-\delta_{k-1}-\delta_k+\delta_{k+1}}t_{w_1}$ has the
right form, by induction hypothesis, we get

\medskip

\no(g)
$t_{w_j}(t_{w-\delta_{k-1}+\delta_{k+1}}+t_{w-\delta_k+\delta_{k+1}}+
t_{w-\delta_{k-1}-\delta_{k}+\delta_{k+1}+\delta_{k+2}})$ has the
right form.

\medskip

Using (b) and (f-g), we see that $t_{w_j}t_w$ has the right form if
all components of $\vare(w)$ are nonnegative.

\medskip

The lemma is proved.

\medskip

\section{The based ring $J_{\Gamma_\lambda\cap\Gamma_\lambda^{-1}}$ and the based ring $J_{\mathbf c}$}

Now we can prove Conjecture 2.3.3.

\bigskip

\no{\bf Theorem 8.2.1} {\sl The map $t_w\to V(\vare(w))$, $w\in\gg$,
 defines a ring isomorphism from $J_{\gg}$ to $R_{F_\l}$.}

\bigskip

\itp Use Theorem 6.4.1, Lemma 8.1.3 and 4.2 (f).

\bigskip

Combining Theorem 2.3.2, Theorem 5.2.6 and Theorem 8.2.1, we see
that Lusztig Conjecture on based ring for type $\tilde A_{n-1}$ is
true.

\bigskip

\no{\bf 8.2.2.} Here  we give a summary of our proof of Lusztig
Conjecture on based ring for type $\tilde A_{n-1}$. Let {\bf c} be
a two-sided  cell of the extended affine Weyl group $W$ associated
with $GL_n(\bbC)$. According to Shi [S] and Lusztig [L3], we have
a corresponding partition $\l$ of $n$. Let $\mu$ be the dual of
$\l$ and $u\in GL_n(\bbC)$ a unipotent element whose Jordan form
has partition $\mu$. Then the centralizer of $u$ in $GL_n(\bbC)$
is connected and so is a maximal reductive subgroup $F_\l$ of the
centralizer. Thus Lusztig Conjecture on the based ring $J_{\bold
c}$ says that $J_{\bold c}$ is isomorphic to the $n_\mu\times
n_\mu$ matrix algebra  over the representation ring $R_{F_\l}$ of
$F_\l$, where $n_\mu$ is the number of left cells of $W$ in $\bold
c$.

\medksip

To prove the required result we first show that $J_{\bold c}$ is
isomorphic to the $n_\mu\times n_\mu$ matrix algebra  over the
based ring of the intersection of a left cell in {\bf c} and its
inverse, see Theorem 2.3.2. Then we establish a bijection between
the intersection and the set Irr$F_\l$ of isomorphism classes of
rational irreducible modules of $F_\l$ in Chapter 5, see Theorem
5.2.6. In Chapters 6-8 we show that the bijection leads to an
isomorphism between the based ring of the intersection and
$R_{F_\l}$. This completes our proof of Lusztig Conjecture for
type $A_{n-1}$.

\bigskip

\no{\bf 8.2.3.} Let $u$ be a unipotent element in $GL_n(\bbC)$
whose Jordan blocks are determined by the dual partition of a
partition $\l$ of $n$. Let $\mathcal B_u$ be the variety
consisting of all Borel subgroups of $GL_n(\bbC)$ that contain u.
In [X4] we show that there is an equivariant $F_\l$-partition of
$\mathcal B_u$. Using the partition we determine the equivariant
$K$-group $K^{F_\l}(\mathcal B_u\times\mathcal B_u)$. According to
3.2 (b) in [X4] and Theorems 2.3.2, 5.2.6, 8.2.1, we know that as
$R_{F_\l}$-modules, the equivariant $K$-group $K^{F_\l}(\mathcal
B_u\times\mathcal B_u)$ is isomorphic to $J_{\bold c}$, where
$\bold c$ is the two-sided cell of $W$ corresponding to $\l$.

\medskip

One can define a convolution on $K^{F_\l}(\mathcal B_u\times\mathcal
B_u)$. In [L10] Lusztig found some canonical bases for some
equivariant $K$-groups. It is likely that

\no(1) under the convolution the equivariant $K$-group $K^{F_\l}(\mathcal
B_u\times\mathcal B_u)$ becomes an associative ring with 1,

\no(2) there exists a canonical $\bbZ$-basis of
 $K^{F_\l}(\mathcal B_u\times\mathcal B_u)$ whose elements are one to
 one corresponding to the elements of the two-sided cell $\bold c$,

\no(3) the bijection between the canonical $\bbZ$-basis in (2) and
$\bold c$ leads to the ring isomorphism between $K^{F_\l}(\mathcal
B_u\times\mathcal B_u)$ and $J_{\bold c}$.

\no(4) the map from an affine Hecke algebras to  the based
ring of a two-sided cell of the corresponding extended affine Weyl
group defined in [L5] has a natural geometric explanation.

\medskip

In next section we explain that Lusztig Conjecture on based ring can
not be generalized to arbitrary extended affine Weyl groups.

\bigskip

\section{$PGL_n(\bbC)$}

 In this  section  we consider the projective linear group $PGL_n(\bbC)$ of degree
n. The extended affine Weyl group associated with  $PGL_n(\bbC)$
is just an affine Weyl group of type $A_{n-1}$. We shall identify
it with the affine Weyl group $W'$ in 2.1.3 (c). For each
two-sided cell $\bc$ (resp. left cell $\Ga$, right cell $\Phi$) of
the extended affine Weyl group $W$ associated with  $GL_n(\bbC)$,
the intersection $\bc'=\bc\cap W'$ (resp. $\Ga\cap W'$, $\Phi\cap
W'$) is a two-sided cell (resp. left cell, right cell) of $W'$.
Obviously, the based ring $J_{\bc'}$ of $\bc'$ is a subring of
$J_{\bc}$.

Now assume that $n=2$ and $\bc'$ is the lowest two-sided cell of
$W'$, i.e., the two-sided cell containing $s_1, s_0$. The
two-sided cell $\bc'$ contains two left cells. Let $\Ga_i\
(i=0,1)$ be the left cell in $\bc'$ such that $R(\Ga_i)=\{s_i\}$.
Then $$\Ga_1=\{s_1,\ s_0s_1,\ s_1s_0s_1,\ s_0s_1s_0s_1,\ ... \},$$
$$\Ga_0=\{s_0,\ s_1s_0 ,\
 s_0s_1s_0,\  s_1s_0s_1s_0,\ ...
\}.$$
The reductive group corresponding to $\bc'$ is
$F_{\bc'}=PGL_2(\bbC)$.

Suppose that there existed a finite $F_{\bc'}$-set $Y$ and a
bijection $\pi:\bc'\to$ the set of isomorphism classes of
irreducible $F_{\bc'}$ vector bundles on $Y\times Y$ such that the
map $t_w\to \pi(w)$ defines a ring isomorphism (preserving the
unit element) between $J_{\bc'}$ and $K_{F_{\bc'}}(Y\times Y)$.
Since $F_{\bc'}$ is connected, $Y$ must be a trivial
$F_{\bc'}$-set. Moreover, since $\bc'$ contains two left cells,
this forces that $Y$ contains two elements. Thus
$K_{F_{\bc'}}(Y\times Y)$ is isomorphic to the $2\times 2$ matrix
algebra  $M_2(R_{F_{\bc'}})$ over the representation ring
$R_{F_{\bc'}}$ of $F_{\bc'}$. As a consequence there should exist
an element $w$ in $\Ga_0\cap\Ga_1^{-1}$ and an element $u$ in
$\Ga_1\cap\Ga^{-1}_0$ such that $t_wt_u=t_{s_1}$. But this is
impossible since we have $$\Ga_0\cap\Ga_1^{-1}=\{s_1s_0,\
s_1s_0s_1s_0,\ s_1s_0s_1 s_0s_1s_0 ,\ ... \},$$
$$\Ga_1\cap\Ga_0^{-1}=\{s_0s_1,\ s_0s_1s_0s_1 ,\
 s_0s_1s_0s_1s_0s_1 ,\ ...
\},$$ and $$t_{(s_1s_0)^k}t_{(s_0s_1)^l}=\sum_{|k-l|\le
q\le|k+l-1|}t_{s_1(s_0s_1)^q}.$$

Therefore we can not find the required $F_{\bc'}$-set $Y$ and map
$\pi$ for the two-sided cell $\bc'$ of $W'$, so that Lusztig
Conjecture on based ring can not be generalized to arbitrary
extended affine Weyl group. However a weak form of the conjecture
might be true which now we are going to state. Let $W_G$ be the
extended affine Weyl group associated with  a connected reductive
algebraic group $G$ over $\bbC$ and let $\bc$ be a two-sided cell
of $W_G$. Denote by $F_\bc$ the reductive group corresponding to
$\bc$. It is likely that for any left cell $\Ga$
 in
$\bc$ we can find a finite transitive $F_\bc$-set $Y_\Ga$ such that
there exists a bijection $\pi_\Ga:\Ga\to$ the set of isomorphism
classes of irreducible $F_{\bc}$ vector bundles on $Y_\Ga\times Y_\Ga$
such that the map $t_w\to \pi(w)$ defines a ring isomorphism
(preserving the unit element) between $J_{\Ga\cap\Ga^{-1}}$ and
$K_{F_{\bc'}}(Y_\Ga\times Y_\Ga)$.

\medskip

Now we show that the weak form of Lusztig Conjecture on based ring
is true for $W'$. Let $\bc$, $\Ga_i$, $\Ga_\l$, $A_{ij}$ be as in
section 2.3. Then $\Ga'_i=\Ga_i\cap W'$ is a left cell of $W'$.
Clearly we have

\bigskip

\no(a) The map $\phi_{ii}: A_{ii}\to A_{11}$ in section 2.3 induces a bijection
$\phi'_{ii}: \Ga'_i\cap{\Ga'}_i^{-1}\to \Ga'_1\cap{\Ga'}_1^{-1}$.

\bigskip

Using Theorem 2.3.2 (a-b) we get

\bigskip

\no(b) The map $t_w\to t_{\phi'_{ii}(w)}$ induces a ring isomorphism from
$J_ {\Ga'_i\cap\ {\Ga'}_i^{-1}}$ to $J_{\Ga'_\l\cap{\Ga'}_\l^{-1}}$,
where $\Ga'_\l=\Ga_1\cap W'=\Ga_\l\cap W'$ is the left cell of $W'$
containing $w_\l$.

\medskip

\no(c) The based ring $J_{\Ga'_\l\cap{\Ga'}_\l^{-1}}$ is commutative.

\bigskip

Recall that we also regard $W$ as a permutation group of $\bbZ$ (see
2.1). Then $W'$ is a permutation group of $\bbZ$ consisting of all the
permutations $\sigma$ that satisfy (1) $\sigma(i+n)=\sigma(i)+n $ for
all $i\in\bbZ$ and (2) $\sum_{i=1}^n(\sigma(i)-i)=0.$ Thus we have

\bigskip

\no(d) Let $w$ be in $\Ga_\l\cap\Ga_\l^{-1}$. Then $w$ is in
$\Ga'_\l\cap{\Ga'}_\l^{-1}$ if and only if the sum of all
$\vare_{ij}(w)$ ($1\le i\le p,\ 1\le j\le n_i)$ is 0.

\bigskip

Let $\bar F_\l$ be the quotient group of $F_\l$ moduloing the center
of $GL_n(\bbC)$. Then $\bar F_\l$ is isomorphic to the reductive group
corresponding to the two-sided cell of $W'$ containing $w_\l$. Let
Dom$(\bar F_\l)$ be the subset of Dom$(F_\l)$ consisting of all
elements $(\vare_{ij})$ in Dom$(F_\l)$ that satisfy $
\displaystyle\sum_
{\st{\sc 1\le i\le p } { 1\le j\le n_i} }
\vare_{ij}=0$. Then we have

\bigskip

\no(e) The map $\vare:\Ga'_\l\cap{\Ga'}_\l^{-1}\to$Dom$(\bar F_\l)$
is a well defined bijection.

\medskip

\no(f) The set of isomorphism classes of  rational
representations of $\bar F_\l$ is one to one corresponding to the set
Dom$(\bar F_\l)$.

\bigskip

For $\vare\in$Dom$(\bar F_\l)$ we shall use $V(\vare)$ to denote a
rational representation of $\bar F_\l$ with highest weight $\vare$.
According to the above discussion and using Theorem 8.2.1 we get

\bigskip

\no{\bf Theorem 8.3.1.} {\sl The map $t_w\to V(\vare(w))$,
 $w\in\Ga'_\l\cap{\Ga'}_\l^{-1} $,
 defines a ring isomorphism from $J_{\Ga'_\l\cap{\Ga'}_\l^{-1}}$
 to $R_{\bar F_\l}$.}

\bigskip

In next section we will show that Lusztig Conjecture on based ring
is true for the extended affine Weyl group associated with  the
special linear group $SL_n(\bbC)$.

\bigskip


\def\om{\omega}

\section{$SL_n(\bbC)$}

\noindent{\bf 8.4.1.} In this section we show that Lusztig
Conjecture on based ring is true for the extended affine Weyl
group $\bar W$ associated with  the special linear group
$SL_n(\bbC)$ of degree $n$. It is known that $\bar W$ is a
quotient group of the extended affine Weyl group $W$ associated
with $GL_n(\bbC)$. The kernel is generated by $\om^n$. For an
element $w$ or a subset $K$ in $W$ we shall use $\bar w$ or $\bar
K$ to denote its image in $\bar W$. The following property is
obvious.

\bigskip

\noindent(a) Let $w$ and $u$ be in $W$. Then $w$ and $u$ are
contained in the same left cell (resp. right cell, two-sided cell)
of $W$ if and only if $\bar w$ and $\bar u$ are contained in the
same left cell (resp. right cell, two-sided cell) of $\bar W$.

\bigskip

 Let $\bar H$ be the Hecke algebra of $\bar W$ over $\A=\bbZ[q,q^{-1}]$
with parameter $q^2$. Then $\bar H$ is a quotient algebra of the Hecke
algebra $H$ of $W$. The kernel is generated by all $T_{\om^nw}-T_w$.
Let $\bar w,\ \bar u,\ \bar v$ be in $\bar W$. As for $W$, we can
define the Laurant polynomial $h_{\bar w,\bar u,\bar v}$ by means of
the Kazhdan-Lusztig basis of $\bar H$ and define the integer
$\gamma_{\bar w,\bar u,\bar v}$ using the Laurant polynomial $h_{\bar
w,\bar u,\bar v}$.

\bigskip

Let $w,u,v$ be in $W$. Then we have

\medskip

\noindent(b) $h_{w,u,v}=h_{\bar w,\bar u,\bar v}$.

\medskip

\noindent(c) $\gamma_{w,u,v}=\gamma_{\bar w,\bar u,\bar v}$.

\bigskip

Let $\bar J$ be the based ring of $\bar W$ with basis elements
$t_{\bar w},\
\bar w\in \bar W$ and structure constants $\gamma_{\bar w,\bar u,\bar
v}$. Then $\bar J$ is a quotient ring of the based ring $J$ of $W$.
The kernel is generated by all $t_{\om^nw}-t_w$. Let {\bf c} be a
two-sided cell of $W$ and $\bar {\mathbf c}$ the two-sided cell of
$\bar W$ corresponding to $\mathbf c$. Then $\bar J_{\bar{\mathbf c}}$
is a quotient ring of the based ring $J_{\mathbf c}$ of $\mathbf c$,
the kernel is generated by all $t_{\om^nw}-t_w$, $w\in\mathbf c$.
Similar conclusion is true for the based ring $J_{\bar
\Gamma\cap\bar\Gamma^{-1}}$, where $\bar\Gamma$ is a left cell in
$\bar\mathbf c$.

\medskip

Let $\Gamma_\lambda,\ \Gamma_i,\ A_{ij}$ and $\phi_{ij}$ be as in
section 2.3. The image in $\bar W$ of $A_{ij}$ is denoted by $\bar
A_{ij}$. Then $\phi_{ij}$ induces a bijection from $\bar A_{ij}$ to
$\bar A_{11}$. We denote the induced map also by $\phi_{ij}$.
According to the above discussion and Theorem 2.3.2 we have

\bigskip
\def\bc{\mathbf c}
\def\no{\noindent}
\def\Ga{\Gamma}
\def\l{\lambda}

\noindent{\bf Theorem 8.4.2.} {\sl Let $\bar\bc$ be the two-sided
cell
of $\bar W$ corresponding to a partition $\l$ of $n$ and $\mu$ the
dual partition of $\l$ and $\bc$ the corresponding two-sided
cell in
$W$.}

\smallskip

\no(a) {\sl The map $t_{\bar w}\to t_{\phi_{ii}(\bar w)}$ induces
a ring isomorphism
 from $J_{\bar\Ga_i\cap\bar
 \Ga_i^{-1}}$ to  $J_{\bar\Ga_\l\cap\bar\Ga_\l^{-1}}$.}

\no(b) {\sl The based ring $J_{\bar\Ga_\l\cap\bar\Ga_\l^{-1}}$ is commutative.}

\no(c) {\sl The map $t_{\bar w}\to E(t_{\phi_{ij}(\bar w)},i,j), \
\bar w\in \bar A_{ij}$ defines an isomorphism from the based ring
$J_{\bar\bc}$ to $M_{n_\mu}(J_{\bar\Ga_\l^{-1}\cap\bar\Ga_\l})$,
the $n_\mu\times n_\mu$ matrix algebra over the ring
$J_{\bar\Ga_\l^{-1}\cap\bar\Ga_\l}$.}

\bigskip

\noindent{\bf 8.4.3.} Let $u$ and $F_\l$ be as in Chapter 4. Then $u$ is in $SL_n(\bbC)$.
Let $F'_\l=F_\l\cap SL_n(\bbC)$. Then $F'_\l$ is a maximal reductive
subgroup of the centralizer in $SL_{n}(\bbC)$ of $u$. Since $F'_\l$
contains the derived group of $F_\l$, we have

\bigskip

\no(a) The restriction to $F'_\l$ of an irreducible rational representation
of $F_\l$ is irreducible. Any irreducible representation of $F'_\l$
can be obtained in this way.

\bigskip

Let Dom$(F_\l$) be as in section 4.1 and $r_i$ be as in section 5.1.
The following result is obvious from the definition of $F_\l$ and
$F'_\l$.

\bigskip
\def\vare{\varepsilon}

\no(b) Let $\vare=(\vare_{11},...,\vare_{pn_p})$ be a dominant weight
in Dom$(F_\l)$ and $V(\vare)$ an irreducible representation of $F_\l$
of highest weight $\vare$. Then the restriction to $F'_\l$ of
$V(\vare)$ is trivial if and only if there exists an integer $k$ such
that $\vare_{ij}=kr_i$ for all $i=1,2,...,p$ and $j=1,...,n_i$.

\bigskip

We introduce an equivalence relation on Dom$(F_\l)$. Let
$\vare=(\vare_{ij})$ and $\xi=(\xi_{ij})$ be in Dom$(F_\l)$. We say
that $\vare$ and $\xi$ are equivalent if there exists an integer $k$
such that $\vare_{ij}-\xi_{ij}=kr_i$ for all $i,j$. The equivalence
class containing $\vare$ will be denoted by $\bar\vare$ and the set of
all equivalence classes in Dom$(F_\l)$ is denoted by Dom$(F'_\l)$. By
(b) we get

\bigskip

\no(c) For each element $\bar\vare$ in Dom$(F'_\l)$ we have an
irreducible rational representation $V(\bar\vare)$ of $F'_\l$ with
highest weight $\bar\vare$. The map $\bar\vare\to V(\bar\vare)$ is a
bijection from Dom$(F'_\l)$ to the set of isomorphism classes of
irreducible rational representations of $F'_\l$.

\bigskip

 Let $\Ga_\l$ be as in section 2.3 and let $w$ and $u$ be elements
 in $\Ga_\l\cap\Ga^{-1}_\l$. It is easy to see
 that $\bar w=\bar u$ if and only if $\vare(w)$ and $\vare(u)$ are
 equivalent. Thus for each $\bar w$ in $\bar\Ga_\l\cap\bar
 \Ga^{-1}_\l$ we
 have a well defined element $\vare (\bar w)=\overline{\vare(w)}$ in
 Dom$(F'_\l)$, where $w$ is a preimage in $W$ of $\bar w$. Using
 Theorem 5.2.6 (b) we get

\bigskip

\no{\bf Theorem 8.4.4.} {\sl The map
 $\vare: \bar w\to \vare(\bar
w)$ defines a bijection from $\bar\Ga_\l\cap\bar\Ga^{-1}_\l$ to
Dom$(F'_\l)$.}

\bigskip

Using 8.4.3 (a), 8.4.1 (c), Theorem 8.2.1 and Theorem 8.4.4 we get

\bigskip

\no{\bf Theorem 8.4.5.} {\sl The map $t_{\bar w}\to V(\vare(\bar
w))$ defines a ring isomorphism between the based ring
$J_{\bar\Ga_\l\cap\bar\Ga^{-1}_\l}$ and the representation ring
$R_{F'_\l}$ of $F'_\l$.}

\bigskip

Combining Theorems 8.4.2, 8.4.4 and 8.4.5 we see that Lusztig
Conjecture on based ring is true for the extended affine Weyl
group associated with  $SL_n(\bbC)$.

\medskip

It is expected that the explicit knowledge on the based rings will
have applications to understand the representations of Hecke algebras
of $W$. Also we can compute some $\mu(y,w)$ (the coefficient of
$q^{\frac12(l(w)-l(y)-1)}$ in the Kazhdan-Lusztig polynomial
$P_{y,w}$) using the explicit knowledge on the based rings, the
details will appear elsewhere.

\backmatter

\bibliographystyle{amsalpha}

\chapter*{Index}

\def\qq{\qquad}

$a$-function\qq1.2

affine Weyl group, extended affine Weyl group\qq1.3

based ring\qq1.5

 cell (left, right, two-sided)\qq1.2

complete d-chain set, complete d-antichain set\qq2.4

complete r-chain set, complete r-antichain set\qq2.4

complete d-antichain family, complete r-antichain family\qq2.4.5

d-chain, d-antichain\qq2.2

d-chain family, d-antichain family\qq2.2

d-chain family set, d-antichain family set \qq2.2

distinguished involution\qq1.3

equivalent antichain\qq2.4.3

Hecke algebra\qquad 1.1

Kazhdan-Lusztig basis\qq1.1

Kazhdan-Lusztig polynomial\qq1.1

level\qq6.3.3

 Lusztig Conjecture on based ring\qq1.5

r-chain, r-antichain\qq2.2

r-chain family, r-antichain family\qq2.2

r-chain family set, r-antichain family set\qq2.2

saturated d-antichain, saturated r-antichain\qq6.3

star operation (left, right)\qq1.4

$\l$-admissible\qq5.2

\vfill\eject

\chapter*{Notation}

\def\qq{\qquad}

\no\S1.1

$\mathcal A$

$C_w$

$l(w)$

$P_{y,w}$

 $T_w$

$y-w$

$y\prec w$

$\mu(y,w)$

\no\S1.2

$a(v)$

 $h_{w,u,v}$

$h'_{w,u,v}$

$L(w)$

$R(w)$

$\tilde T_w$

 $w\ll u$, $w\rl u, \ w\lrl u$

$w\el u$, $w\er u$, $w\elr u$

\no\S1.3

$\mathcal D$

$P$

$R$

$w_0$

 $W_0$

$X$

$\gamma_{w,u,v}$

 $\Omega$

\no\S1.4

$D_L(s,t)$, $D_R(s,t)$

${}^*w$, $w^*$

${}^*w^\star$

$<s,t>$

\no\S1.5

$F_{\bold c}$

Irr$F_{\bold c}$

$J$, $J_{\bold c}$, $\jg$

$R_{F_{\bold c}}$

$t_w$

\no\S2.1

 $s_i$

 $\tau_1,...,\tau_i,...,\tau_n$

$\om$

$\Om$

\no\S2.2

$n_\mu$

 $w_\l$

$\Gamma_\l$

 $\l(w)$

$\mu(w)$

\no\S2.3

$A_{ij}$

$F_\l$

 $M_{n_\mu}(J_{\gg})$

\no\S3.1

$x_1,...,x_i,...,x_n$

\no\S3.2

$x_I$

 $x_\nu$

$X^+$

 $\Gamma_{\bold c}$

$\Phi_{\bold c}$

\no\S3.3

$m_I$

 $m_x$

$N$

\no\S4.1

Dom$(F_\l$)

 $F_\l$

$n_k$

$p$

 $\bbZ_{\text{dom}}^n$

\no\S4.2

$V(x)$

\no\S5.1

$a_{ij}$

 $e_i$

$n_i$

$p$

 $r_i$

$\L_j$

\no\S5.2

 $Z_w$

 $\vare_{k,i,j}$, $\vare_{k,i,j}(w)$

$\vare(w)$

 $\vare_{ij}(w)$

$\vare_{ij}(Z)$

\no\S6.1

$f_{u,v,w}$

$W_\l$

$\tilde T(w)$

\no\S6.2

$l_c$

$m_c$

$\vare_i(w)$

\no\S6.4

$x*y$

$\tau_{il}$

\no\S8.1

$w_1,...,w_i,...,w_{n_i}$

\end{document}